\documentclass[1p]{elsarticle}

\usepackage{lineno}
\usepackage[colorlinks,linkcolor={blue}]{hyperref}
\modulolinenumbers[5]

\usepackage{amsfonts}
\usepackage{graphicx}
\usepackage{amsmath}
\usepackage{graphicx}
\usepackage{amssymb}

\newtheorem{theorem}{Theorem}

\newtheorem{corollary}[theorem]{Corollary}

\newtheorem{definition}[theorem]{Definition}

\newtheorem{Example}[theorem]{Example}

\newtheorem{lemma}[theorem]{Lemma}

\newtheorem{proposition}[theorem]{Proposition}

\newtheorem{remark}[theorem]{Remark}
\newenvironment{proof}[1][Proof]{\textbf{#1.} }{\ \rule{0.5em}{0.5em}}

\journal{Journal of Differential Equations}

\begin{document}

\begin{frontmatter}

\title{Geometric proof for normally hyperbolic invariant manifolds}

\author{Maciej J. Capi\'nski\footnote{Research supported by the Polish National Science Center Grant 2012/05/B/ST1/00355}}
\ead{mcapinsk@agh.edu.pl}
\address{AGH University of Science and Technology, al. Mickiewicza 10, 30-059 Krak\'ow, Poland}

\author{Piotr Zgliczy\'nski\footnote{Research supported by the Polish National Science Center Grant 2011/03/B/ST1/04780}}
\ead{umzglicz@cyf-kr.edu.pl}
\address{Jagiellonian University, ul. prof. Stanis\l awa \L ojasiewicza 6, 
30-348 Krak\'ow, Poland}




\begin{abstract}
We present a new proof of the existence of normally hyperbolic manifolds and their whiskers for maps. Our result is not perturbative. Based on the bounds on the map and its derivative, we establish the existence of the manifold within a given neighbourhood. Our proof follows from a graph transform type method and is performed in the state space of the system. We do not require the map to be invertible. From our method follows also the smoothness of the established manifolds, which depends on the smoothness of the map, as well as rate conditions, which follow from bounds on the derivative of the map. Our method is tailor made for rigorous, interval arithmetic based, computer assisted validation of the needed assumptions.
\end{abstract}


\end{frontmatter}

%\tableofcontents
%\linenumbers

%TCIDATA{Version=5.00.0.2606}
%TCIDATA{LaTeXparent=0,0,MMFedit.tex}

\section{Introduction}

The goal of our paper is to present a geometric proof of the existence of
normally hyperbolic invariant manifolds (NHIMs) for maps, in a vicinity of
an approximate invariant manifold. There are four important features of our
approach: 1) we do not assume that the given map is a perturbation of some
other map for which we have a normally hyperbolic invariant manifold, 2) we
do not require that the map is invertible, 3) the assumptions can be
rigorously checked with computer assistance if our approximation of the
invariant manifold is good enough 4) our method does not require high order
smoothness. From our proof follows the high order smoothness of the
manifolds (provided that the map is suitably smooth), but it is enough to
consider $C^{1}$ bounds for the proof of their existence.

In the standard approach to the proof of various invariant manifold
theorems, all considerations are done in suitable function spaces or
sequences spaces. Moreover the existence of the invariant manifold for
nearby map (or ODE) is usually assumed, see for example \cite{Ch, HPS, Wi}
and the references given there. Typically these proofs do not give any
computable bounds for the size of perturbation for which the invariant
manifold exists.

Our result is in similar spirit to a number of results for establishing
invariant manifolds that have recently emerged, which assume that there
exists a manifold that is `approximately' invariant, and provide conditions
that ensure the existence of a true invariant manifold within a given
neighborhood. In \cite{BLZ5} Bates, Lu and Zeng present such approach within
a context of semiflows, which makes their method general and applicable to
infinite dimensional systems and PDEs. Compared to \cite{BLZ5} our results
is more explicit. Contrary to \cite{BLZ5}, where main theorems about NHIM
require that some constants are sufficiently small depending on other
constants, in our main theorem we just have several explicit inequalities
between pairs of constants. In \cite{CCL, FH,FLS,HL1,HL2} Calleja, Celletti,
Haro, de la Llave, Figueras, Fontich and Sire provide a framework and
results of establishing existence of whiskered tori with quasi periodic
dynamics, which is suitable for computer assisted validation. Our approach
however allows for more general dynamics. All above proofs are based on
constructions in suitable function spaces.

In contrast to the above mentioned approach, in our method the whole proof
is made in the phase space. This method is not entirely new. For example, a
similar approach is adapted in the proof of Jones \cite{J} in the context of
slow-fast system of ODEs. Jones though considered a perturbation of a
normally hyperbolic invariant manifold. In \cite{Ca,CZ} an approach in the
same spirit as in this paper has been applied to establish existence of
topologically normally hyperbolic invariant manifolds. These results are
based on topological arguments and do not establish the smoothness and the
foliations of the invariant manifolds. Similar approach has been applied by
Berger and Bounemoura \cite{BB}, where persistence and smoothness of
invariant manifolds is established using geometric and topological methods.
The result relies though on a perturbation of a normally hyperbolic
invariant manifold.

The method in this paper is based on two types of conditions. The first are
the topological conditions, which we refer to as `covering relations'. These
ensure that we have good topological alignment of the coordinates of the set
within which we establish the existence of the manifold. The second type of
conditions are based on the first derivative of the map and we refer to
these as the `rate conditions'. Our rate conditions are in the same spirit
to those of Fenichel \cite{Fen1,Fen2}. They measure the strength of the
hyperbolic contraction and expansion (within a neighborhood in which we
search for our manifold), in comparison to the dynamics on the normal
coordinates. The stronger the hyperbolicity is, the higher is the order of
the smoothness that can be established.

Our construction of the manifolds follows from a graph transform type
method. We prove that the manifolds emerge from passing to the limit of
graphs in appropriate coordinates. This construction follows primarily from
the covering conditions. To prove that the manifolds are Lipschitz, we show
that our graphs are contained in cones (this is the approach that was taken
in \cite{CZ}). The novelty of this paper lies in the proof of the higher
order smoothness. In our proof, this follows from establishing appropriate
cone conditions for the graphs. We define higher order cones, which span
around Taylor expansions of the graphs. We prove that these cones are
preserved as we iterate the graphs. (Verification of this fact follows from
our rate conditions.) We then show that higher order cone conditions imply
higher order smoothness of the graphs, and that this smoothness is preserved
as we pass to the limit.

We emphasize that in order to apply our method it is sufficient have a good
guess on the position of the manifold and good estimates of the first
derivative of the map. We do not require any estimates on its higher order
derivatives. It is sufficient that we know that the map is appropriately
smooth, and that the first derivative implies our rate conditions.

We believe that this approach is very well suited for computer assisted
(rigorous, interval arithmetic based) validation of the needed assumptions.
Similar approach has already been successfully applied in \cite{Ca,CZ} in
the setting of the rotating H\'{e}non map, in \cite{CR} to establish the
center manifold in the restricted three body problem, in \cite{CS} in the
setting of a driven logistic map or in \cite{Wil} to establish a hyperbolic
attractor in the Kuznetsov system. All these results follow from
verification of cone conditions based on the estimates of the derivative. We
believe that such estimates also imply rate conditions, hence the method
from this paper can easily be used to establish smoothness and fibration of
the manifolds. At present moment it appears that of approaches to NHIMs
mentioned earlier in the introduction only the one based on the
parameterization method \cite{CCL, FH,FLS,HL1,HL2} are ready for computer
assisted proof. This method is however restricted by the requirement of the
quasi-periodic dynamics on the invariant torus.

The paper is organized as follows. After preliminaries introducing basic
notations in Sections~\ref{sec:main-results} we state our main result for
the case of the torus. Sections~\ref{sec:jet-evol}--\ref{sec:stb-fibers}
contain the proof of our main result for the torus. In Section~\ref%
{sec:bundles} we show to how our construction can be carried over from the
torus to arbitrary compact manifold. We decided to work first with the torus
rather then a general manifold, because in that case we can have a global
coordinate chart and the main ideas are not mixed with the technicalities
connected with different charts. In Section \ref{sec:num} we apply our
method to the rotating H\'enon map.

%TCIDATA{Version=5.00.0.2606}
%TCIDATA{LaTeXparent=0,0,MMFedit.tex}

\section{Preliminaries}

\subsection{Notations}

For a point $p=(x,y)$ we shall use $\pi _{x}(p)=x$ to denote the projection
of $p$ onto the $x$ coordinate. We use a notation $B_{k}(p,R)$ for a ball of
radius $R$ in $\mathbb{R}^{k}$, centered at $p$. To simplify notations we
shall write $B_{k}(R)=B_{k}(0,R)$. For a set $A\subset \mathbb{R}^{k}$ we
shall write $\overline{A}$ for closure of $A$ and $\partial A$ for the
boundary of $A$. Throughout the work, the notation $\left\Vert \,\right\Vert 
$ will stand for the Euclidean norm, unless explicitly stated otherwise. For
a set $U\subset \mathbb{R}^{n}$ and a continuous function (homotopy) $h:%
\left[ 0,1\right] \times U\rightarrow \mathbb{R}^{n}$, for $\alpha \in \left[
0,1\right] $ we shall write $h_{\alpha }(x)$ for $h(\alpha ,x)$.

\begin{definition}
Let $f:\mathbb{R}^{n}\rightarrow \mathbb{R}^{k}$ be a $C^{1}$ function. We
define the interval enclosure of the derivative $Df$ on $U\subset \mathbb{R}%
^{n}$, as a set $[Df(U)]\subset \mathbb{R}^{k\times n}$, defined as%
\begin{equation*}
\lbrack Df(U)]=\left\{ A=\left( a_{ij}\right) _{\substack{ i=1,\ldots ,k  \\ %
j=1,\ldots ,n}}:a_{ij}\in \left[ \inf_{x\in U}\frac{\partial f_{j}}{\partial
x_{j}}(x),\sup_{x\in U}\frac{\partial f_{j}}{\partial x_{j}}(x)\right]
\right\} .
\end{equation*}
\end{definition}

\begin{definition}
Let $A:\mathbb{R}^{n}\rightarrow \mathbb{R}^{n}$ be a linear map. Let $\Vert
x\Vert $ be any norm on $\mathbb{R}^{n}$, then we define 
\begin{equation*}
m(A)=\max \left\{ L\in \mathbb{R}:\Vert Ax\Vert \geq L\Vert x\Vert \text{
for all }x\in \mathbb{R}^{n}\right\} .
\end{equation*}%
For an interval matrix $\mathbf{A}\subset \mathbb{R}^{k\times n}$ we set 
\begin{equation*}
m(\mathbf{A})=\inf_{A\in \mathbf{A}}m(A).
\end{equation*}
\end{definition}

\subsection{Taylor formula}

In this section we quickly set up the notations for the Taylor formula. Let%
\begin{equation*}
f:\mathbb{R}^{n}\rightarrow \mathbb{R}^{m}
\end{equation*}%
The $k$-th derivative of $f$ at $p$ is a symmetric $k$-linear operator. On
the basis it is defined as%
\begin{equation*}
D^{k}f(p)\left( e_{i_{1}},\ldots ,e_{i_{k}}\right) =\left( \frac{\partial
^{k}f_{1}(p)}{\partial x_{i_{1}}\ldots \partial x_{i_{k}}},\ldots ,\frac{%
\partial ^{k}f_{m}(p)}{\partial x_{i_{1}}\ldots \partial x_{i_{k}}}\right) .
\end{equation*}%
Using the following multi-index notation $j=\left( j_{1},\ldots
,j_{n}\right) \in \mathbb{N}^{n},$ $h=\left( h_{1},\ldots ,h_{n}\right) \in 
\mathbb{R}^{n}$%
\begin{equation*}
\begin{array}{lll}
\left\vert j\right\vert =j_{1}+\ldots +j_{n},\medskip  & \qquad  & 
h^{j}=h_{1}^{j_{1}}h_{2}^{j_{2}}\ldots h_{n}^{j_{n}}, \\ 
j!=j_{1}!j_{1}!\ldots j_{n}!, & \qquad  & \frac{\partial ^{j}g}{\partial
x^{j}}=\frac{\partial ^{\left\vert j\right\vert }g}{\partial
x_{1}^{j_{1}}\ldots \partial x_{n}^{j_{n}}},%
\end{array}%
\end{equation*}%
we can write out the value of $D^{k}f(p)$ on the diagonal as%
\begin{equation*}
D^{k}f(p)(h^{[k]})=D^{k}f(p)(\underset{k}{\underbrace{h,\ldots ,h}})=\left( 
\begin{array}{c}
\sum_{\left\vert j\right\vert =k}\frac{k!}{j!}h^{j}\frac{\partial ^{j}f_{1}}{%
\partial x^{j}}(p) \\ 
\vdots  \\ 
\sum_{\left\vert j\right\vert =k}\frac{k!}{j!}h^{j}\frac{\partial ^{j}f_{m}}{%
\partial x^{j}}(p)%
\end{array}%
\right) .
\end{equation*}%
The above formula is convenient to formulate the multi-dimensional version
of the Taylor formula:%
\begin{equation*}
f\left( p+h\right) =f(p)+T_{f,m,p}(h)+R_{f,m,p}(h),
\end{equation*}%
where $T_{f,m,p}$ stands for the Taylor expansion of order $m$%
\begin{equation*}
T_{f,m,p}(h)=\sum_{k=1}^{m}\frac{D^{k}f(p)}{k!}(h^{[k]}),
\end{equation*}%
and the reminder $R_{f,m,p}(h)$ can be computed in the integral form 
\begin{equation*}
R_{f,m,p}(h)=\int_{0}^{1}\frac{\left( 1-t\right) ^{m}}{m!}%
D^{m+1}f(p+th)\left( h^{[m+1]}\right) dt.
\end{equation*}

For $f:\mathbb{R}^{n}\supset \mathrm{dom}(f)\rightarrow \mathbb{R}^{m}$ and
a set $A\subset \mathrm{dom}(f)$ we define%
\begin{eqnarray*}
\left\Vert D^{k}f(p)\right\Vert &=&\sup \left\{ \left\Vert D^{k}f(p)\left(
h^{[k]}\right) \right\Vert :\left\Vert h\right\Vert =1\right\} , \\
\left\Vert f\right\Vert _{C^{m}} &=&\sup_{p\in \mathrm{dom}%
(f)}\max_{\left\vert k\right\vert \leq m}\left\Vert \frac{\partial ^{k}f}{%
\partial x^{k}}(p)\right\Vert , \\
\left\Vert f(A)\right\Vert _{C^{m}} &=&\left\Vert f|_{A}\right\Vert _{C^{m}}.
\end{eqnarray*}

%TCIDATA{Version=5.00.0.2606}
%TCIDATA{LaTeXparent=0,0,MMFedit.tex}

\section{Main results\label{sec:main-results}}

The goal of this section is to set up the structure for our NHIM, which will
be diffeomorphic with a manifold $\Lambda $. To make the setup as simple as
possible we will focus on the special case where $\Lambda $ is a torus. This
will simplify notations in many of the arguments, since we will not need to
work with various local charts. We shall prepare the setup though in a way
that will allow for a straightforward generalization to an arbitrary
manifold without boundary. This will be done in section \ref{sec:bundles}.

\subsection{Definitions and setup}

In the simple situation when $\Lambda $ is an $c$-dimensional torus, we are
in a convenient situation, since we have a covering 
\begin{equation*}
\varphi :\mathbb{R}^{c}\rightarrow \Lambda =\left( \mathbb{R}/\mathbb{Z}%
\right) ^{c},
\end{equation*}%
which gives us the set of charts being the restriction of $\varphi $ to
balls $B$ in $\mathbb{R}^{c}$, which are small enough so that $\varphi
:B\rightarrow \Lambda $ is a homeomorphism on its image. We introduce a
notation $R_{\Lambda }>0$ for a radius such that $\varphi _{|B(\lambda
,R_{\Lambda })}$ is a homeomorphism onto its image. When $\Lambda $ is a
torus, we can simply take $R_{\lambda }=\frac{1}{2}.$ Introducing the
notation $R_{\Lambda }$ here though will simplify our future discussion in
section \ref{sec:bundles}, where we generalize the results.

Let $R<\frac{1}{2}R_{\Lambda }$ and denote by $D$ the set 
\begin{equation*}
D=\Lambda \times \overline{B}_{u}(R)\times \overline{B}_{s}(R),
\end{equation*}%
where $\overline{B}_{n}(R)$ stands for a closed ball of radius $R$, centered
at zero, in $\mathbb{R}^{n}$. We consider a $C^{k+1}$ map, for $k\geq 1$, 
\begin{equation*}
f:D\rightarrow \Lambda \times \mathbb{R}^{u}\times \mathbb{R}^{s}.
\end{equation*}

Throughout the paper we shall use the notation $z=(\lambda ,x,y)$ to denote
points in $D$. This means that notation $\lambda $ will stand for points on $%
\Lambda $, notation $x$ for points in $\mathbb{R}^{u}$, and $y$ for points
in $\mathbb{R}^{s}$. We will write $f$ as $(f_{\lambda },f_{x},f_{y})$ ,
where $f_{\lambda },f_{x},f_{y}$ stand for projections onto $\Lambda $, $%
\mathbb{R}^{u}$ and $\mathbb{R}^{s}$, respectively. On $\mathbb{R}^{c}\times 
\mathbb{R}^{u}\times \mathbb{R}^{s}$ we will use the Euclidian norm.

In view of the further generalization to arbitrary manifold let us stress
that our set $D$ can be thought as a subset of the trivial vector bundle $%
\mathbb{T}_{c}\times \mathbb{R}^{u}\times \mathbb{R}^{s}$.

\begin{definition}
The set of points which are in the same good chart with point $q\in D$ will
be denoted by 
\begin{equation*}
P(q)=\{z\in D\ |\ \Vert \pi _{\lambda }z-\pi _{\lambda }q\Vert \leq
R_{\Lambda }/2\}.
\end{equation*}
\end{definition}

Let $L\in \left( \frac{2R}{R_{\Lambda }},1\right) $, and let us define 
\begin{align*}
\mu _{s,1}& =\sup_{z\in D}\left\{ \left\Vert \frac{\partial f_{y}}{\partial y%
}\left( z\right) \right\Vert +\frac{1}{L}\left\Vert \frac{\partial f_{y}}{%
\partial (\lambda ,x)}(z)\right\Vert \right\} , \\
\mu _{s,2}& =\sup_{z\in D}\left\{ \left\Vert \frac{\partial f_{y}}{\partial y%
}\left( z\right) \right\Vert +L\left\Vert \frac{\partial f_{\left( \lambda
,x\right) }}{\partial y}(z)\right\Vert \right\} ,
\end{align*}%
\begin{align*}
\xi _{u,1}& =\inf_{z\in D}\left\{ m\left( \frac{\partial f_{x}}{\partial x}%
(z)\right) -\frac{1}{L}\left\Vert \frac{\partial f_{x}}{\partial \left(
\lambda ,y\right) }(z)\right\Vert \right\} , \\
\xi _{u,1,P}& =\inf_{z\in D}m\left[ \frac{\partial f_{x}}{\partial x}(P(z))%
\right] -\frac{1}{L}\sup_{z\in D}\left\Vert \frac{\partial f_{x}}{\partial
\left( \lambda ,y\right) }(z)\right\Vert , \\
\xi _{u,2}& =\inf_{z\in D}\left\{ m\left( \frac{\partial f_{x}}{\partial x}%
\left( z\right) \right) -L\left\Vert \frac{\partial f_{(\lambda ,y)}}{%
\partial x}(z)\right\Vert \right\} ,
\end{align*}%
\begin{align}
\mu _{cs,1}& =\sup_{z\in D}\left\{ \left\Vert \frac{\partial f_{\left(
\lambda ,y\right) }}{\partial \left( \lambda ,y\right) }(z)\right\Vert
+L\left\Vert \frac{\partial f_{\left( \lambda ,y\right) }}{\partial x}%
(z)\right\Vert \right\} ,  \notag \\
\mu _{cs,2}& =\sup_{z\in D}\left\{ \left\Vert \frac{\partial f_{\left(
\lambda ,y\right) }}{\partial \left( \lambda ,y\right) }(z)\right\Vert +%
\frac{1}{L}\left\Vert \frac{\partial f_{x}}{\partial \left( \lambda
,y\right) }(z)\right\Vert \right\} ,  \notag
\end{align}%
\begin{align*}
\xi _{cu,1}& =\inf_{z\in D}\left\{ m\left( \frac{\partial f_{(\lambda ,x)}}{%
\partial (\lambda ,x)}(z)\right) -L\left\Vert \frac{\partial f_{(\lambda ,x)}%
}{\partial y}(z)\right\Vert \right\} , \\
\xi _{cu,1,P}& =\inf_{z\in D}m\left[ \frac{\partial f_{(\lambda ,x)}}{%
\partial (\lambda ,x)}(P(z))\right] -L\sup_{z\in D}\left\Vert \frac{\partial
f_{(\lambda ,x)}}{\partial y}(z)\right\Vert , \\
\xi _{cu,2}& =\inf_{z\in D}\left\{ m\left( \frac{\partial f_{(\lambda ,x)}}{%
\partial \left( \lambda ,x\right) }(z)\right) -\frac{1}{L}\left\Vert \frac{%
\partial f_{y}}{\partial (\lambda ,x)}(z)\right\Vert \right\} .
\end{align*}

\begin{remark}
Throughout the work, the $L\in (\frac{2R}{R_{\Lambda }},1)$ is a fixed
constant. We shall later see that $L$ is associated with Lipschitz bounds on
the established manifolds (hence the choice of notation).
\end{remark}

The key to the naming of the constants is the following:

\begin{itemize}
\item $\xi _{u,\cdot }$, $\xi _{cu,\cdot }$ - the constants describing lower
bound on the expansion in the unstable or center-unstable directions.

\item $\mu _{s,\cdot }$, $\mu _{cs,\cdot }$ - the constants describing upper
bound for contraction constant in the stable or center-stable direction.

\item The number $1$ or $2$ as second lower index is used according to the
following rule: $1$, when both partial derivatives are of the same component
of $f$, for example $f_{(\lambda ,x)}$ in $\mu _{cs,1}$, while $2$ is used
the differentiation is done with respect to the same block of variables of
various components of $f$.

\item $\xi _{u,1}$, $\xi _{u,2}$, $\xi _{cu,1}$, $\xi _{cu,2}$ are the
expansion bounds and $\mu _{s.1}$, $\mu _{s,2}$, $\mu _{cs,1}$, $\mu _{cs,2}$
are the contraction bounds, that are used for the establishing of smoothness
of invariant manifolds and their fibres.

\item $\xi _{u,1,P}$, $\xi _{cu,1,P}$ are more stringent bounds (i.e. $\xi
_{\cdot ,1,P}\leq \xi _{\cdot ,1}$). They are used to ensure lower bounds on
the expansion on the $x$ and $\left( \lambda ,x\right) $ coordinates.
\end{itemize}

\begin{definition}
\label{def:rate-conditions}We say that $f$ satisfies rate conditions of
order $k\geq 1$ if $\xi _{u,1},$ $\xi _{u,1,P},$ $\xi _{u,2},$ $\xi _{cu,1},$
$\xi _{cu,1,P},$ $\xi _{cu,2}$ are strictly positive, and for all $k\geq
j\geq 1$ holds%
\begin{equation}
\mu _{s,1}<1<\xi _{u,1,P},  \label{eq:rate-cond-0}
\end{equation}%
\begin{align}
\frac{\mu _{cs,1}}{\xi _{u,1,P}}& <1,\qquad \frac{\mu _{s,1}}{\xi _{cu,1,P}}%
<1,  \label{eq:rate-cond-1} \\
\frac{\left( \mu _{cs,1}\right) ^{j+1}}{\xi _{u,2}}& <1,\qquad \frac{\mu
_{s,2}}{(\xi _{cu,1})^{j+1}}<1,  \label{eq:rate-cond-2} \\
\frac{\mu _{cs,2}}{\xi _{u,1}}& <1,\qquad \frac{\mu _{s,1}}{\xi _{cu,2}}<1.
\label{eq:rate-cond-3}
\end{align}%
We say that $f$ satisfies rate conditions of order zero, if only (\ref%
{eq:rate-cond-0})--(\ref{eq:rate-cond-1}) are satisfied.
\end{definition}

\begin{figure}[tbp]
\begin{center}
\includegraphics[height=4cm]{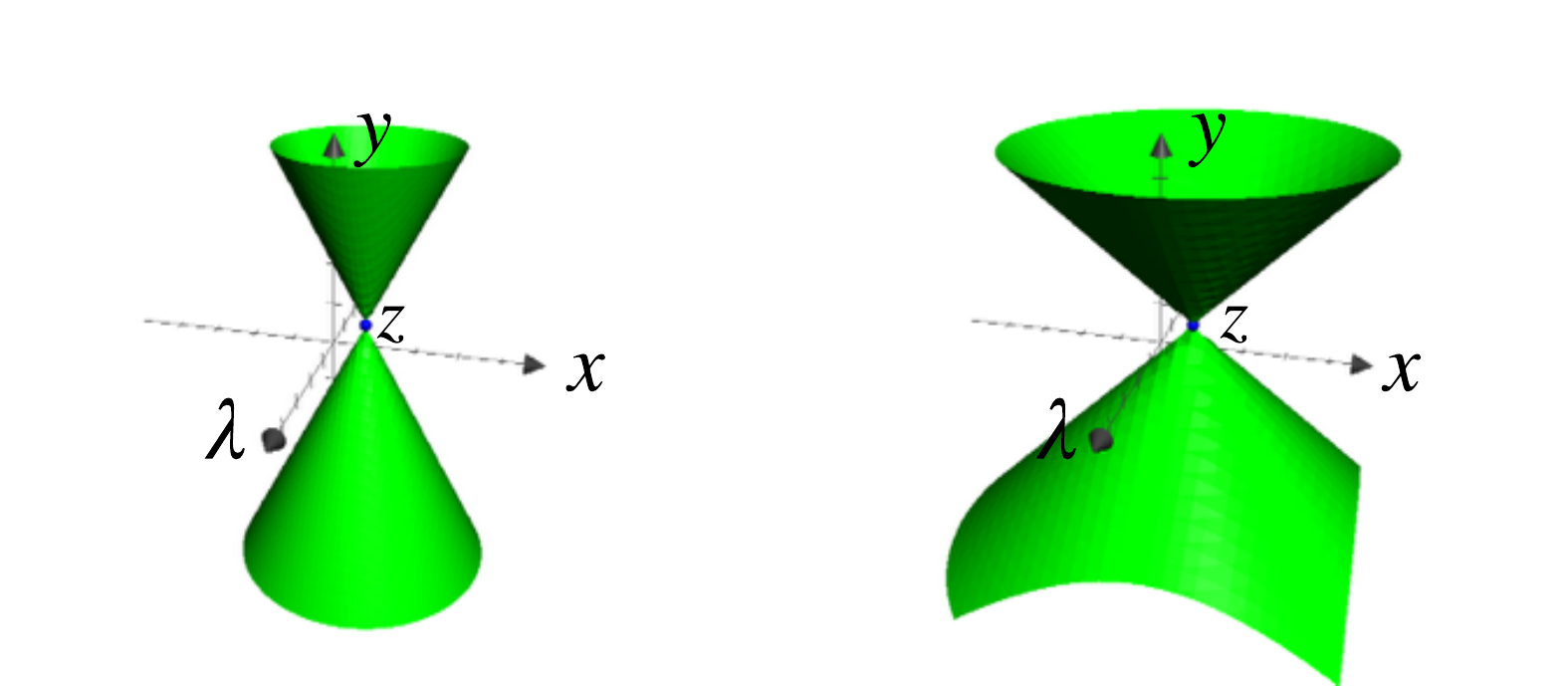}
\end{center}
\caption{The stable cone $J_{s}(z,M)$ for $M=\frac{1}{2}$ on the left, and $%
M=1$ on the right.}
\label{fig:st-cones}
\end{figure}

\begin{figure}[tbp]
\begin{center}
\includegraphics[height=4cm]{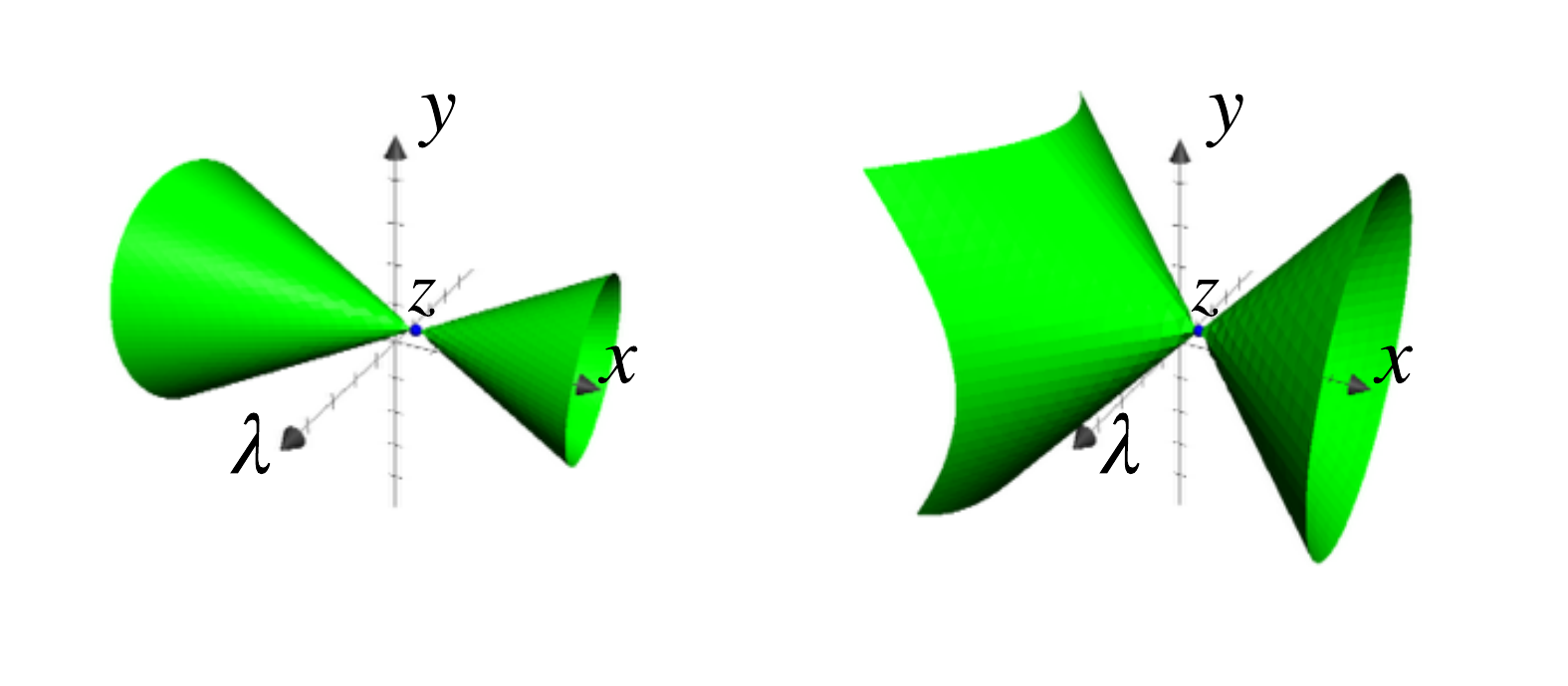}
\end{center}
\caption{The stable cone $J_{u}(z,M)$ for $M=\frac{1}{2}$ on the left, and $%
M=1$ on the right.}
\label{fig:unst-cones}
\end{figure}

We introduce the following notation: 
\begin{eqnarray}
J_{s}(z,M) &=&\left\{ \left( \lambda ,x,y\right) :\left\Vert \left( \lambda
,x\right) -\pi _{\lambda ,x}z\right\Vert \leq M\left\Vert y-\pi
_{y}z\right\Vert \right\} ,  \label{eq:stable-cone-3d} \\
J_{u}\left( z,M\right) &=&\left\{ \left( \lambda ,x,y\right) :\left\Vert
\left( \lambda ,y\right) -\pi _{\lambda ,y}z\right\Vert \leq M\left\Vert
x-\pi _{x}z\right\Vert \right\} .  \label{eq:unstable-cone-3d}
\end{eqnarray}%
We shall refer to $J_{s}(z,M)$ as a stable cone of slope $M$ at $z$, and to $%
J_{u}(z,M)$ as an unstable cone of slope $M$ at $z$. The cones are depicted
in Figures \ref{fig:st-cones}, \ref{fig:unst-cones}.

\begin{remark}
\label{rem:cone-in-chart}For any $z^{\ast }\in D$ and $z\in J_{u}(z^{\ast
},M)$ with $M\leq \frac{1}{L}$ we see that 
\begin{equation*}
\Vert \pi _{\lambda }(z^{\ast }-z)\Vert \leq \Vert \pi _{(\lambda
,y)}(z-z^{\ast })\Vert \leq 1/L\left\Vert \pi _{x}\left( z-z^{\ast }\right)
\right\Vert \leq 2R/L<R_{\Lambda }.
\end{equation*}%
This means that 
\begin{equation*}
J_{u}(z^{\ast },M)\cap D\subset \overline{B}_{c}(\lambda ^{\ast },R_{\Lambda
})\times \overline{B}_{u}(R)\times \overline{B}_{s}(R),
\end{equation*}%
for $\lambda =\pi _{\lambda }z^{\ast }$. Similarly, for $M\leq \frac{1}{L}$ 
\begin{equation*}
J_{s}(z^{\ast },M)\cap D\subset \overline{B}_{c}(\lambda ^{\ast },R_{\Lambda
})\times \overline{B}_{u}(R)\times \overline{B}_{s}(R).
\end{equation*}%
In other words, intersections of unstable (stable) cones with $D$ are
contained in sets on which we can use a single chart $P(z^{\ast })$.
\end{remark}

\begin{definition}
We say that a sequence $\left\{ z_{i}\right\} _{i=-\infty }^{0}$ is a (full)
backward trajectory of a point $z$ if $z_{0}=z,$ and $f\left( z_{i-1}\right)
=z_{i}$ for all $i\leq 0.$
\end{definition}

\begin{definition}
We define the center-stable set in $D$ as%
\begin{equation*}
W^{cs}=\left\{ z:f^{n}(z)\in D\text{ for all }n\in \mathbb{N}\right\} .
\end{equation*}
\end{definition}

\begin{definition}
We define the center-unstable set in $D$ as%
\begin{equation*}
W^{cu}=\{z:\text{there is a full backward trajectory of }z\text{ in }D\}.
\end{equation*}
\end{definition}

\begin{definition}
We define the maximal invariant set in $D$ as 
\begin{equation*}
\Lambda ^{\ast }=\{z:\text{there is a full trajectory of }z\text{ in }D\}.
\end{equation*}
\end{definition}

\begin{definition}
\label{def:Ws-fiber} Assume that $z\in W^{cs}$. We define the stable fiber
of $z$ as 
\begin{equation*}
W_{z}^{s}=\left\{ p\in D:f^{n}\left( p\right) \in J_{s}\left(
f^{n}(z),1/L\right) \cap D\text{ for all }n\in \mathbb{N}\right\} .
\end{equation*}
\end{definition}

\begin{definition}
\label{def:Wu-fiber} Assume that $z\in W^{cu}$. We define the unstable fiber
of $z$ as%
\begin{eqnarray*}
W_{z}^{u} &=&\{p\in D:\exists \text{ backward trajectory }\left\{
p_{i}\right\} _{i=-\infty }^{0}\text{ of }p\text{ in }D, \\
&&\quad \text{for any such backward trajectory} \\
&&\quad \text{and any backward trajectory }\left\{ z_{i}\right\} _{i=-\infty
}^{0}\text{ of }z\text{ in }D \\
&&\quad \text{holds }\left. p_{i}\in J_{u}\left( z_{i},1/L\right) \cap
D\right. \}.
\end{eqnarray*}
\end{definition}

The definitions of $W_{z}^{s}$ and $W_{z}^{u}$ are related to cones, which
is a nonstandard approach, the standard one is through convergence rates. We
will show that our definition implies the convergence rate as in the
standard theory.

Under our assumptions it will turn out that $f$ is injective on $W^{cu}$.
Therefore the backward orbit in the definition of $W^u_q$ is unique.

\begin{definition}
\label{def:back-cc}We say that $f$ satisfies backward cone conditions if the
following condition is fulfilled:

If $z_{1},z_{2},f(z_{1}),f(z_{2})\in D$ and $f(z_{1})\in J_{s}\left(
f(z_{2}),1/L\right) $ then%
\begin{equation*}
z_{1}\in J_{s}\left( z_{2},1/L\right) .
\end{equation*}
\end{definition}

\begin{remark}
The assumption that $f$ satisfies backward cone conditions will turn out to
be necessary in order to ensure that the established NHIM is a graph over $%
\Lambda $. After formulating our main Theorem \ref{th:main}, we follow up
with Examples \ref{exmpl:Mobius1}, \ref{exmpl:Mobius2}, in which we
demonstrate that without backward cone conditions the result cannot be
obtained.
\end{remark}

For $\lambda \in \Lambda $ we define the following sets: 
\begin{eqnarray*}
D_{\lambda } &=&\overline{B}_{c}\left( \lambda ,R_{\Lambda }\right) \times 
\overline{B}_{u}(R)\times \overline{B}_{s}(R), \\
D_{\lambda }^{+} &=&\overline{B}_{c}\left( \lambda ,R_{\Lambda }\right)
\times \overline{B}_{u}(R)\times \partial B_{s}(R), \\
D_{\lambda }^{-} &=&\overline{B}_{c}\left( \lambda ,R_{\Lambda }\right)
\times \partial \overline{B}_{u}(R)\times B_{s}(R).
\end{eqnarray*}

\begin{definition}
\label{def:covering}We say that $f$ satisfies covering conditions if for any 
$z\in D$ there exists a $\lambda ^{\ast }\in \Lambda $, such that the
following conditions hold:

For $U=J_{u}(z,1/L)\cap D$, there exists a homotopy $h$ 
\begin{equation*}
h:\left[ 0,1\right] \times U\rightarrow B_{c}\left( \lambda ^{\ast
},R_{\Lambda }\right) \times \mathbb{R}^{u}\times \mathbb{R}^{s},
\end{equation*}%
and a linear map $A:\mathbb{R}^{u}\rightarrow \mathbb{R}^{u}$ which satisfy:

\begin{enumerate}
\item $h_{0}=f|_{U},$

\item for any $\alpha \in \left[ 0,1\right] $, 
\begin{eqnarray}
h_{\alpha }\left( U\cap D_{\pi _{\theta }z}^{-}\right) \cap D_{\lambda
^{\ast }} &=&\emptyset ,  \label{eq:homotopy-exit} \\
h_{\alpha }\left( U\right) \cap D_{\lambda ^{\ast }}^{+} &=&\emptyset ,
\label{eq:homotopy-enter}
\end{eqnarray}

\item $h_{1}\left( \lambda ,x,y\right) =\left( \lambda ^{\ast },Ax,0\right) $%
,

\item $A\left( \partial B_{u}(R)\right) \subset \mathbb{R}^{u}\setminus 
\overline{B}_{u}(R).$
\end{enumerate}
\end{definition}

In the above definition a reasonable choice for $\lambda^{\ast}$ will be $%
\lambda^{\ast}=\pi_\lambda f(z)$. In fact any point sufficiently close to $%
\pi_\lambda f(z)$ will be also good.

\subsection{The main theorem}

\begin{theorem}
\label{th:main}(Main result) Let $k\geq 1$ and $f:D\rightarrow \Lambda
\times \mathbb{R}^{u}\times \mathbb{R}^{s}$ be a $C^{k+1}$ map. If $f$
satisfies covering conditions, rate conditions of order $k$ and backward
cone conditions, then $W^{cs},W^{cu}$ and $\Lambda ^{\ast }$ are $C^{k}$
manifolds, which are graphs of $C^{k}$ functions 
\begin{align*}
w^{cs}& :\Lambda \times \overline{B}_{s}(R)\rightarrow \overline{B}_{u}(R),
\\
w^{cu}& :\Lambda \times \overline{B}_{u}(R)\rightarrow \overline{B}_{s}(R),
\\
\chi & :\Lambda \rightarrow \overline{B}_{u}(R)\times \overline{B}_{s}(R),
\end{align*}%
meaning that 
\begin{align*}
W^{cs}& =\left\{ \left( \lambda ,w^{cs}(\lambda ,y),y\right) :\lambda \in
\Lambda ,y\in \overline{B}_{s}(R)\right\} , \\
W^{cu}& =\left\{ \left( \lambda ,x,w^{cu}(\lambda ,y)\right) :\lambda \in
\Lambda ,x\in \overline{B}_{u}(R)\right\} , \\
\Lambda ^{\ast }& =\left\{ \left( \lambda ,\chi (\lambda )\right) :\lambda
\in \Lambda \right\} .
\end{align*}%
Moreover, $f_{|W^{cu}}$ is an injection, $w^{cs}$ and $w^{cu}$ are Lipschitz
with constants $L$, and $\chi $ is Lipschitz with the constant $\frac{\sqrt{2%
}L}{\sqrt{1-L^{2}}}$. The manifolds $W^{cs}$ and $W^{cu}$ intersect
transversally, and $W^{cs}\cap W^{cu}=\Lambda ^{\ast }$.

The manifolds $W^{cs}$ and $W^{cu}$ are foliated by invariant fibers $%
W_{z}^{s}$ and $W_{z}^{u}$. The $W_{z}^{s}$ and $W_{z}^{u}$ are graphs of $%
C^{k}$ functions%
\begin{eqnarray*}
w_{z}^{s} &:&\overline{B}_{s}(R)\rightarrow \Lambda \times \overline{B}%
_{u}(R), \\
w_{z}^{u} &:&\overline{B}_{u}(R)\rightarrow \Lambda \times \overline{B}%
_{s}(R),
\end{eqnarray*}%
meaning that%
\begin{eqnarray*}
W_{z}^{s} &=&\left\{ \left( w_{z}^{s}\left( y\right) ,y\right) :y\in 
\overline{B}_{s}(R)\right\} , \\
W_{z}^{u} &=&\left\{ \left( \pi _{\lambda }w_{z}^{u}\left( x\right) ,x,\pi
_{y}w_{z}^{u}\left( x\right) \right) :x\in \overline{B}_{u}(R)\right\} .
\end{eqnarray*}%
The functions $w_{z}^{s}$ and $w_{z}^{u}$ are Lipschitz with constants $1/L$%
. Moreover, 
\begin{eqnarray*}
W_{z}^{s} &=&\{p\in D:f^{n}(p)\in D\text{ for all }n\geq 0,\text{ and} \\
&&\left. \exists n_{0},\exists C>0\text{ (which can depend on }p\text{)}%
\right. \\
&&\left. \text{s.t. for }n\geq n_{0},\text{ $f^{n}(p),f^{n}(z)$ are in the
same chart and }\right. \\
&&\left. \left\Vert f^{n}(p)-f^{n}(z)\right\Vert \leq C\mu
_{s,1}^{n}\right\} ,
\end{eqnarray*}%
and if $\{z_{i}\}_{i=-\infty }^{0}$ is the unique backward trajectory of $z$
in $D$, then 
\begin{eqnarray*}
W_{z}^{u} &=&\{p\in W^{cu}:\text{such that the unique backward trajectory }%
\{p_{i}\}_{i=-\infty }^{0} \\
&&\text{ of }p\text{ in }D\text{ satisfies the following condition} \\
&&\left. \exists n_{0}\geq 0,\exists \text{ }C>0\text{ (which can depend on }%
p\text{)}\right. \\
&&\left. \text{s.t. for }n\geq n_{0},\text{ $p_{-n},z_{-n}$ are in the same
chart and }\right. \\
&&\left. \left\Vert p_{-n}-z_{-n}\right\Vert \leq C\xi _{u,1,P}^{-n}\right\}
.
\end{eqnarray*}
\end{theorem}

Observe that bound on $L\in \left( \frac{2R}{R_{\Lambda }},1\right) $ gives
us lower bounds for the Lipschitz constants for functions $w^{cu}$, $w^{cs}$%
, $w^{u}$, $w^{s}$, which is clearly an overestimate for the case when $%
\mathbb{T}\times \{0\}\times \{0\}$ is our NHIM. This lower bound is a
consequence of choices we have made when formulating Theorem \ref{th:main},
as we did not want to introduce different constants for each type of cones,
plus several inequalities between them. However, below we give conditions
which allow to obtain better Lipschitz constants.

\begin{theorem}
\label{th:wsz}Let $M\in \left( 0,1/L\right) $ and 
\begin{eqnarray}
\mu &=&\sup_{z\in D}\left\{ \left\Vert \frac{\partial f_{y}}{\partial y}%
\left( z\right) \right\Vert +M\left\Vert \frac{\partial f_{y}}{\partial
(\lambda ,x)}(z)\right\Vert \right\} ,  \label{eq:wsz-mu} \\
\xi &=&\inf_{z\in D}m\left( \left[ \frac{\partial f_{(\lambda ,x)}}{\partial
\left( \lambda ,x\right) }(P(z))\right] \right) -\frac{1}{M}\sup_{z\in
D}\left\Vert \frac{\partial f_{(\lambda ,x)}}{\partial y}(z)\right\Vert .
\label{eq:wsz-xi}
\end{eqnarray}%
If assumptions of Theorem \ref{th:main} hold true and also $\frac{\xi }{\mu }%
>1$, then the function $w_{z}^{s}$ from Theorem \ref{th:main} is Lipschitz
with constant $M.$
\end{theorem}

\begin{theorem}
\label{th:wuz}Let $M\in(0,1/L)$ and 
\begin{eqnarray*}
\xi &=&\inf_{z\in D}m\left( \left[ \frac{\partial f_{x}}{\partial x}(P(z))%
\right] \right) -M\sup_{z\in D}\left\Vert \frac{\partial f_{x}}{\partial
\left( \lambda ,y\right) }(z)\right\Vert , \\
\mu &=&\sup_{z\in D}\left\{ \left\Vert \frac{\partial f_{\left( \lambda
,y\right) }}{\partial \left( \lambda ,y\right) }(z)\right\Vert +\frac{1}{M}%
\left\Vert \frac{\partial f_{\left( \lambda ,y\right) }}{\partial x}%
(z)\right\Vert \right\} .
\end{eqnarray*}%
If assumptions of Theorem \ref{th:main} hold true and also $\frac{\xi }{\mu }%
>1$, then the function $w_{z}^{u}$ from Theorem \ref{th:main} is Lipschitz
with constant $M.$
\end{theorem}

\begin{theorem}
\label{th:wcu}Let $M\in (0,L)$ and 
\begin{eqnarray*}
\xi &=&\inf_{z\in D}m\left[ \frac{\partial f_{(\lambda ,x)}}{\partial
(\lambda ,x)}(P(z))\right] -M\sup_{z\in D}\left\Vert \frac{\partial
f_{(\lambda ,x)}}{\partial y}(z)\right\Vert , \\
\mu &=&\sup_{z\in D}\left\{ \left\Vert \frac{\partial f_{y}}{\partial y}%
\left( z\right) \right\Vert +\frac{1}{M}\left\Vert \frac{\partial f_{y}}{%
\partial (\lambda ,x)}(z)\right\Vert \right\} .
\end{eqnarray*}%
If assumptions of Theorem \ref{th:main} hold true and also $\frac{\xi }{\mu }%
>1$, then the function $w^{cu}$ from Theorem \ref{th:main} is Lipschitz with
constant $M$.
\end{theorem}

\begin{theorem}
\label{th:wcs}Let $M\in (0,L)$ and 
\begin{eqnarray*}
\xi &=&\inf_{z\in D}m\left[ \frac{\partial f_{x}}{\partial x}(P(z))\right] -%
\frac{1}{M}\sup_{z\in D}\left\Vert \frac{\partial f_{x}}{\partial \left(
\lambda ,y\right) }(z)\right\Vert , \\
\mu &=&\sup_{z\in D}\left\{ \left\Vert \frac{\partial f_{\left( \lambda
,y\right) }}{\partial \left( \lambda ,y\right) }(z)\right\Vert +M\left\Vert 
\frac{\partial f_{\left( \lambda ,y\right) }}{\partial x}(z)\right\Vert
\right\}
\end{eqnarray*}%
If assumptions of Theorem \ref{th:main} hold true and also $\frac{\xi }{\mu }%
>1$, then the function $w^{cs}$ from Theorem \ref{th:main} is Lipschitz with
constant $M$.
\end{theorem}

\subsection{Comments on the inequalities and examples}

Let $J_{s}^{c}(z,M)$ and $J_{u}^{c}(z,M)$ stand for the complements of $%
J_{s}(z,M)$ and $J_{u}(z,M),$ respectively. We now comment about what
various inequalities in Definition \ref{def:rate-conditions} of rate
conditions mean and what they are needed for:

\begin{itemize}
\item $\mu _{cs,1}<\xi _{u,1,P}$: the forward invariance of $J_{u}(z,1/L)$
(Corollary \ref{cor:unstable-lip}). $\xi _{u,1,P}>1$: the expansion in $%
J_{u}(z,1/L)$ for $x$ - coordinate (Lemma \ref{lem:Ju-expansion}). This is
needed for the proof of the existence of $W^{cs}$ (Section \ref%
{sec:wcs-exists}).

\item $\xi _{cu,1,P}>\mu _{s,1}$: the forward invariance of $%
J_{s}^{c}(z,1/L) $ (Corollary \ref{cor:stable-lip}). $\mu _{s,1}<1$: the
contraction in $y$-direction in $J_{s}(z,1/L)$ (Lemma \ref%
{lem:Js-contraction}). This is needed for the proof of the existence of $%
W^{cu}$ (Section \ref{sec:wcu-exists}).

\item $\frac{\mu _{s,2}}{(\xi _{cu,1})^{j+1}}<1$, $j=1,\dots ,k$: the $C^{k}$%
-smoothness of $W^{cu}$ (Lemma \ref{lem:cu-smooth}).

\item $\frac{\left( \mu _{cs,1}\right) ^{j+1}}{\xi _{u,2}}<1$, $j=1,\dots ,k$
: the $C^{k}$-smoothness of $W^{cs}$ (Lemma \ref{lem:cs-smooth}).

\item $\frac{\mu _{cs,1}}{\xi _{u,1,P}}<1$: the existence of fibers $%
W_{q}^{u}$ (Lemma \ref{lem:wuq-conv}). $\frac{\mu _{cs,2}}{\xi _{u,1}}<1$:
the $C^{k}$ smoothness of $W_{q}^{u}$ (Lemma \ref{lem:Wuz-Ck}).

\item $\frac{\mu _{s,1}}{\xi _{cu,1,P}}<1$: the existence of fibers $%
W_{q}^{s}$ (Lemma \ref{lem:Wsz-lem1}). $\frac{\mu _{s,1}}{\xi _{cu,2}}<1$:
the $C^{k}$ smoothness of $W_{q}^{s}$ (Lemma \ref{lem:Wsz-lem2}).
\end{itemize}

\begin{figure}[tbp]
\begin{center}
\includegraphics[height=1.1in]{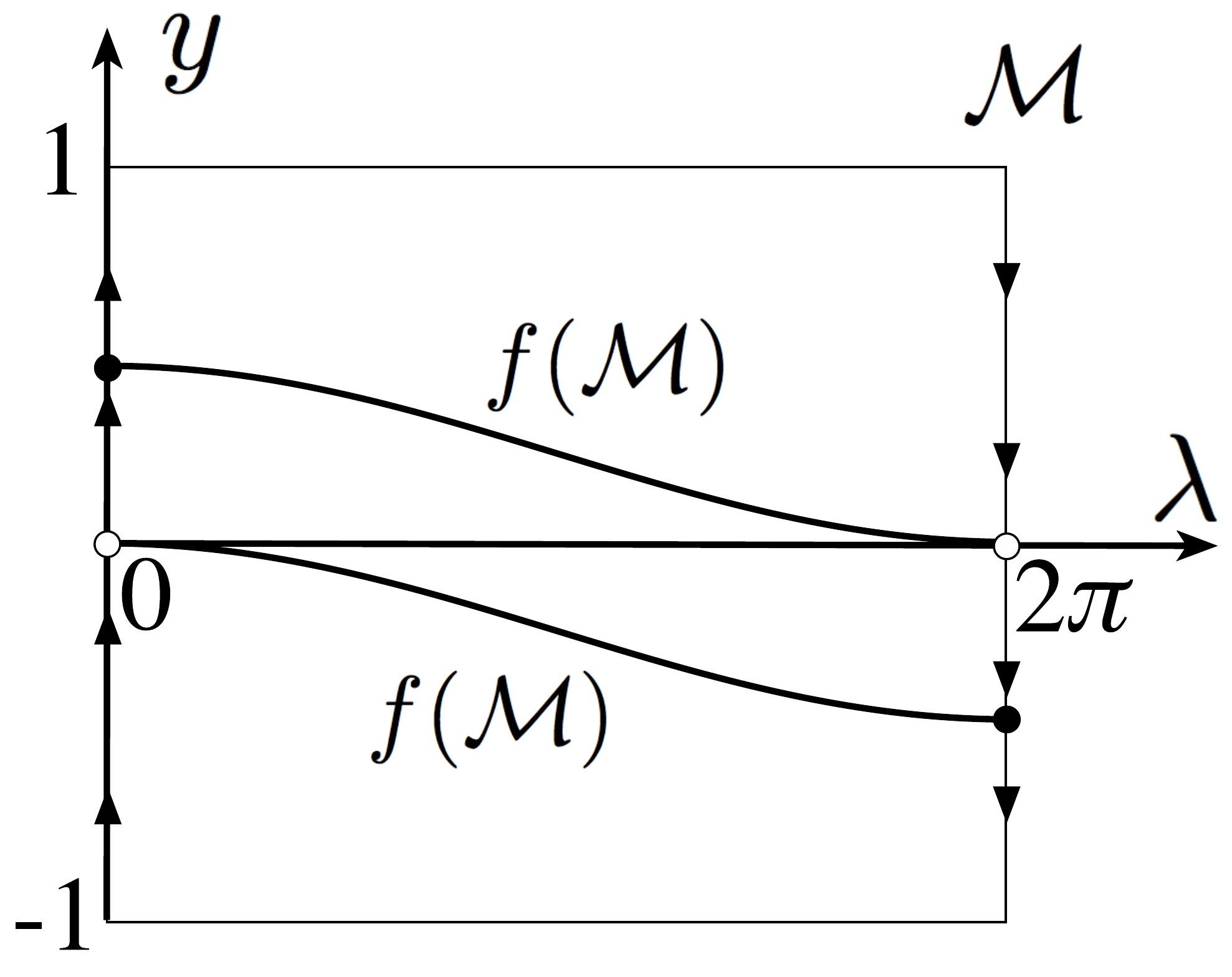}
\end{center}
\caption{The M\"{o}bius strip from Example \protect\ref{exmpl:Mobius1}.}
\label{fig:Mobius1}
\end{figure}

We now give two examples which show that in the absence of the backward cone
condition, the invariant set might not be a graph over $\Lambda $.

\begin{Example}
\label{exmpl:Mobius1}Consider a M\"{o}bius strip $\mathcal{M}$ depicted in
Figure \ref{fig:Mobius1}. The M\"{o}bius strip is parameterised by $\left(
\lambda ,y\right) $, with $\lambda \in \lbrack 0,2\pi )$ and $y\in \overline{%
B}_{s}(R)=\left[ -1,1\right] $. The two vertical edges which are glued
together are depicted with arrows.

Let $\xi >2$ be a constant. We consider a map $f:\mathcal{M}\times \overline{%
B}_{u}(R)\rightarrow \mathcal{M}\times \overline{B}_{u}(R),$%
\begin{equation*}
f(\left( \lambda ,y\right) ,x)=\left( \left( 2\lambda ,\frac{1}{4}+\frac{1}{4%
}\cos \lambda \right) ,\xi x\right) .
\end{equation*}%
On the unstable coordinate $x$, $f$ is simply a linear expansion. The stable
coordinate $y$ is the vertical coordinate on $\mathcal{M}$. The coordinates $%
\left( \lambda ,y\right) $ and $x$ are decoupled. Intuitively, on $\left(
\lambda ,y\right) $ the map does the following. It projects $\mathcal{M}$
into a horizontal circle, and then stretches it and wraps twice around $%
\mathcal{M}$ as in Figure \ref{fig:Mobius1}. For such a map all assumptions
of Theorem \ref{th:main} are fulfilled, except for the backward cone
conditions. We see that in the absence of the backward cone conditions, the
invariant manifold can be a set which is not a graph over $\Lambda $.
\end{Example}

\begin{figure}[tbp]
\begin{center}
\includegraphics[height=1.1in]{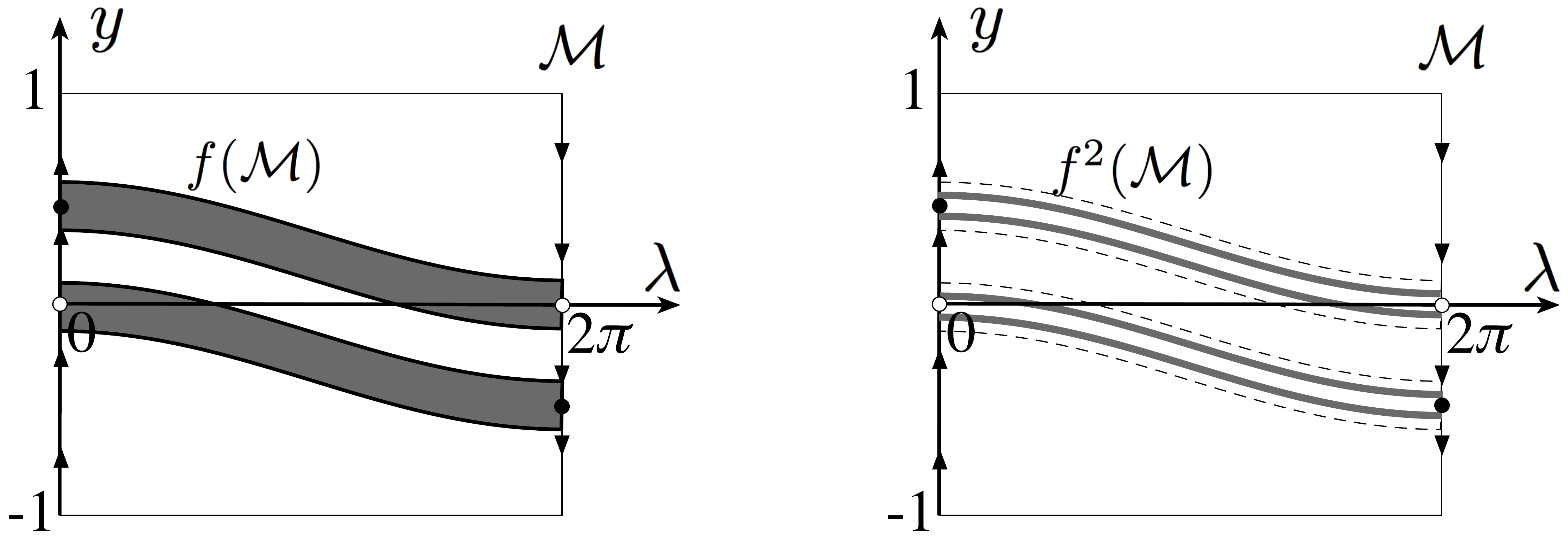}
\end{center}
\caption{The M\"{o}bius strip from Example \protect\ref{exmpl:Mobius2}.}
\label{fig:Mobius2}
\end{figure}

\begin{Example}
\label{exmpl:Mobius2}We can modify Example \ref{exmpl:Mobius1} slightly to
obtain a more interesting result. Assume that $\left\vert \mu \right\vert <%
\frac{1}{4},$ $\xi >2,$ and consider%
\begin{equation*}
f(\left( \lambda ,y\right) ,x)=\left( \left( 2\lambda ,\frac{1}{4}+\frac{1}{4%
}\cos \lambda +\mu y\right) ,\xi x\right) .
\end{equation*}%
The difference is that instead of collapsing $\mathcal{M}$ completely, we
contract in the $y$ coordinate. Then $f(\mathcal{M})$ will be the set
depicted on the left plot of Figure \ref{fig:Mobius2}. The second iterate is
shown in the right plot of Figure \ref{fig:Mobius2}. We thus see that the
invariant set has a Cantor structure.
\end{Example}

Above examples are artificial. Similar features though can be found for
instance in the Kuznetzov system (see \cite{Kuz},\cite{Wil}), where we have
a hyperbolic invariant set in $\mathbb{R}^{3}$, which has a Cantor set
structure. By adding the assumption that $f$ satisfies backward cone
conditions we rule out such cases, and establish NHIMs that are graphs over $%
\Lambda $.

%TCIDATA{Version=5.00.0.2606}
%TCIDATA{LaTeXparent=0,0,MMFedit.tex}

\section{Cone evolution\label{sec:jet-evol}}

In this section we introduce the notion of "higher order cones". These will
be used to control the smoothness of established manifolds. The section
contains auxiliary results. The construction of the manifolds is performed
in Sections \ref{sec:wcu-exists}, \ref{sec:wcs-exists} and \ref%
{sec:nhim-exists}.

\subsection{Unstable cones\label{sec:unstable-cones}}

In this section we introduce the cones. We formulate the results in a
setting where we have two coordinates $\mathrm{x}$ and $\mathrm{y}$, instead
of the three coordinates $\lambda ,$ $x,$ $y$ from Section \ref%
{sec:main-results}. This is because the results are formulated in more
general terms. Later, we shall apply these taking $\mathrm{x}=\left( \lambda
,x\right) $ and $\mathrm{y}=y$ (or, in other instances, $\mathrm{x}=x$ and $%
\mathrm{y}=(\lambda ,y)$) in our construction of the manifolds. Thus, the
subtle change of font in $\mathrm{x}$ and $\mathrm{y}$ plays an important
role.

Let 
\begin{equation*}
\mathcal{P}_{m}:\mathbb{R}^{u}\rightarrow\mathbb{R}^{s},\qquad\mathcal{P}%
_{m}(0)=0.
\end{equation*}
be a polynomial of degree $m$.

\begin{definition}
We define an unstable cone of order $m$ at $z,$ spanned on $\mathcal{P}_{m},$
with a bound $M>0,$ as a set of the form%
\begin{equation}
J_{u}(z,\mathcal{P}_{m},M)=\left\{ z+\left( \mathrm{x},\mathrm{y}+\mathcal{P}%
_{m}(\mathrm{x})\right) :\quad \left\Vert \mathrm{y}\right\Vert \leq
M\left\Vert \mathrm{x}\right\Vert ^{m+1}\right\} .  \label{eq:Jet-u-m}
\end{equation}
\end{definition}

\begin{remark}
We emphasize that the index $m$ in $J_{u}(z_{0},\mathcal{P}_{m},M)$ is
important since it stands for the order $m$ of the cone. Cones of order $m$
are always associated with polynomials of degree $m$. Let us also observe
that if we take a polynomial (of degree zero) $\mathcal{P}_{0}=0$, then for $%
\mathrm{x}=x$ and $\mathrm{y}=\left( \lambda ,y\right) $ the cones defined
in (\ref{eq:unstable-cone-3d}) and (\ref{eq:Jet-u-m}) are the same:%
\begin{equation*}
J_{u}(z,\mathcal{P}_{0}=0,M)=J_{u}\left( z,M\right) .
\end{equation*}
\end{remark}

For $\delta>0$ we define 
\begin{equation*}
J_{u}(z_{0},\mathcal{P}_{m},M,\delta)=J_{u}(z_{0},\mathcal{P}_{m},M)\cap%
\overline{B}(0,\delta).
\end{equation*}

The above defined cones are devised to control higher order derivatives of
functions. The following lemmas explain this relation.

\begin{lemma}
\label{lem:func-in-unstb-jet} Assume that $g:\mathbb{R}^{u}\supset \mathrm{%
dom}(g)\rightarrow \mathbb{R}^{s}$ is a $C^{m+1}$ function. Let $\mathrm{x}%
_{0}\in \mathrm{dom}(g)$, $M>\Vert D^{m+1}g(\mathrm{x}_{0})\Vert $. Then
there exists a $\delta >0$, such that 
\begin{equation}
\{(\mathrm{x},g(\mathrm{x}))\ |\ \Vert \mathrm{x}-\mathrm{x}_{0}\Vert \leq
\delta \}\subset J_{u}(z_{0},\mathcal{P}_{m},M/\left( m+1\right) !)
\label{eq:func-in-unstb-jet}
\end{equation}%
for $z_{0}=\left( \mathrm{x}_{0},g(\mathrm{x}_{0})\right) $ and $\mathcal{P}%
_{m}(\mathrm{x})=T_{g,m,\mathrm{x}_{0}}(\mathrm{x}).$
\end{lemma}

\begin{proof}
The proof is given in \ref{app:func-in-unstb-jet}.
\end{proof}

The crucial property of $J_{u}$ is that Lemma \ref{lem:func-in-unstb-jet}
can be reversed to give bounds on the higher order derivatives:

\begin{lemma}
\label{lem:inv-func-in-unstb-jet} Assume that $g:\mathbb{R}^{u}\supset 
\mathrm{dom}(g)\rightarrow \mathbb{R}^{s}$ is a $C^{m+1}$ function. Let $%
\mathrm{x}_{0}\in \mathrm{dom}(g)$ and assume that there exists $\delta >0$,
such that 
\begin{equation}
\{(\mathrm{x},g(\mathrm{x}))\ |\ \Vert \mathrm{x}-\mathrm{x}_{0}\Vert \leq
\delta \}\subset J_{u}(z_{0},\mathcal{P}_{m},M),
\label{eq:inv-func-in-unstb-jet}
\end{equation}%
where $z_{0}=(\mathrm{x}_{0},g(\mathrm{x}_{0}))$ and $\mathcal{P}_{m}(%
\mathrm{x})=T_{g,m,\mathrm{x}_{0}}(\mathrm{x}).$ Then there exists a
constant $C$ (which depends only on $m$ and $s$), such that for any $%
j_{1},\ldots,j_{m+1}\in \{1,\ldots ,u\}$ 
\begin{equation*}
\left\Vert \frac{\partial ^{m+1}g(\mathrm{x}_{0})}{\partial \mathrm{x}%
_{i_{1}}\ldots \partial \mathrm{x}_{i_{m+1}}}\right\Vert \leq CM.
\end{equation*}
\end{lemma}

\begin{proof}
See \ref{app:inv-func-in-unstb-jet}.
\end{proof}

We now show that when $f$ satisfies certain conditions, unstable cones are
mapped into themselves. We start with a simple case of cones of order zero.

\begin{theorem}
\label{thm:Lip-unstable-jet}Let $U\subset \mathbb{R}^{u}\times \mathbb{R}%
^{s} $ be a convex neighborhood of zero and assume that $f:U\rightarrow 
\mathbb{R}^{u}\times \mathbb{R}^{s}$ is a $C^{1}$ map satisfying $f(0)=0$.
If for $M>0$ 
\begin{eqnarray}
m\left[ \frac{\partial f_{\mathrm{x}}}{\partial \mathrm{x}}(U)\right]
-M\sup_{\mathrm{x}\in U}\left\Vert \frac{\partial f_{\mathrm{x}}}{\partial 
\mathrm{y}}(z)\right\Vert &\geq &\xi ,  \label{eq:Lip-jet-cond-1} \\
\sup_{z\in U}\left\{ \left\Vert \frac{\partial f_{\mathrm{y}}}{\partial 
\mathrm{y}}(z)\right\Vert +\frac{1}{M}\left\Vert \frac{\partial f_{\mathrm{y}%
}}{\partial \mathrm{x}}(z)\right\Vert \right\} &\leq &\mu ,
\label{eq:Lip-jet-cond-2}
\end{eqnarray}%
and%
\begin{equation}
\frac{\xi }{\mu }>1,  \label{eq:Lip-jet-cond-3}
\end{equation}%
then 
\begin{equation*}
f(J_{u}(0,\mathcal{P}_{0}=0,M)\cap U)\subset \mathrm{int}J_{u}(0,\mathcal{R}%
_{0}=0,M)\cup \{0\}.
\end{equation*}
\end{theorem}

\begin{proof}
See \ref{app:Lip-unstable-jet}.
\end{proof}

The following theorem shows that, under appropriate assumptions, cones of
order $m$ map to other cones, with the same bound $M$.

\begin{theorem}
\label{thm:unstb-jet-propagation} Let $D\subset \mathbb{R}^{u}\times \mathbb{%
R}^{s}$ be a convex bounded neighborhood of zero and assume that $%
f:D\rightarrow \mathbb{R}^{u}\times \mathbb{R}^{s}$ is a $C^{m+1}$ map
satisfying $f(0)=0$ and $\Vert f(D)\Vert _{C^{m+1}}\leq C$. Assume that we
have two polynomials $\mathcal{P}_{m},\mathcal{R}_{m}:\mathbb{R}%
^{u}\rightarrow \mathbb{R}^{s}$ with coefficients bounded by $C$, such that%
\begin{equation}
\mathrm{graph}(T_{\pi _{y}f\circ (\mathrm{id},\mathcal{P}_{m}),m,0})\subset 
\mathrm{graph}\left( \mathcal{R}_{m}\right) .
\label{eq:full-unstb-graph-cond}
\end{equation}

If for $\xi >0$, and $\rho <1$%
\begin{equation}
\begin{array}{rr}
m\left( \frac{\partial f_{\mathrm{x}}}{\partial \mathrm{x}}(0)+\frac{%
\partial f_{\mathrm{x}}}{\partial \mathrm{y}}(0)D\mathcal{P}_{m}(0)\right)
\geq & \xi ,\medskip \\ 
\left\Vert \frac{\partial f_{\mathrm{x}}}{\partial \mathrm{y}}(0)\right\Vert
\leq & B,\medskip \\ 
\left\Vert \frac{\partial f_{\mathrm{y}}}{\partial \mathrm{y}}(0)-D\mathcal{R%
}_{m}(0)\frac{\partial f_{\mathrm{x}}}{\partial \mathrm{y}}(0)\right\Vert
\leq & \mu ,%
\end{array}
\label{eq:full-unstb-cond-1}
\end{equation}%
and%
\begin{equation}
\frac{\mu }{\xi ^{m+1}}<\rho ,  \label{eq:full-unstb-rate-cond}
\end{equation}%
then there exists a constant $M^{\ast }=M^{\ast }\left( C,B,1/\xi ,\rho
\right) $, such that for any $M>M^{\ast }$ there exists a $%
\delta=\delta(M,C,B,1/\xi)$ such that 
\begin{equation*}
f(J_{u}(0,\mathcal{P}_{m},M,\delta )\cap D)\subset J_{u}(0,\mathcal{R}%
_{m},M).
\end{equation*}%
Moreover, if for some $K>0$ holds $C,B,\frac{1}{\xi }\in \left[ 0,K\right] $%
, then $M^{\ast }$ depends only on $K$ and $\rho $.
\end{theorem}

\begin{proof}
See \ref{app:unstb-jet-propagation}.
\end{proof}

\subsection{Stable cones\label{sec:stable-cones}}

Let 
\begin{equation*}
\mathcal{Q}_{m}:\mathbb{R}^{s}\rightarrow \mathbb{R}^{u},\qquad \mathcal{Q}%
_{m}(0)=0.
\end{equation*}%
be a polynomial of degree $m$.

\begin{definition}
We define a stable cone of order $m$ at $z_{0},$ spanned on $\mathcal{Q}%
_{m}, $ with a bound $M>0,$ as a set of the form%
\begin{equation*}
J_{s}(z_{0},\mathcal{Q}_{m},M)=\left\{ z_{0}+\left( \mathrm{x}+\mathcal{Q}%
_{m}(\mathrm{y}),\mathrm{y}\right) :\quad \left\Vert \mathrm{x}\right\Vert
\leq M\left\Vert \mathrm{y}\right\Vert ^{m+1}\right\} .
\end{equation*}
\end{definition}

For $\delta >0$ we define 
\begin{equation*}
J_{s}(z_{0},\mathcal{Q}_{m},M,\delta )=J_{s}(z_{0},\mathcal{Q}_{m},M)\cap 
\overline{B}(0,\delta ),
\end{equation*}%
and we also denote complements of the cones as%
\begin{eqnarray*}
J_{s}^{c}(z_{0},\mathcal{Q}_{m},M) &=&\mathbb{R}^{u}\times \mathbb{R}%
^{s}\setminus J_{s}(z_{0},\mathcal{Q}_{m},M), \\
J_{s}^{c}(z_{0},\mathcal{Q}_{m},M,\delta ) &=&\overline{B}(0,\delta
)\setminus J_{s}(z_{0},\mathcal{Q}_{m},M,\delta ).
\end{eqnarray*}

Mirror results to Lemmas \ref{lem:func-in-unstb-jet}, \ref%
{lem:inv-func-in-unstb-jet} can be formulated for stable cones:

\begin{lemma}
\label{lem:func-in-stb-jet} Assume that $g:\mathbb{R}^{s}\supset \mathrm{dom}%
(g)\rightarrow \mathbb{R}^{u}$ is a $C^{m+1}$ function. Let $\mathrm{y}%
_{0}\in \mathrm{dom}(g)$, $M>\Vert D^{m+1}g(\mathrm{y}_{0})\Vert $. Then
there exists $\delta >0$, such that 
\begin{equation*}
\{(g(\mathrm{y}),\mathrm{y})\ |\ \Vert \mathrm{y}-\mathrm{y}_{0}\Vert \leq
\delta \}\subset J_{s}(z_{0},\mathcal{P}_{m},M/\left( m+1\right) !)
\end{equation*}%
for $z_{0}=\left( g(\mathrm{y}_{0}),\mathrm{y}_{0}\right) $ and $\mathcal{P}%
_{m}(\mathrm{y})=T_{g,m,\mathrm{y}_{0}}(\mathrm{y}).$
\end{lemma}

\begin{lemma}
\label{lem:inv-func-in-stb-jet} Assume that $g:\mathbb{R}^{s}\supset \mathrm{%
dom}(g)\rightarrow \mathbb{R}^{u}$ is a $C^{m+1}$ function. Let $\mathrm{y}%
_{0}\in \mathrm{dom}(g)$ and assume that there exists $\delta >0$, such that 
\begin{equation*}
\{(g(\mathrm{y}),\mathrm{y})\ |\ \Vert \mathrm{y}-\mathrm{y}_{0}\Vert \leq
\delta \}\subset J_{u}(z_{0},\mathcal{P}_{m},M),
\end{equation*}%
where $z_{0}=(g(\mathrm{y}_{0}),\mathrm{y}_{0})$ and $\mathcal{P}_{m}(%
\mathrm{y})=T_{g,m,\mathrm{y}_{0}}(\mathrm{y}).$ Then there exists a
constant $C$ (which depends only on $m$), such that for any $j_{1},\ldots
,j_{m+1}\in \{1,\ldots ,u\}$ 
\begin{equation*}
\left\Vert \frac{\partial ^{m+1}g(\mathrm{y}_{0})}{\partial \mathrm{y}%
_{i_{1}}\ldots \partial \mathrm{y}_{i_{m+1}}}\right\Vert \leq CM.
\end{equation*}
\end{lemma}

Since proofs of Lemmas \ref{lem:func-in-stb-jet}, \ref%
{lem:inv-func-in-stb-jet} follow from mirror arguments to the proofs of
Lemmas \ref{lem:func-in-unstb-jet}, \ref{lem:inv-func-in-unstb-jet}, we omit
their proofs.

We now give the following theorems, which are in similar spirit to Theorem~%
\ref{thm:Lip-unstable-jet}, \ref{thm:unstb-jet-propagation}. The difference
is that they concern images of \textbf{complements} of cones (and not images
of the cones themselves, as is the case in Theorems \ref%
{thm:Lip-unstable-jet}, \ref{thm:unstb-jet-propagation}.)

\begin{theorem}
\label{thm:Lip-stable-jet}Let $U\subset \mathbb{R}^{u}\times \mathbb{R}^{s}$
be a convex neighborhood of zero and assume that $f:U\rightarrow \mathbb{R}%
^{u}\times \mathbb{R}^{s}$ is a $C^{1}$ map satisfying $f(0)=0$. Assume that
for $M>0$ 
\begin{eqnarray}
m\left( \left[ \frac{\partial f_{\mathrm{x}}}{\partial \mathrm{x}}(U)\right]
\right) -\frac{1}{M}\sup_{z\in U}\left\Vert \frac{\partial f_{\mathrm{x}}}{%
\partial \mathrm{y}}(z)\right\Vert &\geq &\xi ,  \label{eq:Lip-jet-cond-4} \\
\sup_{z\in U}\left\{ \left\Vert \frac{\partial f_{\mathrm{y}}}{\partial 
\mathrm{y}}(z)\right\Vert +M\left\Vert \frac{\partial f_{\mathrm{y}}}{%
\partial \mathrm{x}}(z)\right\Vert \right\} &\leq &\mu ,
\label{eq:Lip-jet-cond-5}
\end{eqnarray}%
and%
\begin{equation}
\frac{\xi }{\mu }>1,  \label{eq:Lip-jet-cond-6}
\end{equation}%
then 
\begin{equation*}
f(\overline{J_{s}^{c}(0,\mathcal{Q}_{0}=0,M)}\cap U)\subset J_{s}^{c}(0,%
\mathcal{R}_{0}=0,M)\cup \{0\}.
\end{equation*}
\end{theorem}

\begin{proof}
The result follows from Theorem \ref{thm:Lip-unstable-jet}. Details are
given in \ref{app:Lip-stable-jet}.
\end{proof}

\begin{theorem}
\label{thm:stb-jet-propagation} Let $D\subset \mathbb{R}^{u}\times \mathbb{R}%
^{s}$ be a convex bounded neighborhood of zero and assume that $%
f:D\rightarrow \mathbb{R}^{u}\times \mathbb{R}^{s}$ is a $C^{m+1}$ map
satisfying $f(0)=0$ and $\Vert f(D)\Vert _{C^{m+1}}\leq C$. Assume that we
have two polynomials $\mathcal{Q}_{m},\mathcal{R}_{m}:\mathbb{R}%
^{s}\rightarrow \mathbb{R}^{u}$ with coefficients bounded by $C$, such that%
\begin{equation*}
\mathrm{graph}(T_{\pi _{x}f\circ (\mathcal{Q}_{m},\mathrm{id}),m,0})\subset 
\mathrm{graph}\left( \mathcal{R}_{m}\right) .
\end{equation*}

If for $\xi >0$, and $\rho <1$%
\begin{equation*}
\begin{array}{rr}
m\left( \frac{\partial f_{\mathrm{x}}}{\partial \mathrm{x}}(0)+D\mathcal{P}%
_{m}(0)\frac{\partial f_{\mathrm{y}}}{\partial \mathrm{x}}(0)\right) \geq & 
\xi ,\medskip \\ 
\left\Vert \frac{\partial f_{\mathrm{y}}}{\partial \mathrm{x}}(0)\right\Vert
\leq & B,\medskip \\ 
\left\Vert \frac{\partial f_{\mathrm{y}}}{\partial \mathrm{y}}(0)-\frac{%
\partial f_{\mathrm{y}}}{\partial \mathrm{x}}(0)D\mathcal{Q}%
_{m}(0)\right\Vert \leq & \mu ,%
\end{array}%
\end{equation*}%
and%
\begin{equation}
\frac{\mu ^{m+1}}{\xi }<\rho ,  \label{eq:full-stb-rate-cond}
\end{equation}%
then there exists a constant $M^{\ast }=M^{\ast }\left( C,B,1/\xi ,\rho
\right) $, such that for any $M>M^{\ast }$ there exists $\delta =\delta (M)$
such that\textbf{\ }$\delta =\delta (M,C,B,1/\xi )$ 
\begin{equation*}
f(J_{s}^{c}(0,\mathcal{Q}_{m},M,\delta )\cap U)\subset J_{s}^{c}(0,\mathcal{R%
}_{m},M).
\end{equation*}%
Moreover, if for some $K>0$ holds $C,B,\frac{1}{\xi }\in \left[ 0,K\right] $%
, then $M^{\ast }$ depends only on $K$ and $\rho $.
\end{theorem}

\begin{proof}
The proof is given in \ref{app:stb-jet-propagation}.
\end{proof}

\subsection{Center-stable and center-unstable cones}

We now return to the setting in which we have three coordinates $\left(
\lambda ,x,y\right) $. Recall that in these coordinates stable cones $%
J_{s}\left( z,M\right) $ and unstable cones $J_{u}\left( z,M\right) $ were
defined using (\ref{eq:stable-cone-3d}--\ref{eq:unstable-cone-3d}). In
addition we define center-stable and center-unstable cones as%
\begin{eqnarray*}
J_{cs}\left( z,M\right) &=&\left\{ \left( \lambda ,x,y\right) :\left\Vert
x-\pi _{x}z\right\Vert <M\left\Vert \left( \lambda ,y\right) -\pi _{\lambda
,y}z\right\Vert \right\} \cup \{z\}, \\
J_{cu}\left( z,M\right) &=&\left\{ \left( \lambda ,x,y\right) :\left\Vert
y-\pi _{y}z\right\Vert <M\left\Vert \left( \lambda ,x\right) -\pi _{\lambda
,x}z\right\Vert \right\} \cup \{z\},
\end{eqnarray*}%
respectively.

Observe that $J_{cs}(z,M)=J_{u}^{c}(z,1/M)\cup \{z\}$ and $%
J_{cu}(z,M)=J_{s}^{c}(z,1/M)\cup \{z\}$. We see that $J_{cs}(z,M)$ and $%
J_{cu}(z,M)$ as defined above are not contained in domain of single good
chart. However we will always take intersections of these cones with the
domain of a good chart.

As in section \ref{sec:main-results}, we consider a $C^{k+1}$ map, with $%
k\geq 0$, 
\begin{equation*}
f:D\rightarrow \Lambda \times \mathbb{R}^{u}\times \mathbb{R}^{s},
\end{equation*}%
where $D=\Lambda \times \overline{B}_{u}(R)\times \overline{B}_{s}(R)$. We
rewrite some of the results from sections \ref{sec:unstable-cones}, \ref%
{sec:stable-cones} in terms of coordinates $(\lambda ,x,y)$, formulating
them as corollaries.

\begin{corollary}
\label{cor:unstable-lip}If $f$ satisfies the rate conditions of order $k=0$
(see Definition \ref{def:rate-conditions}) then for any $z\in D$%
\begin{equation*}
f\left( J_{u}\left( z,1/L\right) \cap D\right) \subset \mathrm{int}%
J_{u}\left( f(z),1/L\right) \cup \{f(z)\}.
\end{equation*}%
In alternative notation, $f\left( J_{cs}^{c}\left( z,L\right) \cap D\right)
\subset \mathrm{int}J_{cs}^{c}\left( f(z),L\right) \cup \{f(z)\}.$
\end{corollary}

\begin{proof}
This follows from Theorem \ref{thm:Lip-unstable-jet}, taking coordinates $%
\mathrm{x}=x,$ $\mathrm{y}=(\lambda ,y)$ and constants $M=1/L,$ $\xi =\xi
_{u,1,P},$ $\mu =\mu _{cs,1}$. The assumption (\ref{eq:Lip-jet-cond-3}) of
Theorem \ref{thm:Lip-unstable-jet} follows from the rate condition (\ref%
{eq:rate-cond-1}).
\end{proof}

\begin{corollary}
\label{cor:stable-lip}If $f$ satisfies the rate conditions of order $k=0$
then for any $z\in D$%
\begin{equation*}
f\left( \overline{J_{cu}\left( z,L\right) }\cap \overline{B}_{c}(\pi
_{\lambda }z,R_{\Lambda })\times \overline{B}_{u}(R)\times \overline{B}%
_{s}(R)\right) \subset J_{cu}\left( f(z),L\right) .
\end{equation*}%
In alternative notation, 
\begin{equation*}
f(\overline{J_{s}^{c}(z,1/L)}\cap \overline{B}_{c}(\pi _{\lambda
}z,R_{\Lambda })\times \overline{B}_{u}(R)\times \overline{B}_{s}(R))\subset
J_{s}^{c}(f(z),1/L)\cup \left\{ f(z)\right\} .
\end{equation*}
\end{corollary}

\begin{proof}
This follows from Theorem \ref{thm:Lip-stable-jet}, taking coordinates $%
\mathrm{x}=(\lambda ,x),$ $\mathrm{y}=y$ and constants $M=1/L,$ $\xi =\xi
_{cu,1,P},$ $\mu =\mu _{s,1}$.
\end{proof}

\begin{lemma}
\label{lem:Ju-expansion}If $f$ satisfies the rate conditions of order $k=0$
and two points $z_{1},z_{2}\in D$ satisfy $z_{1}\in J_{u}\left(
z_{2},1/L\right) $, then 
\begin{equation*}
\left\Vert \pi _{x}\left( f(z_{1})-f(z_{2})\right) \right\Vert \geq \xi
_{u,1,P}\left\Vert \pi _{x}(z_{1}-z_{2})\right\Vert .
\end{equation*}
\end{lemma}

\begin{proof}
See \ref{app:Ju-expansion}.
\end{proof}

\begin{lemma}
\label{lem:Js-contraction}If $f$ satisfies the rate conditions of order $k=0$
and two points $z_{1},z_{2}\in D$, satisfy $z_{1}\in J_{s}\left(
z_{2},1/L\right) $ and $f(z_{1})\in J_{s}\left( f(z_{2}),1/L\right) $, then%
\begin{equation*}
\left\Vert \pi _{y}\left( f(z_{1})-f(z_{2})\right) \right\Vert \leq \mu
_{s,1}\left\Vert \pi _{y}(z_{1}-z_{2})\right\Vert .
\end{equation*}
\end{lemma}

\begin{proof}
See \ref{app:Js-contraction}.
\end{proof}

\begin{lemma}
\label{lem:cu-cone-xi-cu}Assume that $z_{1},z_{2}$ are in the same chart. If 
$f$ satisfies the rate conditions of order $k=0$ and $z_{1}\in J_{cu}\left(
z_{2},L\right) $, then%
\begin{equation*}
\left\Vert \pi _{\left( \lambda ,x\right) }\left( f(z_{1})-f(z_{2})\right)
\right\Vert \geq \xi _{cu,1,P}\left\Vert \pi _{\left( \lambda ,x\right)
}\left( z_{1}-z_{2}\right) \right\Vert .
\end{equation*}
\end{lemma}

\begin{proof}
See \ref{app:cu-cone-xi-cu}.
\end{proof}

\begin{lemma}
\label{lem:cs-cone-mu-cs}Assume that $z_{1},z_{2}$ and $f(z_{1}),f(z_{2})$
are in the same charts. If $f$ satisfies the rate conditions of order $k=0$
and $z_{1}\in J_{cs}(z_{2},L)$, then%
\begin{equation*}
\left\Vert \pi _{\left( \lambda ,y\right) }\left( f(z_{1})-f(z_{2})\right)
\right\Vert \leq \mu _{cs,1}\left\Vert \pi _{\left( \lambda ,y\right)
}\left( z_{1}-z_{2}\right) \right\Vert .
\end{equation*}
\end{lemma}

\begin{proof}
See \ref{app:cs-cone-mu-cs}.
\end{proof}

%TCIDATA{Version=5.00.0.2606}
%TCIDATA{LaTeXparent=0,0,MMFedit.tex}

\section{Discs}

In this section we introduce the notion of discs. These will be the building
blocks for the construction of our invariant manifolds.

\begin{definition}
We say that a continuous function $b:\overline{B}_{u}(R)\rightarrow D$ is a
horizontal disc if for any $x\in \overline{B}_{u}(R)$%
\begin{equation}
\pi _{x}b(x)=x\qquad \text{and\qquad }b\left( \overline{B}_{u}(R)\right)
\subset J_{u}\left( b(x),1/L\right) .  \label{eq:hor-disc-prop}
\end{equation}
\end{definition}

\begin{definition}
We say that a continuous function $b:\overline{B}_{s}(R)\rightarrow D$ is a
vertical disc if for any $y\in \overline{B}_{s}(R)$%
\begin{equation}
\pi _{y}b(y)=y\qquad \text{and\qquad }b\left( \overline{B}_{s}(R)\right)
\subset J_{s}\left( b(y),1/L\right) .  \label{eq:ver-disc-prop}
\end{equation}
\end{definition}

By Remark \ref{rem:cone-in-chart}, we see that any horizontal or vertical
disc can be contained in a set on which we can use a single chart. This fact
will prove important in Section \ref{sec:bundles} where we reformulate our
results for more general $\Lambda $.

In our former works \cite{CZ,Z} the disks as defined above where said to
satisfy cone conditions.

\begin{definition}
We say that a continuous function $b:\Lambda \times \overline{B}%
_{u}(R)\rightarrow D$ is a center-horizontal disc if for any $\left( \lambda
,x\right) \in \Lambda \times \overline{B}_{u}(R)$%
\begin{equation*}
\pi _{\left( \lambda ,x\right) }b(\lambda ,x)=\left( \lambda ,x\right)
\end{equation*}%
and 
\begin{equation}
b\left( \overline{B}_{c}(\lambda ,R_{\Lambda })\times \overline{B}%
_{u}(R)\right) \subset J_{cu}\left( b(\lambda ,x),L\right) .
\label{eq:ch-cone-cond}
\end{equation}
\end{definition}

\begin{definition}
We say that a continuous function $b:\Lambda \times \overline{B}%
_{s}(R)\rightarrow D$ is a center-vertical disc if for any $\left( \lambda
,y\right) \in \Lambda \times \overline{B}_{s}(R)$%
\begin{equation*}
\pi _{\left( \lambda ,y\right) }b(\lambda ,y)=\left( \lambda ,y\right)
\end{equation*}%
and%
\begin{equation}
b\left( \overline{B}_{c}(\lambda ,R_{\Lambda })\times \overline{B}%
_{s}(R)\right) \subset J_{cs}\left( b(\lambda ,y),L\right).
\label{eq:cv-cone-cond}
\end{equation}
\end{definition}

\begin{lemma}
\label{lem:unstable-disc}Assume that $b:\overline{B}_{u}(R)\rightarrow D$ is
a horizontal disc. If $f$ satisfies the covering conditions and the rate
conditions of order $l = 0$ , then there exists a horizontal disc $b^{\ast }:%
\overline{B}_{u}(R)\rightarrow D$ such that $f\circ b(\overline{B}%
_{u}(R))\cap D= b^{\ast }(\overline{B}_{u}(R))$. Moreover, if $f$ and $b$
are $C^{k} $, then so is $b^{\ast }$.
\end{lemma}

\begin{proof}
The proof is given in appendix \ref{app:unstable-disc}.
\end{proof}

The disc $b^{\ast }$ from Lemma \ref{lem:unstable-disc} is a graph transform
of $b$. From now on we shall use the notation $\mathcal{G}_{h}(b)$ instead
of $b^{\ast }$.

\begin{lemma}
\label{lem:center-unstable-disc}Assume that $b:\Lambda \times \overline{B}%
_{u}(R)\rightarrow D$ is a center-horizontal disc. If $f$ satisfies the
covering conditions, backward cone conditions and the rate conditions of
order $l = 0$, then there exists a center-horizontal disc $b^{\ast }:\Lambda
\times \overline{B}_{u}(R)\rightarrow D$ such that 
\begin{equation*}
f\circ b(\Lambda \times \overline{B}_{u}(R))\cap D= b^{\ast }(\Lambda \times 
\overline{B}_{u}(R)).
\end{equation*}%
Moreover, if $f$ and $b$ are $C^{k}$, then so is $b^{\ast }$.
\end{lemma}

\begin{proof}
The proof is given in appendix \ref{app:center-unstable-disc}.
\end{proof}

From now on we shall use the notation $\mathcal{G}_{ch}(b)$ instead of $%
b^{\ast }$ for the disc from Lemma \ref{lem:center-unstable-disc}.

%TCIDATA{Version=5.00.0.2606}
%TCIDATA{LaTeXparent=0,0,MMFedit.tex}

\section{Center-unstable manifold\label{sec:wcu-exists}}

In this section we prove the existence and smoothness the manifold $W^{cu}$
from Theorem \ref{th:main}. The proof follows from a graph transform type
method, in which we take successive iterates of center-horizontal discs, and
these converge to the center unstable manifold.

We start with the following lemma, which establishes the existence of $%
W^{cu} $.

\begin{lemma}
\label{lem:w-cu}Assume that $f$ satisfies covering conditions, backward cone
conditions and rate conditions of order $l\geq 0$. Let $b_{i}$ be the
sequence of center-horizontal discs defined as $b_{0}(\lambda ,x)=(\lambda
,x,0)$, $b_{i+1}=\mathcal{G}_{ch}(b_{i})$ for $i>0$. Then $b_{i}$ converge
uniformly to a center-horizontal disc $\left( \lambda ,x\right) \rightarrow
\left( \lambda ,x,w^{cu}(\lambda ,x)\right) $, where%
\begin{equation*}
w^{cu}:\Lambda \times \overline{B}_{u}(R)\rightarrow \overline{B}_{s}(R).
\end{equation*}%
Moreover 
\begin{equation*}
W^{cu}=\left\{ \left( \lambda ,x,w^{cu}(\lambda ,x)\right) :\Lambda \times 
\overline{B}_{u}(R)\right\} .
\end{equation*}
\end{lemma}

\begin{proof}
We use a notation $\theta =(\lambda ,x)$. We will show that $\pi_{y}b_{i}$
is a Cauchy sequence in the supremum norm, which converges to $W^{cu}$.

Let us fix $\theta \in \Lambda \times \overline{B}_{u}(R)$. For any $k\in 
\mathbb{N}$, since $b_{k}$ is center-horizontal disk, there exists a finite
backward orbit $\{q_{i}^{k}\}_{i=-k,\dots ,0}$, such that 
\begin{equation*}
q_{0}^{k}=b_{k}(\theta ),\quad \pi _{\theta }(q_{0}^{k})=\theta .
\end{equation*}%
From the backward cone condition it follows that for any $i<0$ points $%
\{q_{i}^{k}\}$ for $k\geq |i|$ are in the same chart and 
\begin{equation*}
q_{i}^{k_{1}}\in J_{s}(q_{i}^{k_{2}},1/L),\qquad k_{1},k_{2}\geq |i|.
\end{equation*}%
Therefore we have 
\begin{equation}
\Vert q_{i}^{k_{1}}-q_{i}^{k_{2}}\Vert \leq \left( 1+1/L\right) \Vert \pi
_{y}(q_{i}^{k_{1}}-q_{i}^{k_{2}})\Vert .  \label{eq:dist-in-Js}
\end{equation}%
From Lemma~\ref{lem:Js-contraction} it follows that for $j\in \mathbb{Z}%
_{-}\cup \{0\},$ and $k_{2}>k_{1}\geq |j|$ holds%
\begin{eqnarray}
\Vert \pi _{y}q_{j}^{k_{1}}-\pi _{y}q_{j}^{k_{2}}\Vert &=&\Vert \pi
_{y}f^{k_{1}+j}(q_{-k_{1}}^{k_{1}})-\pi
_{y}f^{k_{1}+j}(q_{-k_{1}}^{k_{2}})\Vert \leq  \label{eq:Wcu-back-orb-estm}
\\
&\leq &(\mu _{s,1})^{k_{1}+j}\Vert \pi
_{y}(q_{-k_{1}}^{k_{1}}-q_{-k_{1}}^{k_{2}})\Vert \leq 2R(\mu
_{s,1})^{k_{1}+j}.  \notag
\end{eqnarray}%
From (\ref{eq:Wcu-back-orb-estm}) and (\ref{eq:dist-in-Js}) it follows that
for each $j\in \mathbb{Z}_{-}\cup \{0\}$ holds 
\begin{equation}
\Vert q_{j}^{k_{1}}-q_{j}^{k_{2}}\Vert \leq (1+1/L)2R(\mu
_{s,1})^{k_{1}+j},\qquad k_{2}>k_{1}\geq |j|.  \label{eq:Wcu-back-orb-dist}
\end{equation}%
Since $q_{0}^{k}=b_{k}(\theta )$ condition (\ref{eq:Wcu-back-orb-dist})
establishes uniform convergence of $b_{k}$ to $b^{\ast }$, moreover also the
backward orbits form a Cauchy sequence and converge to full backward orbit
of $b^{\ast }(\theta )$.

From the above it follows also that $\pi _{y}b_{k}$ converge uniformly to a
continuous function $w^{cu}(\theta )=\pi _{y}b^{\ast }(\theta )$.

Assume now that we have a $z\in D$ that has a full backward trajectory $%
\{z_{k}\}_{k=-\infty }^{0}$ in $D$. We need to show that $z=\left( \theta
^{\ast },w^{cu}(\theta ^{\ast })\right) $ for some $\theta ^{\ast }\in
\Lambda \times \overline{B}_{u}(R)$. Let $z^{\ast }=\left( \theta ^{\ast
},w^{cu}(\theta ^{\ast })\right) $ for $\theta ^{\ast }=\pi _{\theta }z$. We
will show that $z=z^{\ast }$. Since $\pi _{\theta }z=\pi _{\theta }z^{\ast }$%
, 
\begin{equation*}
z\in J_{s}(z^{\ast },1/L).
\end{equation*}%
Then for the backward trajectory $\{z_{k}^{\ast }\}_{k=-\infty }^{0}$ of $%
z^{\ast }$, by the backward cone conditions,%
\begin{equation*}
z_{k}\in J_{s}(z_{k}^{\ast },1/L)\qquad \text{for }k=0,-1,-2,\ldots \text{ .}
\end{equation*}%
By Lemma \ref{lem:Js-contraction}, this implies that 
\begin{equation*}
\left\Vert \pi _{y}(z^{\ast }-z)\right\Vert \leq \left( \mu _{s,1}\right)
^{k}\left\Vert \pi _{y}\left( z_{k}^{\ast }-z_{k}\right) \right\Vert \leq
2\left( \mu _{s,1}\right) ^{k}.
\end{equation*}%
Since $\mu _{s,1}<1,$ we see that $z^{\ast }=z.$

Passing to the limit in the cone condition (\ref{eq:ch-cone-cond}) for $b_{k}
$ one can see that%
\begin{equation*}
b^{\ast }\left( \overline{B}_{c}\left( \lambda ,R_{\Lambda }\right) \times 
\overline{B}_{u}\left( R\right) \right) \subset \overline{J_{cu}\left(
b\left( \lambda ,x\right) \right) }.
\end{equation*}%
Since $W^{cu}$ is invariant under $f$, by Corollary \ref{cor:stable-lip}%
\label{proof-tu-zmieniony-cor} we obtain (\ref{eq:ch-cone-cond}) for $%
b^{\ast }.$ Thus $b^{\ast }$ is a center-horizontal disc.
\end{proof}

\begin{lemma}
\label{lem:Cm-norm-bound}Assume that $f$ is $C^{k+1}$ and satisfies covering
conditions, backward cone conditions and rate conditions of order $l\geq 0$.
Let $m\leq k$. Let $b_{i}$ be the sequence of center-horizontal discs
defined as $b_{0}(\lambda ,x)=(\lambda ,0,0)$, $b_{i+1}=\mathcal{G}%
_{ch}(b_{i})$ for $i=0,1,2,\ldots $. Assume that $b_{i}$ are $C^{m}$ and
that for any $i,$ $\left\Vert \pi _{y}b_{i}\right\Vert _{C^{m}}<c_{m},$ with 
$c_{m}$ independent of $i$. If the order $l$ of the rate conditions is
greater or equal to $m$, then $\left\Vert \pi _{y}b_{i}\right\Vert
_{C^{m+1}}<c_{m+1} $ for a constant independent of $i$.
\end{lemma}

\begin{proof}
Let us fix $i\in \mathbb{N}$. Our aim will be to show that $\left\Vert \pi
_{y}b_{i}\right\Vert _{C^{m+1}}$ is bounded and that the bound is
independent of $i$. Let $\theta _{i}$ be any chosen point from $\Lambda
\times \overline{B}_{u}(R)$ and let $\theta _{0},\ldots ,\theta _{i}\in
\Lambda \times \overline{B}_{u}(R)$ be a sequence such that%
\begin{equation*}
b_{l+1}(\theta _{l+1})=f(b_{l}(\theta _{l})),
\end{equation*}%
for $l=0,\ldots ,i-1$. Note that%
\begin{equation*}
\theta _{l+1}=\pi _{\theta }f(b_{l}(\theta _{l})).
\end{equation*}%
For $l=0,\ldots ,i$ let $\mathcal{P}_{m}^{l}:\mathbb{R}^{c+u}\supset B\left(
0,\delta \right) \rightarrow \mathbb{R}^{s}$ be a polynomial of degree $m$,
defined as 
\begin{equation*}
\mathcal{P}_{m}^{l}=\pi _{y}b_{l}(\theta _{l})+T_{\pi _{y}b_{l},m,\theta
_{l}}.
\end{equation*}%
Observe that since $\left\Vert \pi _{y}b_{l}\right\Vert _{C^{m}}<c_{m}$ for $%
c_{m}$ independent from $l$, the polynomials $\mathcal{P}_{m}^{l}$ have a
uniform bounds for their coefficients, which is independent from $l$ and $i$%
. Since $b_{l+1}=\mathcal{G}_{ch}(b_{l})$ we also see that for $l=0,\ldots
,i-1$ 
\begin{eqnarray*}
\mathrm{graph}(T_{\pi _{y}f\circ (\mathrm{id},\mathcal{P}_{m}^{l}),m,0}) &=&%
\mathrm{graph}(T_{\pi _{y}f\circ b_{l},m,0}) \\
&=&\mathrm{graph}(T_{\pi _{y}b_{l+1},m,0}) \\
&=&\mathrm{graph}\left( \mathcal{P}_{m}^{l+1}\right) .
\end{eqnarray*}

Since $\pi _{y}b_{l}$ are Lipschitz with a constant $L$%
\begin{equation*}
\left\Vert D\mathcal{P}_{m}^{l}(\theta )\right\Vert \leq L.
\end{equation*}%
Let us consider coordinates $\mathrm{x}=(\lambda ,x)$ $\mathrm{y}=y$ and
let, $\xi =\xi _{cu,1}$, $\mu =\mu _{s,2}$ and $B=\left\Vert f\right\Vert
_{C^{1}}$. Then, from Theorem~\ref{thm:unstb-jet-propagation}, for
sufficiently large $M$ and sufficiently small $\delta $%
\begin{equation}
f(J_{cu}(b_{l}(\theta _{l}),\mathcal{P}_{m}^{l},M,\delta ))\subset
J_{cu}(b_{l+1}(\theta _{l+1}),\mathcal{P}_{m}^{l+1},M).
\label{eq:cone-alignment-cu}
\end{equation}%
Note that the choice of $M$ and $\delta $ does not depend on $l$. Since $%
b_{0}$ is a flat disc, we have 
\begin{equation*}
b_{0}(B_{c+u}(\theta _{0},\delta ))\subset J_{cu}(b_{0}(\theta _{0}),%
\mathcal{P}_{m}^{0}=0,M,\delta ).
\end{equation*}%
This by (\ref{eq:cone-alignment-cu}) implies that 
\begin{equation}
b_{l}(B_{c+u}(\theta _{l},\delta ))\subset J_{cu}(b_{l}(\theta _{l}),%
\mathcal{P}_{m}^{l},M,\delta ).  \label{eq:bk-cu-in-cone}
\end{equation}

From (\ref{eq:bk-cu-in-cone}), by Lemma \ref{lem:inv-func-in-unstb-jet}, we
obtain a uniform bound 
\begin{equation*}
\left\Vert \frac{\partial ^{m+1}\pi _{y}b_{i}(\theta _{i})}{\partial \mathrm{%
x}_{i_{1}}\ldots \partial \mathrm{x}_{i_{m+1}}}\right\Vert \leq CM,
\end{equation*}%
where $C$ is independent of $\theta _{i}$. This means that%
\begin{equation*}
\left\Vert \pi _{y}b_{i}\right\Vert _{C^{m+1}}\leq c_{m+1},
\end{equation*}%
where $c_{m+1}$ depends on $CM$ and $c_{m}$, but is independent of $i$.
\end{proof}

\begin{lemma}
\label{lem:cu-smooth}If $f$ satisfies the assumptions of Theorem \ref%
{th:main}, then the manifold $W^{cu}$ is $C^{k}.$
\end{lemma}

\begin{proof}
Let $b_{i}$ be the sequence of center-horizontal discs defined as $%
b_{0}(\lambda ,x)=(\lambda ,0,0)$, $b_{i+1}=\mathcal{G}_{ch}(b_{i})$, for $%
i>0$. By Lemma \ref{lem:center-unstable-disc}, we know that $b_{i}$ are $%
C^{k+1}$. By Lemma \ref{lem:w-cu} we know that they converge uniformly to 
\begin{equation*}
W^{cu}=\{\left( \theta ,w^{cu}(\theta )\right) :\theta \in \Lambda \times 
\overline{B}_{u}\}.
\end{equation*}%
We need to show that $C^{k}$ smoothness is preserved as we pass to the limit.

Since $\pi _{y}b_{i}$ are Lipschitz with a constant $L$, we see that $%
\left\Vert \pi _{y}b_{i}\right\Vert _{C^{1}}\leq c_{1}$, where $c_{1}$ is
independent of $i$. Rate conditions of order $k$, imply rate conditions of
order $m$ for $m\leq k$; in particular for $m=1$. Hence, by Lemma \ref%
{lem:Cm-norm-bound} we obtain that $\left\Vert \pi _{y}b_{i}\right\Vert
_{C^{2}}\leq c_{2}$.

Applying Lemma \ref{lem:Cm-norm-bound} inductively we obtain that $%
\left\Vert \pi _{y}b_{i}\right\Vert _{C^{k+1}}\leq c_{k+1}$, with $c_{k+1}$
independent of $i$. This implies that derivatives of $\pi _{y}b_{i}$ of
order smaller or equal to $k$ are uniformly bounded and uniformly
equicontinuous. This by the Arzela Ascoli theorem implies that $\pi
_{y}b_{i} $ and their derivatives of order smaller or equal to $k$ converge
uniformly. Thus $w^{cu}$ is $C^{k}$, as required.
\end{proof}

\begin{lemma}
\label{lem:fWu-injective}If $f$ satisfies the assumptions of Theorem \ref%
{th:main}, then $f|_{W^{cu}}$ is injective.
\end{lemma}

\begin{proof}
If $p_{1},p_{2}\in D$ and $f(p_{1})=f(p_{2})$ then $f(p_{1})\in
J_{s}(f(p_{2}),1/L)$ and by the backward cone conditions $p_{1}\in
J_{s}(p_{2},1/L)$ hence by Remark \ref{rem:cone-in-chart}, $p_{1}$ and $%
p_{2} $ are in the same chart. This means that it is enough to show that for
any $p_{1},p_{2}\in W^{cu}$ which are on the same chart, we can not have $%
p_{1}\neq p_{2}$ and $f(p_{1})=f(p_{2})$.

Let $\theta _{1}\in \Lambda \times \overline{B}_{u}(R)$, $\theta _{2}\in 
\overline{B}_{c}(\pi _{\lambda }\theta ,R_{\Lambda })\times \overline{B}%
_{u}(R)$ and $\theta _{1}\neq \theta _{2}$. By Corollary~\ref%
{lem:Ju-expansion} it follows that 
\begin{equation*}
\left\Vert \pi _{\theta }f\left( \theta _{1},w^{cu}\left( \theta _{1}\right)
\right) -\pi _{\theta }f\left( \theta _{2},w^{cu}\left( \theta _{2}\right)
\right) \right\Vert \geq \xi _{cu,1,P}\left\Vert \theta _{1}-\theta
_{2}\right\Vert .
\end{equation*}%
This implies that $f\left( \theta _{1},w^{cu}\left( \theta _{1}\right)
\right) \neq f\left( \theta _{2},w^{cu}\left( \theta _{2}\right) \right) ,$
as required.
\end{proof}

Lemmas \ref{lem:w-cu}, \ref{lem:cu-smooth} and \ref{lem:fWu-injective}
combined, prove the assertion about $W^{cu}$ from Theorem \ref{th:main}.

We finish the section by proving Theorem \ref{th:wcu}.

\begin{proof}[Proof of Theorem \protect\ref{th:wcu}]
\label{proof-lip-M-bound}The result follows from showing that%
\begin{equation}
\left\Vert \pi _{y}\left( b_{i}\left( \theta _{1}\right) -b_{i}\left( \theta
_{2}\right) \right) \right\Vert \leq M\left\Vert \pi _{\theta }\left(
b_{i}\left( \theta _{1}\right) -b_{i}\left( \theta _{2}\right) \right)
\right\Vert .  \label{eq:lip-M-bound-for-Wcu}
\end{equation}%
By definition of $b_{0},$ (\ref{eq:lip-M-bound-for-Wcu}) clearly holds for $%
i=0.$ To prove (\ref{eq:lip-M-bound-for-Wcu}) for all $i\in \mathbb{N}$, one
can inductively apply the same argument as the one from the proof of Theorem %
\ref{thm:Lip-stable-jet} (page \pageref{app:Lip-stable-jet}). Passing with $%
i $ to infinity we obtain our claim.
\end{proof}

%TCIDATA{Version=5.00.0.2606}
%TCIDATA{LaTeXparent=0,0,MMFedit.tex}

\section{Center-stable manifold \label{sec:wcs-exists}}

The goal of this section is to establish the existence of the center stable
manifold $W^{cs}$ from Theorem \ref{th:main}.

We will represent $W^{cs}$ as a limit of graphs of smooth functions. Here we
take the first step in this direction. For any $i\in \mathbb{Z}_{+}$ and $%
(\lambda ,y)\in \Lambda \times \overline{B}_{s}(R)$ we consider the
following problem: Find $x$ such that 
\begin{equation}
\pi _{x}f^{i}(\lambda ,x,y)=0  \label{eq:dk-cond}
\end{equation}%
under the constraint 
\begin{equation}
f^{l}(\lambda ,x,y)\in D,\quad l=0,1,\dots ,i.  \label{eq:dk-constr}
\end{equation}%
From Lemma \ref{lem:unstable-disc} it follows immediately that this problem
has a unique solution $x_{i}(\lambda ,y)$ which is as smooth as $f$.

\begin{lemma}
\label{lem:wcs-graph-lim} Let $b_{i}:\Lambda \times \overline{B}%
_{s}(R)\rightarrow D$ be given by $b_{i}(\lambda ,y)=(\lambda ,x_{i}(\lambda
,y),y)$. Then $b_{i}$ is a center vertical disc and the sequence $b_{i}$
converges uniformly to $W^{cs}$. Moreover, $W^{cs}$ is a center vertical
disk in $D$, such that 
\begin{equation}
\pi _{x}W^{cs}\subset B_{u}(R).  \label{eq:pxWcs}
\end{equation}
\end{lemma}

\begin{proof}
To show that $b_{i}$ is a center vertical disc, we have to prove that if $%
\lambda _{1}\in B_{c}(\lambda _{2},R_{\Lambda })$, then $b_{i}(\lambda
_{1},y_{1})\in J_{cs}\left( b_{i}(\lambda _{2},y_{2}),L\right) $. We will
argue by the contradiction. Assume that $b_{i}(\lambda _{1},y_{1})\notin
J_{cs}\left( b_{i}(\lambda _{2},y_{2}),L\right) $, which implies $%
b_{i}(\lambda _{1},y_{1})\in J_{u}\left( b_{i}(\lambda
_{2},y_{2}),1/L\right) $. Then from Lemma \ref{lem:Ju-expansion}, applied
inductively, it follows that 
\begin{equation*}
\Vert \pi _{x}(f^{i}(b_{i}(\lambda _{1},y_{1}))-f^{i}(b_{i}(\lambda
_{2},y_{2})))\Vert \geq \xi _{u,1,P}^{i}\Vert \pi _{x}\left( b_{i}(\lambda
_{2},y_{2})-b_{i}(\lambda _{2},y_{2})\right) \Vert >0.
\end{equation*}%
This contradicts (\ref{eq:dk-cond}). This establishes that $b_{i}$ are
center vertical discs.

To prove the uniform convergence of $b_{i}$ we show the Cauchy condition for
this sequence. Let $i,j\in \mathbb{Z}_{+}$. We have $b_{i}(\lambda ,y)\in
J_{u}\left( b_{i+j}(\lambda ,y),1/L\right) $, hence by Lemma \ref%
{lem:Ju-expansion} 
\begin{equation*}
\Vert \pi _{x}(f^{i}(b_{i}(\lambda ,y))-f^{i}(b_{i+j}(\lambda ,y)))\Vert
\geq \xi _{u,1,P}^{i}\Vert \pi _{x}\left( b_{i}(\lambda ,y)-b_{i+j}(\lambda
,y)\right) \Vert .
\end{equation*}%
Since 
\begin{equation*}
\Vert \pi _{x}(f^{i}(b_{i}(\lambda ,y))-f^{i}(b_{i+j}(\lambda ,y)))\Vert
=\left\Vert \pi _{x}f^{i}(b_{i+j}(\lambda ,y))\right\Vert \leq R
\end{equation*}%
we obtain 
\begin{equation*}
\left\Vert \pi _{x}\left( b_{i}(\lambda ,y)-b_{i+j}(\lambda ,y)\right)
\right\Vert \leq \frac{R}{\xi _{u,1,P}^{i}}.
\end{equation*}%
This, since $\xi _{u,1,P}>1$, proves uniform convergence of $b_{i}$ to some
disk $b$. Observe that 
\begin{equation*}
\left\{ b(\lambda ,y):\left( \lambda ,y\right) \in \Lambda \times \overline{B%
}_{s}(R)\right\} \subset W^{cs},
\end{equation*}%
because for each $(\lambda ,y)\in \Lambda \times \overline{B}_{s}(R),$ $%
b(\lambda ,y)=\lim_{i\rightarrow \infty }b_{i}(\lambda ,y)$ and 
\begin{equation}
f^{l}(b_{i}(\lambda ,y))\in D,\quad l=0,\dots ,i.  \label{eq:wcs-gtr-approx}
\end{equation}%
Fixing $l$ in (\ref{eq:wcs-gtr-approx}) and passing to the limit with $i$,
we obtain that for all $l\in \mathbb{Z}_{+},$ $f^{l}(b(\lambda ,y))\in D$.

We now need to show that we can not have a point $z\in W^{cs}$ such that $%
z\neq b\left( \pi _{\left( \lambda ,y\right) }z\right) .$ Since $b\left( \pi
_{\left( \lambda ,y\right) }z\right) \in J_{u}\left( z,1/L\right) ,$ by
Lemma \ref{lem:Ju-expansion}, for all $i\geq 0$%
\begin{equation*}
\Vert \pi _{x}(f^{i}(b(\pi _{\left( \lambda ,y\right) }z))-f^{i}(z))\Vert
\geq \xi _{u,1,P}^{i}\Vert \pi _{x}\left( b(\pi _{\left( \lambda ,y\right)
}z)-z\right) \Vert .
\end{equation*}%
Since $\Vert \pi _{x}(f^{i}(b(\pi _{\left( \lambda ,y\right)
}z))-f^{i}(z))\Vert \leq 2R$ and since $\xi _{u,1,P}>1,$ we see that $\pi
_{x}\left( b(\pi _{\left( \lambda ,y\right) }z)-z\right) =0$, which implies
that $b(\pi _{\left( \lambda ,y\right) }z)=z$.

Condition (\ref{eq:pxWcs}) is an immediate consequence of (\ref%
{eq:homotopy-exit}), since if we had $z\in W^{cs}$ with $\left\Vert \pi
_{x}z\right\Vert =R$ then $f(z)\notin D$.

We finish by showing that $b$ is a center-vertical disc. We have already
established that $b_{i}$ are center-vertical discs. Passing to the limit,
for any $\left( \lambda ,y\right) \in \Lambda \times B_{s}(R)$, 
\begin{equation}
b\left( \overline{B}_{c}(\lambda ,R_{\Lambda })\times \overline{B}%
_{s}(R)\right) \subset \overline{J_{cs}\left( b(\lambda ,y),L\right) }.
\label{eq:cent-stab-cond-in-proof}
\end{equation}%
The condition (\ref{eq:cv-cone-cond}) follows from Corollary \ref%
{cor:unstable-lip} by the following argument. If we had a point in $\left(
\lambda ^{\ast },y^{\ast }\right) \neq (\lambda ,y)$ such that%
\begin{equation*}
b\left( \lambda ^{\ast },y^{\ast }\right) \in \partial J_{cs}\left(
b(\lambda ,y),L\right) 
\end{equation*}%
then $b\left( \lambda ^{\ast },y^{\ast }\right) \in J_{cs}^{c}\left(
b(\lambda ,y),L\right) $ and by Corollary \ref{cor:unstable-lip}, 
\begin{equation}
f\left( b\left( \lambda ^{\ast },y^{\ast }\right) \right) \in \mathrm{int}%
J_{cs}^{c}\left( f(b(\lambda ,y)),L\right) .
\label{eq:cent-stab-cond-in-proof-2}
\end{equation}%
Since $f(W^{cu})=W^{cu},$ (\ref{eq:cent-stab-cond-in-proof-2}) contradicts (%
\ref{eq:cent-stab-cond-in-proof}).
\end{proof}

\begin{lemma}
\label{lem:Cm-norm-bound-cs}Let $m\leq k.$ Let $b_{i}$ be the sequence of
center-horizontal discs defined in Lemma \ref{lem:wcs-graph-lim}. Assume
that $b_{i}$ are $C^{m}$ and that for any $i,$ $\left\Vert \pi
_{x}b_{i}\right\Vert _{C^{m}}<c_{m},$ with $c_{m}$ independent of $i$. If $f$
satisfies rate conditions of order $m$, then $\left\Vert \pi
_{x}b_{i}\right\Vert _{C^{m+1}}<c_{m+1}$ for a constant independent of $i$.
\end{lemma}

\begin{proof}
The proof goes along the same lines as the proof of Lemma \ref%
{lem:Cm-norm-bound}. We shall write $\theta =\left( \lambda ,y\right) $.
Since $b_{i}$ follows from the solution of problem (\ref{eq:dk-cond})%
\begin{equation}
f(b_{l+1}(\Lambda \times \overline{B}_{s}(R)))\subset b_{l}(\Lambda \times 
\overline{B}_{s}(R)).  \label{eq:wcs-bi-image}
\end{equation}

Let us fix $i\in \mathbb{N}$. Our aim will be to show that $\left\Vert \pi
_{x}b_{i}\right\Vert _{C^{m+1}}$ is bounded and that the bound is
independent of $i$. Let $\theta _{i}$ be any chosen point from $\Lambda
\times \overline{B}_{s}(R)$ and let $\theta _{i-1},\ldots ,\theta _{0}\in
\Lambda \times \overline{B}_{u}(R)$ be a sequence defined as%
\begin{equation*}
\theta _{l}=\pi _{\theta }f(b_{l+1}(\theta _{l+1})),
\end{equation*}%
for $l=0,\ldots ,i-1$. By (\ref{eq:wcs-bi-image}),%
\begin{equation*}
b_{l}(\theta _{l})=f(b_{l+1}(\theta _{l+1})).
\end{equation*}%
For $l=0,\ldots ,i$, let $\mathcal{P}_{m}^{l}:\mathbb{R}^{c+s}\supset
B\left( 0,\delta \right) \rightarrow \mathbb{R}^{s}$ be a polynomial of
degree $m$, defined as 
\begin{equation*}
\mathcal{P}_{m}^{l}=\pi _{x}b_{l}\left( \theta _{l}\right) +T_{\pi
_{x}b_{l},m,\theta _{l}}.
\end{equation*}%
Since $\left\Vert \pi _{y}b_{l}\right\Vert _{C^{m}}<c_{m}$ for $c_{m}$
independent from $l$, the polynomials $\mathcal{P}_{m}^{l}$ have a uniform
bounds for their coefficients, which is independent from $l$ and $i$. Since $%
b_{l}$ are center vertical discs, $\pi _{x}b_{l}$ are Lipschitz with a
constant $L,$ hence%
\begin{equation*}
\left\Vert D\mathcal{P}_{m}^{l}(0)\right\Vert \leq L.
\end{equation*}

Let us consider coordinates $\mathrm{x}=x,$ $\mathrm{y}=(\lambda ,y)$ and
constants $\xi =\xi_{u,2}$, $\mu =\mu _{cs,1}$ and $B=\left\Vert
f\right\Vert _{C^{1}}$. By (\ref{eq:wcs-bi-image}) we see that for $%
k=0,\ldots ,i-1$ 
\begin{eqnarray*}
\mathrm{graph}(T_{\pi _{x}f\circ (\mathrm{id},\mathcal{P}_{m}^{l+1}),m,0})
&=&\mathrm{graph}(T_{\pi _{x}f\circ b_{l+1},m,0}) \\
&\subset &\mathrm{graph}(T_{\pi _{x}b_{l},m,0}) \\
&=&\mathrm{graph}\left( \mathcal{P}_{m}^{l}\right) .
\end{eqnarray*}%
From Theorem \ref{thm:stb-jet-propagation}, for sufficiently large $M$ and
sufficiently small $\delta $%
\begin{equation}
f(J_{cs}^{c}(b_{l+1}(\theta _{l+1}),\mathcal{P}_{m}^{l+1},M,\delta ))\subset
J_{cs}^{c}(b_{l}(\theta _{l}),\mathcal{P}_{m}^{l},M).
\label{eq:cone-alignment-cs}
\end{equation}%
Note that the choice of $M$ and $\delta $ does not depend on $l$. Since $%
b_{0}(\lambda ,y)=\left( \lambda ,0,y\right) $ is a flat disc, we have 
\begin{equation*}
b_{0}(B_{c+s}(\theta _{0},\delta ))\subset J_{cu}(b_{0}(\theta _{0}),%
\mathcal{P}_{m}^{0}=0,M,\delta ).
\end{equation*}%
This by (\ref{eq:cone-alignment-cs}) implies that 
\begin{equation}
b_{l}(B_{c+u}(\theta _{k},\delta ))\subset J_{cu}(b_{l}(\theta _{l}),%
\mathcal{P}_{m}^{l},M,\delta ).  \label{eq:bk-cs-in-cone}
\end{equation}

From (\ref{eq:bk-cs-in-cone}), by Lemma \ref{lem:inv-func-in-stb-jet}, we
obtain a uniform bound 
\begin{equation*}
\left\Vert \frac{\partial ^{m+1}\pi _{x}b_{i}(\theta _{i})}{\partial \mathrm{%
y}_{i_{1}}\ldots \partial \mathrm{y}_{i_{m+1}}}\right\Vert \leq CM,
\end{equation*}%
where $C$ is independent of $\theta _{i}$. This means that%
\begin{equation*}
\left\Vert \pi _{x}b_{i}\right\Vert _{C^{m+1}}\leq c_{m+1},
\end{equation*}%
where $c_{m+1}$ depends on $CM$ and $c_{m}$, but is independent of $i$.
\end{proof}

\begin{lemma}
\label{lem:cs-smooth} If $f$ satisfies the assumptions from Theorem \ref%
{th:main}, then the manifold $W^{cs}$ is $C^{k}.$
\end{lemma}

\begin{proof}
The functions $\pi _{x}b_{i}$ are $C^{k+1}$ and uniformly Lipschitz with
constant $L$. The fact that $C^{k}$ smoothness is preserved as we pass to
the limit follows from Lemma \ref{lem:Cm-norm-bound-cs} and mirror arguments
to the proof of Lemma \ref{lem:cu-smooth}.
\end{proof}

We finish the section by proving Theorem \ref{th:wcs}:

\begin{proof}[Proof of Theorem \protect\ref{th:wcs}]
\label{proof:wcs}We shall write $\theta =\left( \lambda ,y\right) .$ Our
first aim is to show that for $\theta _{1}\neq \theta _{2}$%
\begin{equation}
\left\Vert \pi _{x}\left( b_{i}\left( \theta _{1}\right) -b_{i}\left( \theta
_{2}\right) \right) \right\Vert \leq M\left\Vert \theta _{1}-\theta
_{2}\right\Vert .  \label{eq:tmp-lip-wcu-cond}
\end{equation}%
Let $z_{1}=b_{i}\left( \theta _{1}\right) $, $z_{2}=b_{i}\left( \theta
_{2}\right) $ and suppose that (\ref{eq:tmp-lip-wcu-cond}) does not hold.
Then $z_{1}\in J_{u}\left( z_{2},1/M\right) .$ From Theorem \ref%
{thm:Lip-unstable-jet} (taking $\mathrm{x}=x$, $\mathrm{y}=\theta $ and $1/M$
in place of $M$) we see that for $l=0,\ldots i,$ $f^{l}(z_{1})\in
J_{u}\left( f^{l}(z_{2}),1/M\right) .$ By the same argument as the one in
the proof of Lemma \ref{lem:Ju-expansion} (on page \pageref{app:Ju-expansion}%
) we have%
\begin{equation*}
\left\Vert \pi _{x}\left( f^{i}\left( z_{1}\right) -f^{i}\left( z_{2}\right)
\right) \right\Vert \geq \xi ^{i}\left\Vert \pi _{x}\left(
z_{1}-z_{2}\right) \right\Vert >0,
\end{equation*}%
which contradicts the fact that by definition of $b_{i}$, $\pi
_{x}f^{i}\left( z_{1}\right) =\pi _{x}f^{i}\left( z_{2}\right) =0$. Thus we
have proven (\ref{eq:tmp-lip-wcu-cond}).

The claim follows by passing with $i$ to infinity in (\ref%
{eq:tmp-lip-wcu-cond}).
\end{proof}

%TCIDATA{Version=5.00.0.2606}
%TCIDATA{LaTeXparent=0,0,MMFedit.tex}

\section{Normally hyperbolic manifold\label{sec:nhim-exists}}

In this section we establish the existence of the normally hyperbolic
invariant manifold from Theorem \ref{th:main}. Throughout the section we
assume that assumptions of Theorem \ref{th:main} are satisfied.

\begin{lemma}
\label{lem:Wcu-Wcs-intersect}For any $\lambda ^{\ast }\in \Lambda $ there
exists a point $p^{\ast }\in W^{cu}\cap W^{cs}$ with $\pi _{\lambda }p^{\ast
}=\lambda ^{\ast }$.
\end{lemma}

\begin{proof}
Let $G:\overline{B}_{c}(\lambda ^{\ast },R_{\Lambda })\times \overline{B}%
_{u}(R)\times \overline{B}_{s}(R)\rightarrow \overline{B}_{c}(\lambda ^{\ast
},R_{\Lambda })\times \overline{B}_{u}(R)\times \overline{B}_{s}(R)$ be
defined as%
\begin{equation*}
G\left( \lambda ,x,y\right) =\left( \lambda ^{\ast },w^{cs}(\lambda
,y),w^{cu}(\lambda ,x)\right) .
\end{equation*}%
By the Brouwer fixed point theorem, there exists a $p^{\ast }$ such that $%
G\left( p^{\ast }\right) =p^{\ast }$. We see that $\pi _{\lambda }p^{\ast
}=\lambda ^{\ast }$. Let $x^{\ast }=\pi _{x}p^{\ast }$ and $y^{\ast }=\pi
_{y}p^{\ast }$. Since $G\left( p^{\ast }\right) =p^{\ast },$%
\begin{eqnarray*}
w^{cs}(\lambda ^{\ast },y^{\ast }) &=&x^{\ast }, \\
w^{cu}(\lambda ^{\ast },x^{\ast }) &=&y^{\ast },
\end{eqnarray*}%
hence%
\begin{equation*}
p^{\ast }=\left( \lambda ^{\ast },x^{\ast },w^{cu}(\lambda ^{\ast },x^{\ast
})\right) =\left( \lambda ^{\ast },w^{cs}(\lambda ^{\ast },y^{\ast
}),y^{\ast }\right)
\end{equation*}%
clearly lies on $W^{cu}\cap W^{cs}$.
\end{proof}

\begin{lemma}
\label{lem:Wcu-Wcs-intersect-transv}Let $p\in W^{cu}\cap W^{cs}$, then $%
W^{cu}$ and $W^{cs}$ intersect transversally at $p$.
\end{lemma}

\begin{proof}
The manifold $W^{cu}$ is parameterized by $\phi _{cu}:\left( \lambda
,x\right) \rightarrow \left( \lambda ,x,w^{cu}(\lambda ,x)\right) $ and $%
W^{cs}$ is parameterized by $\phi _{cs}:( \lambda ,y) \rightarrow ( \lambda
,w^{cs}(\lambda ,y),y)$. Let 
\begin{equation*}
V=\mathrm{span}\{D\phi _{cu}(p)v+D\phi _{cs}(p)w:v\in \mathbb{R}^{c}\times 
\mathbb{R}^{u},w\in \mathbb{R}^{c}\times \mathbb{R}^{s}\}.
\end{equation*}%
We need to show that%
\begin{equation}
V=\mathbb{R}^{c}\times \mathbb{R}^{u}\times \mathbb{R}^{s}.
\label{eq:V-span}
\end{equation}%
We see that $V$ is equal to the range of the $\left( c+u+s\right) \times
\left( c+c+u+s\right) $ matrix%
\begin{equation*}
A=\left( 
\begin{array}{llll}
\mathrm{id} & \mathrm{id} & 0 & 0 \\ 
0 & \frac{\partial w_{cs}}{\partial \lambda } & \mathrm{id} & \frac{\partial
w_{cs}}{\partial y} \\ 
\frac{\partial w_{cu}}{\partial \lambda } & 0 & \frac{\partial w_{cu}}{%
\partial x} & \mathrm{id}%
\end{array}%
\right) .
\end{equation*}%
We will show that 
\begin{equation*}
B=\left( 
\begin{array}{ll}
\mathrm{id} & \frac{\partial w_{cs}}{\partial y} \\ 
\frac{\partial w_{cu}}{\partial x} & \mathrm{id}%
\end{array}%
\right)
\end{equation*}%
is invertible. If for $q=\left( x,y\right) $, $Bq=0$ then%
\begin{equation*}
x-\frac{\partial w_{cs}}{\partial y}\frac{\partial w_{cu}}{\partial x}x=0,
\end{equation*}%
which since $w_{cs}$ and $w_{cu}$ are Lipschitz with constant $L<1$ implies
that $x=0,$ and in turn that $y=0$. Since $B$ is invertible, it is evident
that the rank of $A$ is $c+u+s,$ which implies (\ref{eq:V-span}).
\end{proof}

\begin{lemma}
The $\Lambda ^{\ast }=W^{cu}\cap W^{cs}$ is a $C^{k}$ manifold, which is a
graph over $\Lambda $ of a function $\chi :\Lambda \rightarrow \overline{B}%
_{u}\times \overline{B}_{s}$, which is Lipschitz with a constant $\frac{%
\sqrt{2}L}{\sqrt{1-L^{2}}}$.
\end{lemma}

\begin{proof}
From Lemma \ref{lem:Wcu-Wcs-intersect} it follows that for every $\lambda
\in \Lambda$ the set $W^{cu} \cap W^{cs} \cap \{p \in D \ | \
\pi_\lambda=\lambda \}$ is nonempty. We will show that this set consists
from one point $\chi(\lambda)$ and $\chi$ is a function satisfying the
Lipschitz condition.

Assume that $p_{1}=( \lambda _{1},x_1,y_1) \in W^{cs} \cap W^{cu}$ and $%
p_{2}=( \lambda _{2},x_2,y_2)) \in W^{cs} \cap W^{cu}$. Moreover, we assume
that $\lambda_1 \in B_c(\lambda_2,R_\Lambda)$ (they are in the same chart).
Therefore we know that 
\begin{eqnarray}
p_{1} &\in &J_{cu}\left( p_{2},L\right) ,  \label{eq:p1-Jcu} \\
p_{1} &\in &J_{cs}\left( p_{2},L\right) .  \label{eq:p1-Jcs}
\end{eqnarray}%
Let $( \lambda ,x,y) =p_{1}-p_{2}=( \lambda _{1}-\lambda_{2},x_1 -
x_2,y_1-y_2) $. By (\ref{eq:p1-Jcu}--\ref{eq:p1-Jcs}) we obtain 
\begin{equation*}
\left\Vert y\right\Vert \leq L\left\Vert \left( \lambda ,x\right)
\right\Vert ,\qquad \qquad \left\Vert x\right\Vert \leq L\left\Vert \left(
\lambda ,y\right) \right\Vert ,
\end{equation*}%
hence 
\begin{align*}
\left\Vert y\right\Vert ^{2}& \leq L^{2}\left( \left\Vert \lambda
\right\Vert ^{2}+\left\Vert x\right\Vert ^{2}\right) , \\
\left\Vert x\right\Vert ^{2}& \leq L^{2}\left( \left\Vert \lambda
\right\Vert ^{2}+\left\Vert y\right\Vert ^{2}\right) .
\end{align*}%
From above%
\begin{equation*}
\left( \left\Vert x\right\Vert ^{2}+\left\Vert y\right\Vert ^{2}\right)
\left( 1-L^{2}\right) \leq 2L^{2}\left\Vert \lambda \right\Vert ^{2},
\end{equation*}%
which gives 
\begin{equation*}
\left\Vert \left( x,y\right) \right\Vert \leq \frac{\sqrt{2}L}{\sqrt{1-L^{2}}%
}\left\Vert \lambda \right\Vert .
\end{equation*}
Observe that this implies that if $\lambda_1=\lambda_2$, then $p_1=p_2$.
This establishes the uniqueness of the intersection of $W^{cu} \cap W^{cs}
\cap \{p \in D \ | \ \pi_\lambda=\lambda \}$, therefore $\chi(\lambda)$ is
well defined.

From the above computations it follows that 
\begin{equation*}
\left\Vert \chi (\lambda _{1})-\chi (\lambda _{2})\right\Vert \leq \frac{%
\sqrt{2}L}{\sqrt{1-L^{2}}}\left\Vert \lambda _{1}-\lambda _{2}\right\Vert.
\end{equation*}

The fact that $\Lambda ^{\ast }$ is a graph of a $C^{k}$ function $\chi
:\Lambda \rightarrow \overline{B}_{u}\times \overline{B}_{s}$ follows from (%
\ref{lem:Wcu-Wcs-intersect-transv}) and the fact that $W^{cu}$ and $W^{cs}$
are $C^{k}$.
\end{proof}

%TCIDATA{Version=5.00.0.2606}
%TCIDATA{LaTeXparent=0,0,MMFedit.tex}

\section{Unstable fibers\label{sec:unstb-fibers}}

The goal of this section is to establish the existence of the foliation of $%
W^{cu}$ into unstable fibers $W_{z}^{u}$ for $z\in W^{cu}$. In this section $%
\theta=(\lambda,y)$. Throughout this section we assume that assumptions of
Theorem \ref{th:main} hold.

For $z\in D$ we define $b_{z}$ as a horizontal disk in $D$ by $%
b_{z}(x)=(\pi_{\lambda}z,x,\pi_{y}z)$.

By Lemma \ref{lem:fWu-injective} we know that $f|_{W^{cu}}$ is injective,
hence for any $z\in W^{cu}$ the backward trajectory is unique and equal to $%
\left\{ \left( f|_{W^{cu}}\right) ^{i}\right\} _{i=-\infty }^{0}$. To
simplify notations we shall denote such backward trajectory by $z_{i}=\left(
f|_{W^{cu}}\right) ^{i}$, for $i=0,-1,\ldots $.

For $z\in W^{cu}$ consider a sequence of horizontal disks in $D$, 
\begin{equation}
d_{n,z}=\mathcal{G}_{h}^{n}(b_{z_{-n}}),\quad n=1,2,\ldots 
\label{eq:wuq-dkq}
\end{equation}%
where $G_{h}$ is the graph transform defined just after Lemma~\ref%
{lem:unstable-disc}. Our aim will be to show that $d_{n,z}$ converge to $%
W_{z}^{u}$, as defined by Definition \ref{def:Wu-fiber}. We start with a
technical lemma.

\begin{lemma}
\label{lem:unstb-fib-conv} Assume that $f^{j}(z)\in D$ for $j=0,1,\dots ,n$.
If $f^{j}(q_{i})\in J_{u}\left( f^{j}(z),1/L\right) \cap D$, for $i=1,2$ and 
$j=0,1,\dots ,n$ and 
\begin{equation*}
f^{n}(q_{1})\notin J_{u}(f^{n}(q_{2}),1/L)
\end{equation*}%
then for $j=0,1,\dots,n$ holds 
\begin{align*}
\Vert \pi _{\theta }(f^{j}(q_{1})-f^{j}(q_{2}))\Vert & \leq \frac{4R}{L}%
\left( \frac{\mu _{cs,1}}{\xi _{u,1,P}}\right)^{j} \frac{1}{(\xi
_{u,1,P})^{n-j}}, \\
\Vert f^{j}(q_{1})-f^{j}(q_{2})\Vert & \leq \left( 1+L\right) \frac{4R}{L}%
\left( \frac{\mu _{cs,1}}{\xi _{u,1,P}}\right)^{j} \frac{1}{(\xi
_{u,1,P})^{n-j}}.
\end{align*}
\end{lemma}

\begin{proof}
Our assumption $f^{j}(q_{i})\in J_{u}\left( f^{j}(z),1/L\right) $, $i=1,2$
and $j=0,1,\dots ,n$ implies that $f^{j}(q_{1}),f^{j}(q_{2}),f^j(z)$ are
contained in the same charts for $j=0,1,\dots ,n$.

By Corollary \ref{cor:unstable-lip}, $f^{j}(q_{1})\notin
J_{u}(f^{j}(q_{2}),1/L)$ for $j=0,1,\dots ,n$, hence 
\begin{equation*}
f^{j}(q_{1})\in J_{cs}(f^{j}(q_{2}),L).
\end{equation*}
By Lemma \ref{lem:cs-cone-mu-cs} this implies for $j=1,2,\dots,n$ 
\begin{equation}
\Vert \pi _{\theta }(f^{j}(q_{1})-f^{j}(q_{2}))\Vert \leq \mu _{cs,1}\Vert
\pi _{\theta }(f^{j-1}(q_{1})-f^{j-1}(q_{2}))\Vert \leq \mu _{cs,1}^{j}\Vert
\pi _{\theta }(q_{1}-q_{2})\Vert .  \label{eq:q+-iter-estm}
\end{equation}%
We estimate $\Vert \pi _{\theta }(q_{1}-q_{2})\Vert $ using the expansion in
the $x$-direction. By Lemma \ref{lem:Ju-expansion} we have for $i=1,2$ 
\begin{equation*}
2R\geq \Vert \pi _{x}(f^{n}(q_{i})-f^{n}(z))\Vert \geq \xi _{u,1,P}\Vert \pi
_{x}(f^{n-1}(q_{i})-f^{n-1}(z))\Vert \geq \xi _{u,1,P}^{n}\Vert \pi
_{x}(q_{i}-z)\Vert ,
\end{equation*}%
hence we obtain 
\begin{equation*}
\Vert \pi _{x}(q_{i}-z)\Vert \leq \frac{2R}{\xi _{u,1,P}^{n}}.
\end{equation*}%
Since $q_{i}\in J_{u}(z,1/L)$ we get 
\begin{equation*}
\Vert \pi _{\theta }(q_{i}-z)\Vert \leq \frac{1}{L}\Vert \pi
_{x}(q_{i}-z)\Vert \leq \frac{2R}{L\xi _{u,1,P}^{n}}.
\end{equation*}%
From the triangle inequality we obtain 
\begin{equation*}
\Vert \pi _{\theta }(q_{1}-q_{2})\Vert \leq \Vert \pi _{\theta
}(q_{1}-z)\Vert +\Vert \pi _{\theta }(q_{2}-z)\Vert \leq \frac{4R}{L\xi
_{u,1,P}^{n}}.
\end{equation*}%
By combining the above inequality with (\ref{eq:q+-iter-estm}) we obtain 
\begin{equation*}
\Vert \pi _{\theta }(f^{j}(q_{1})-f^{j}(q_{2}))\Vert \leq \frac{4R}{L}\left( 
\frac{\mu _{cs,1}}{\xi _{u,1,P}}\right) ^{j} \frac{1}{\xi _{u,1,P}^{n-j}}.
\end{equation*}

Since $f^{n}(q_{1})\notin J_{u}(f^{n}(q_{2}),1/L)$,%
\begin{equation*}
\left\Vert \pi _{\theta }\left( f^{n}(q_{1})-f^{n}(q_{2})\right) \right\Vert
>\frac{1}{L}\left\Vert \pi _{x}\left( f^{n}(q_{1})-f^{n}(q_{2})\right)
\right\Vert ,
\end{equation*}%
hence 
\begin{align*}
\Vert f^{j}(q_{1})-f^{j}(q_{2})\Vert & \leq \left\Vert \pi _{\theta }\left(
f^{j}(q_{1})-f^{j}(q_{2})\right) \right\Vert +\left\Vert \pi _{x}\left(
f^{j}(q_{1})-f^{j}(q_{2})\right) \right\Vert \\
& \leq \left( 1+L\right) \frac{4R}{L} \left( \frac{\mu _{cs,1}}{\xi _{u,1,P}}%
\right)^{j} \frac{1}{(\xi _{u,1,P})^{n-j}}
\end{align*}%
as required.
\end{proof}

\begin{lemma}
\label{lem:wuq-conv}Assume that $z\in W^{cu}$. For any $n\geq 0,$ $d_{n,z}$
is a horizontal disc and the sequence $d_{n,z}$ converges uniformly to a
horizontal disk $d_z$. Moreover, 
\begin{equation*}
W_{z}^{u}=\left\{ \left( \pi _{\lambda }w_{z}^{u}\left( x\right) ,x,\pi
_{y}w_{z}^{u}\left( x\right) \right) :x\in \overline{B}_{u}(R)\right\} ,
\end{equation*}
where $w_{z}^{u}:\overline{B}_{u}(R)\rightarrow \Lambda \times \overline{B}%
_{s}(R)$ and $d_z(x)=\left( \pi _{\lambda }w_{z}^{u}\left( x\right) ,x,\pi
_{y}w_{z}^{u}\left( x\right) \right)$.
\end{lemma}

\begin{proof}
We show first the uniform convergence. Let us fix $x\in \overline{B}_{u}(R)$%
. Our goal is to estimate $\|d_{n+j,z}(x) - d_{n,z}(x)\|$. Observe first
that from the definition of the graph transform $\mathcal{G}_{h}$, it
follows that for each $n\in \mathbb{Z}_{+}$ and for each $x\in \overline{B}%
_{u}(R)$ the point $d_{n,z}(x)$ has a backward orbit $\left\{ p_{i}\right\}
_{i=-n}^{0}$ of length $n+1,$ 
\begin{align*}
p_{0}& =d_{n,z}(x),\quad f(p_{i})=p_{i+1}\quad \text{for }i=-n,-n+1,\dots ,-1
\\
p_{i}& \in J_{u}\left( z_{i},1/L\right) ,\quad p_{i}\in d_{n+i,z_{i}}\left( 
\overline{B}_{u}(R)\right) \quad \text{for }i=-n,-n+1,\dots ,-1,0.
\end{align*}

Let $n,j$ be positive integers. From the above observation we can find
(define) $q_{1}$ and $q_{2}$ as follows. Let $q_{1}$ be such that $%
f^{i}(q_{1})\in J_{u}\left( z_{-n+i},1/L\right) $ for $i=0,1,\dots ,n$ and $%
f^{n}(q_{1})=d_{n,z}(x)$, analogously, let $q_{2}$ be such that $%
f^{i}(q_{2})\in J_{u}\left( z_{-n+i},1/L\right) $ for $i=0,1,\dots ,n$ and $%
f^{n}(q_{2})=d_{n+j,z}(x)$.

Observe that since 
\begin{equation*}
\pi _{x}\left( f^{n}(q_{2})-f^{n}(q_{1})\right) =\pi _{x}\left(
d_{n,z}(x)-d_{n+j,z}(x)\right) =0,
\end{equation*}%
we have $\Vert \pi _{\theta }(f^{n}(q_{2})-f^{n}(q_{1}))\Vert =\Vert
f^{n}(q_{2})-f^{n}(q_{1})\Vert $. Assume that $f^{n}(q_{2})\neq f^{n}(q_{1})$%
, then from Lemma~\ref{lem:unstb-fib-conv} applied to $q_{1}$, $q_{2}$ and $%
z_{-n}$ it follows that 
\begin{equation*}
\Vert d_{n,z}(x)-d_{n+j,z}(x)\Vert =\Vert \pi _{\theta
}(d_{n,z}(x)-d_{n+j,z}(x))\Vert \leq \frac{4R}{L}\left( \frac{\mu _{cs,1}}{%
\xi _{u,1,P}}\right) ^{n}.
\end{equation*}%
Since by our assumptions $\frac{\mu _{cs,1}}{\xi _{u,1,P}}<1$, we see that $%
d_{n,z}$ is a Cauchy sequence. Let us denote the limit by $d_{z}$.

Since $b_{z_{-n}}$ is a horizontal disc, so by Lemma~\ref{lem:unstable-disc}
is $d_{n,z}=\mathcal{G}_{h}^{n}\left( b_{z_{-n}}\right) .$ The properties (%
\ref{eq:hor-disc-prop}) are preserved when passing to the limit, hence $%
d_{z} $ is a horizontal disc.

We show that for all $x\in \overline{B}_{u}(R),$ $d_{z}(x)\in W^{cu}$. For
this we need to construct a full backward orbit through $d_{z}(x)$. Let us
consider backward orbits through $d_{n,z}(x)$ of length $n+1$. From Lemma~%
\ref{lem:unstb-fib-conv} it follows that they converge to full backward
orbit through $d_{z}(x)$. Therefore $d_{z}(x)\in W^{cu}$ for $x\in \overline{%
B}_{u}(R)$. From this reasoning it follows also that for $i<0$ 
\begin{equation}
\mathcal{G}_{h}(d_{z_{i}})=d_{z_{i+1}}.  \label{eq:finv-dq}
\end{equation}

We will now show that $\{d_{z}(x)\ |\ x\in \overline{B}_{u}(R)\}\subset
W_{z}^{u}$. For any $x\in \overline{B}_{u}(R)$ and backward trajectory $%
\left\{ p_{i}\right\} _{i=-\infty }^{0}$ of $d_{z}(x)$, from (\ref%
{eq:finv-dq}) it follows that $p_i \in d_{z_i}$ for $i\leq 0$. Since $%
d_{z_{i}}$ are horizontal discs we infer that $p_i \in J_{u}\left(
z_{i},1/L\right)$, as required.

To show that $W_{z}^{u}\subset \{d_{z}(x)\ |\ x\in \overline{B}_{u}(R)\},$
let us consider $p\in W^{cu}$, with a backward trajectory (note that by
Lemma \ref{lem:fWu-injective} such trajectory is unique) $\left\{
p_{i}\right\} _{i=-\infty }^{0},$ $p_{i}=\left( f|_{W^{cu}}\right) ^{i}(p),$
such that $p_{i}\in J_{u}\left( z_{i},1/L\right) \ $for all $i<0$. Let $%
x=\pi _{x}p$. We will show that $p=d_{z}(x)$. From Lemma~\ref%
{lem:unstb-fib-conv} it follows that 
\begin{equation*}
\Vert p-d_{z}(x)\Vert =\Vert \pi _{\theta }(p-d_{z}(x))\Vert
=\lim_{n\rightarrow \infty }\left\Vert \pi _{\theta
}(p-d_{n,z}(x))\right\Vert \leq \lim_{n\rightarrow \infty }\frac{4R}{L}%
\left( \frac{\mu _{cs,1}}{\xi _{u,1,P}}\right) ^{n}=0.
\end{equation*}%
Therefore $p=d_{z}(x)$.

The function $w_{z}^{u}$ can be defined as $w_{z}^{u}(x)=\pi _{\theta
}d_{z}(x).$
\end{proof}

\begin{lemma}
\label{lem:Cm-norm-bound-fib-u}Let $m\leq k.$ Let $d_{n,z}$ be the sequence
of horizontal discs defined as $d_{n,z}=\mathcal{G}_{h}^{i}(b_{z_{-n}}).$
Assume that $d_{n,z}$ are $C^{m}$ and that for any $i, $ $\left\Vert \pi
_{\theta }d_{n,z}\right\Vert _{C^{m}}<c_{m},$ with $c_{m}$ independent of $n$%
. If $f$ satisfies rate conditions of order $m$, then $\left\Vert \pi
_{\theta }d_{n,z}\right\Vert _{C^{m+1}}<c_{m+1}$ for a constant independent
of $n$.
\end{lemma}

\begin{proof}
The proof follows from identical arguments to the proof of Lemma \ref%
{lem:Cm-norm-bound}. The only difference is that when we apply Theorem \ref%
{thm:unstb-jet-propagation}, we choose coordinates $\mathrm{x}=x,$ $\mathrm{y%
}=\theta =\left( \lambda ,y\right) $ and constants $\xi =\xi_{1,u}$, $\mu
=\mu _{cs,2}$. Note that conditions (\ref{eq:rate-cond-0}), (\ref%
{eq:rate-cond-3}) imply (\ref{eq:full-unstb-rate-cond}) for any $m\geq 0$.
\end{proof}

\begin{lemma}
\label{lem:Wuz-Ck}For any $z\in W^{cu}$ the manifold $W_{z}^{u}$ is $C^{k}$.
\end{lemma}

\begin{proof}
The functions $\pi _{\theta }d_{n,z}$ are $C^{k+1}$ and uniformly Lipschitz
with constant $1/L$. The fact that $C^{k}$ smoothness is preserved as we
pass to the limit follows from Lemma \ref{lem:Cm-norm-bound-fib-u} and
mirror arguments to the proof of Lemma \ref{lem:cu-smooth}.
\end{proof}

\begin{lemma}
\label{lem:unstb-fiber-C-conv} For any $z\in W^{cu}$. If $p\in W_{z}^{u}$,
then for $n\geq 0$ $f_{|W^{cu}}^{-n}(p),f|_{W^{cu}}^{-n}(z)$ are in the same
chart and 
\begin{equation*}
\left\Vert \left( f_{|W^{cu}}\right)^{-n}(p)-\left(f|_{W^{cu}}\right)
^{-n}(z)\right\Vert \leq \left( 1+\frac{1}{L}\right) \Vert \pi
_{x}(p-z)\Vert \xi _{u,1,P}^{-n},\qquad n\geq 0.
\end{equation*}
\end{lemma}

\begin{proof}
We first observe that for any $q_{1},q_{2}\in D,$ such that $q_{1}\in
J_{u}\left( q_{2},1/L\right) ,$ holds%
\begin{equation*}
\left\Vert \pi _{\theta }\left( q_{1}-q_{2}\right) \right\Vert \leq \frac{1}{%
L}\left\Vert \pi _{x}\left( q_{1}-q_{2}\right) \right\Vert ,
\end{equation*}%
\begin{equation}
\left\Vert q_{1}-q_{2}\right\Vert \leq \left\Vert \pi _{\theta }\left(
q_{1}-q_{2}\right) \right\Vert +\left\Vert \pi _{x}\left( q_{1}-q_{2}\right)
\right\Vert \leq \left( 1+\frac{1}{L}\right) \left\Vert \pi _{x}\left(
q_{1}-q_{2}\right) \right\Vert .  \label{eq:tmp-Wuz-conv-1}
\end{equation}

Let $p\in W_{z}^{u}$. By Lemma \ref{lem:fWu-injective}, backward
trajectories of $p$ and $z$ are unique and equal to $p_{-n}=\left(
f|_{W^{cu}}\right) ^{-n}(p),$ $z_{-n}=\left( f|_{W^{cu}}\right) ^{-n}(z)$,
for $n\geq 0.$ Note that by the definition of $W_{z}^{u}$, for $i=0,\ldots
,n $, $p_{-i}\in J_{u}\left( z_{-i},1/L\right) .$ Let us fix $n\geq 0$. From
Lemma~\ref{lem:Ju-expansion} and (\ref{eq:tmp-Wuz-conv-1}), it follows that%
\begin{eqnarray*}
\Vert \pi _{x}(p-z)\Vert &=&\left\Vert \pi _{x}\left(
f^{n}(p_{-n})-f^{n}(z_{-n})\right) \right\Vert \\
&\geq &\xi _{u,1,P}^{n}\left\Vert \pi _{x}\left( p_{-n}-z_{-n}\right)
\right\Vert \\
&\geq &\xi _{u,1,P}^{n}\left( 1+\frac{1}{L}\right) ^{-1}\left\Vert
p_{-n}-z_{-n}\right\Vert .
\end{eqnarray*}%
This proves that for any point in $W_{z}^{u}$ holds 
\begin{equation*}
\left\Vert p_{-n}-z_{-n}\right\Vert \leq \xi _{u,1,P}^{-n}C,
\end{equation*}%
for $C=\left( 1+\frac{1}{L}\right) \Vert \pi _{x}(p-z)\Vert ,$ as required.
\end{proof}

\begin{lemma}
\label{lem:Wuz-convergence}For $z\in W^{cu}$ we define a set $U=U(z)$ as%
\begin{eqnarray*}
U &=&\{p\in D:\exists \text{ backward trajectory }\{p_{i}\}_{i=-\infty }^{0}%
\text{ of $p\in D$, and} \\
&&\text{for any such }\{p_{i}\},\ \exists C>0\text{ (which may depend on $p$%
),}\ \exists n_{0}\geq 0 \\
&&\,\text{s.t. for $n\geq n_{0},$ $p_{-n}$ and $\left( f|_{W^{cu}}\right)
^{-n}(z)$ are in the same good chart } \\
&&\,\text{and }\left\Vert p_{-n}-\left( f|_{W^{cu}}\right)
^{-n}(z)\right\Vert \leq C\xi _{u,1,P}^{-n}\}.
\end{eqnarray*}%
Then%
\begin{equation*}
W_{z}^{u}=U.
\end{equation*}
\end{lemma}

\begin{proof}
Observe that from Lemma~\ref{lem:unstb-fiber-C-conv} we obtain $%
W_{z}^{u}\subset U$. We will show that $U\subset W_{z}^{u}$ by contradiction.

Let $p\in U\setminus W_{z}^{u}$. Obviously $p\in W^{cu}$, hence for $i\leq 0$%
, by Lemma \ref{lem:fWu-injective}, its backward trajectory is uniquely
defined. Let $i^{\ast }=-n_{0}\leq 0$. Then for $n\leq i^{\ast }$ points $%
p_{-n},z_{-n}$ lie in the same chart and 
\begin{equation}
\left\Vert p_{-n}-z_{-n}\right\Vert \leq C\xi _{u,1,P}^{-n}.
\label{eq:u-xi-contr-1}
\end{equation}

Since $p\notin W_{z}^{u}$ then there exists $j_{0}\leq 0$ such that $%
p_{j_{0}}\notin J_{u}\left( z_{j_{0}},1/L\right) $. From the forward
invariance of $J_{u}$'s (see Cor.~\ref{cor:unstable-lip}) it follows that $%
p_{j}\notin J_{u}\left( z_{j},1/L\right) $ for $j\leq j_{0}$. Hence we can
find $i^{\ast \ast }\leq i^{\ast }<0$ such that $p_{i^{\ast \ast }}\notin
J_{u}\left( z_{i^{\ast \ast }},1/L\right) $. For any $i\leq i^{\ast \ast }$
holds 
\begin{equation*}
p_{i}\in J_{cs}(z_{i},L),\quad \mbox{for $i\leq i^{\ast \ast }$}
\end{equation*}%
For $n>\left\vert i^{\ast \ast }\right\vert $ from Lemma~\ref%
{lem:cs-cone-mu-cs} 
\begin{equation*}
\left\Vert \pi _{\theta }\left( f(p_{-n})-f(z_{-n})\right) \right\Vert \leq
\mu _{cs,1}\left\Vert \pi _{\theta }\left( p_{-n}-z_{-n}\right) \right\Vert ,
\end{equation*}%
hence, by the same argument, 
\begin{eqnarray*}
\left\Vert \pi _{\theta }\left( p_{i^{\ast \ast }}-z_{i^{\ast \ast }}\right)
\right\Vert &=&\left\Vert \pi _{\theta }\left( f^{n+i^{\ast \ast
}}(p_{-n})-f^{n+i^{\ast \ast }}(z_{-n})\right) \right\Vert \\
&\leq &\mu _{cs,1}^{n+i^{\ast \ast }}\left\Vert \pi _{\theta }\left(
p_{-n}-z_{-n}\right) \right\Vert \\
&=&\mu _{cs,1}^{n}\mu _{cs,1}^{i^{\ast \ast }}\left\Vert \pi _{\theta
}\left( p_{-n}-z_{-n}\right) \right\Vert ,
\end{eqnarray*}%
and in turn%
\begin{eqnarray}
\left\Vert p_{-n}-z_{-n}\right\Vert &\geq &\left\Vert \pi _{\theta }\left(
p_{-n}-z_{-n}\right) \right\Vert  \label{eq:u-xi-contr-2} \\
&\geq &\mu _{cs,1}^{-n}\left( \mu _{cs,1}^{-i^{\ast \ast }}\left\Vert \pi
_{\theta }\left( p_{i^{\ast \ast }}-z_{i^{\ast \ast }}\right) \right\Vert
\right) .  \notag
\end{eqnarray}%
Since $\xi _{u,1,P}>\mu _{cs,1}$, conditions (\ref{eq:u-xi-contr-1}) and (%
\ref{eq:u-xi-contr-2}) contradict each other. This means that $p\in
W_{z}^{u} $, as required.
\end{proof}

%\begin{remark}
%Let us note that in Lemma \ref{lem:Wuz-convergence} we have established that
%from our assumptions it follows that the constant $C$ from (\ref{eq:Wuz-in-main-thm}) does not depend on the choice of $z$ and $p$.
%\textbf{raczej, ze moze byc wzieta niezaleznie od $z$ i $p$.}
%\end{remark}

Lemmas \ref{lem:wuq-conv}, \ref{lem:Wuz-Ck}, \ref{lem:Wuz-convergence}
combined prove the claims about $W_{z}^{u}$ from Theorem \ref{th:main}. We
now prove Theorem \ref{th:wuz}, which can be used to obtain tighter
Lipschitz bounds on $w_{z}^{u}$.

\begin{proof}[Proof of Theorem \protect\ref{th:wuz}]
From Theorem \ref{thm:Lip-unstable-jet}, taking coordinates $\mathrm{x}=x$, $%
y=\theta $, for $q\in D,$ since $\frac{\xi }{\mu }>1$, 
\begin{equation}
f\left( J_{u}\left( q,M\right) \cap D\right) \subset J_{u}\left( f\left(
q\right) ,M\right) .  \label{eq:Jets-u-M-prop}
\end{equation}

By definition of $d_{0,z}$, it follows that $d_{0,z}(x_{1})\in f\left(
J_{u}\left( d_{0,z}(x_{2}),M\right) \cap D\right) ,$ for any $x_{1},x_{2}\in 
\overline{B}_{u}.$ By (\ref{eq:Jets-u-M-prop}), since $d_{n,z}=\mathcal{G}%
\left( d_{n-1,z}\right) $, we see that 
\begin{equation*}
d_{n,z}(x_{1})\in f\left( J_{u}\left( d_{n,z}(x_{2}),M\right) \cap D\right) ,
\end{equation*}%
for any $x_{1},x_{2}\in \overline{B}_{u}.$ Hence%
\begin{equation*}
\left\Vert \pi _{\theta }\left( d_{n,z}(x_{1})-d_{n,z}(x_{2})\right)
\right\Vert \leq M\left\Vert \pi _{x}\left(
d_{n,z}(x_{1})-d_{n,z}(x_{2})\right) \right\Vert =M\left\Vert
x_{1}-x_{2}\right\Vert ,
\end{equation*}%
and passing with $n$ to infinity gives%
\begin{equation*}
\left\Vert w_{z}^{u}(x_{1})-w_{z}^{u}(x_{2})\right\Vert \leq M\left\Vert
x_{1}-x_{2}\right\Vert ,
\end{equation*}%
as required.
\end{proof}

\begin{proposition}
\label{prop:fiber-inter}Let $z\in W^{cu}$. Then the intersection $%
W_{z}^{u}\cap W^{cs}$ consists of a single point and is transversal. Also
the intersection $W_{z}^{u}\cap \Lambda^*$ consists of a single point.
\end{proposition}

\begin{proof}
The proof follows from similar arguments to the proofs of Lemma \ref%
{lem:Wcu-Wcs-intersect} and Theorem \ref{lem:Wcu-Wcs-intersect-transv}.

First we show that $W_{z}^{u}$ and $W^{cs}$ intersect. By Remark \ref%
{rem:cone-in-chart}, for any point $q\in W_{z}^{u}$ we have%
\begin{equation*}
W_{z}^{u}\subset D_{\pi _{\lambda }q}=\overline{B}_{c}\left( \pi _{\lambda
}q,R_{\Lambda }\right) \times \overline{B}_{u}(R)\times \overline{B}_{s}(R).
\end{equation*}
Let us define the following function $G:D_{\pi _{\lambda }q}\rightarrow
D_{\pi _{\lambda }q},$%
\begin{equation*}
G\left( \lambda ,x,y\right) =\left( \pi _{\lambda }w_{z}^{u}(x),w^{cs}\left(
\lambda ,y\right) ,\pi _{y}w_{z}^{u}(x)\right) .
\end{equation*}%
By the Brouwer theorem we know that there exists a $q^{\ast }=\left( \lambda
^{\ast },x^{\ast },y^{\ast }\right) $ for which $q^{\ast }=G(q^{\ast })$.
This means that%
\begin{equation*}
W_{z}^{u}\ni \left( \pi _{\lambda }w_{z}^{u}(x^{\ast }),x^{\ast },\pi
_{y}w_{z}^{u}(x^{\ast })\right) =q^{\ast }=\left( \lambda ^{\ast
},w^{cs}\left( \lambda ^{\ast },y^{\ast }\right) ,y^{\ast }\right) \in
W^{cs},
\end{equation*}%
hence $q^{\ast }$ $W_{z}^{u}\cap W^{cs}.$

Now we show that the intersection point is unique. Let $q_{1},q_{2}\in
W_{z}^{u}\cap W^{cs}$. Since $W_{z}^{u}$ is a vertical disc, 
\begin{equation}
\left\Vert \pi _{\left( \lambda ,y\right) }\left( q_{1}-q_{2}\right)
\right\Vert \leq 1/L\left\Vert \pi _{x}\left( q_{1}-q_{2}\right) \right\Vert
.  \label{eq:cone-inter-1}
\end{equation}%
Since $W^{cs}$ is a center-vertical disc, if $q_{1}\neq q_{2}$ then%
\begin{equation*}
\left\Vert \pi _{x}\left( q_{1}-q_{2}\right) \right\Vert <L\left\Vert \pi
_{\left( \lambda ,y\right) }\left( q_{1}-q_{2}\right) \right\Vert ,
\end{equation*}%
a contradiction with (\ref{eq:cone-inter-1}), hence $q_{1}=q_{2}.$

Now we prove the transversality of the intersection. This is a similar
argument to the proof of Theorem \ref{lem:Wcu-Wcs-intersect-transv}.  We
first note that since $W_{z}^{u}$ is a center-vertical disc, $w_{z}^{u}$ is
Lipschitz with a constant $\rho <1/L$. The manifold $W_{z}^{u}$ is
parameterized by $\phi _{z,u}:x\rightarrow \left( \pi _{\lambda
}w_{z}^{u}(x),x,\pi _{y}w_{z}^{u}(x)\right) $ and $W^{cs}$ is parameterized
by $\phi _{cs}:(\lambda ,y)\rightarrow (\lambda ,w^{cs}(\lambda ,y),y)$. Let 
\begin{equation*}
V=\mathrm{span}\{D\phi _{z,u}(x^{\ast })v+D\phi _{cs}(\lambda ^{\ast
},y^{\ast })w:v\in \mathbb{R}^{u},w\in \mathbb{R}^{c}\times \mathbb{R}^{s}\}.
\end{equation*}%
We need to show that%
\begin{equation}
V=\mathbb{R}^{c}\times \mathbb{R}^{u}\times \mathbb{R}^{s}.
\label{eq:V-span-2}
\end{equation}%
We see that $V$ is equal to the range of the $\left( c+u+s\right) \times
\left( c+u+s\right) $ matrix%
\begin{equation*}
A=\left( 
\begin{array}{lll}
\frac{\partial \pi _{\lambda }w_{z}^{u}}{\partial x} & \mathrm{id} & 0 \\ 
\mathrm{id} & \frac{\partial w^{cs}}{\partial \lambda } & \frac{\partial
w^{cs}}{\partial y} \\ 
\frac{\partial \pi _{y}w_{z}^{u}}{\partial x} & 0 & \mathrm{id}%
\end{array}%
\right) .
\end{equation*}%
We will show that%
\begin{equation*}
B=\left( 
\begin{array}{ll}
\frac{\partial \pi _{\lambda }w_{z}^{u}}{\partial x} & \mathrm{id} \\ 
\mathrm{id} & \frac{\partial w^{cs}}{\partial \lambda }%
\end{array}%
\right)
\end{equation*}
is invertible. If for $p=\left( \lambda ,x\right) $, $Bp=0$ then%
\begin{equation*}
x-\frac{\partial \pi _{\lambda }w_{z}^{u}}{\partial x}\frac{\partial w^{cs}}{%
\partial \lambda }x=0,
\end{equation*}%
which since $w_{z}^{u}$ and $w^{cs}$ are Lipschitz with constants $\rho <1/L$
and $L<1,$ respectively, implies that $x=0,$ and in turn that $\lambda =0.$
Since $B$ is invertible, it is evident that the rank of $A$ is $c+u+s,$
which implies (\ref{eq:V-span-2}).

Since $W_{z}^{u}\subset W^{cu}$ we see that $q^{\ast }\in W_{z}^{u}\cap
W^{cs}\subset W^{cu}\cap W^{cs}=\Lambda^*,$ hence $q^{\ast }\in
W_{z}^{u}\cap \Lambda^*.$ The fact that the intersection point is unique
follows from the fact that $\Lambda^*\subset W^{cu}$ and already established
uniqueness of the intersection point $W_{z}^{u}\cap W^{cs}$.
\end{proof}

%TCIDATA{Version=5.00.0.2606}
%TCIDATA{LaTeXparent=0,0,MMFedit.tex}

\section{Stable fibers\label{sec:stb-fibers}}

The goal of this section is to establish the existence of the foliation of $%
W^{cs}$ into the stable fibers $W_{z}^{s}$ for $z\in W^{cs}.$ In this
section $\theta =\left( \lambda ,x\right) $. Throughout the section we
assume that the assumptions of Theorem \ref{th:main} hold.

Let us fix a point $z\in W^{cs}$. Let $y\in \overline{B}_{s}$ and consider
the following problem: Find $\theta $ such that%
\begin{equation*}
\pi _{\theta }\left( f^{n}\left( \theta ,y\right) -f^{n}\left( z\right)
\right) =0,
\end{equation*}%
under the constraint%
\begin{equation*}
f^{i}\left( \theta ,y\right) \in D\qquad \text{for }i=1,\ldots ,n.
\end{equation*}%
By taking $b_{y}(\theta )=\left( \theta ,y\right) $ and observing that $%
f^{i}\left( \theta ,y\right) =f^{i}\left( b_{y}(\theta )\right) ,$ from
Lemma \ref{lem:center-unstable-disc} it follows immediately that this
problem has a unique solution $\theta \left( y\right) $ which is as smooth
as $f$.

We define 
\begin{equation}
d_{n,z}\left( y\right) =\left( \theta (y),y\right) .  \label{eq:dnz-def-stbl}
\end{equation}%
Our objective will be to prove that $d_{n,z}\left( y\right)$ are vertical
disks converging uniformly to $W_{z}^{s}$ as $n$ tends to infinity. First we
prove a technical lemma.

\begin{lemma}
\label{lem:stb-fib-conv}Assume that $f^{j}(z)\in D$ for $j=0,1,\dots ,n$. If 
$f^{j}(q_{i})\in J_{s}\left( f^{j}(z),1/L\right) \cap D$ for $i=1,2$, $%
j=1,\ldots ,n$ and if%
\begin{equation*}
q_{1}\notin J_{s}\left( q_{2},1/L\right)
\end{equation*}%
then%
\begin{eqnarray*}
\left\Vert \pi _{\theta }\left( q_{1}-q_{2}\right) \right\Vert &\leq &\frac{%
4R}{L}\left( \frac{\mu _{s,1}}{\xi _{cu,1,P}}\right) ^{n}, \\
\left\Vert q_{1}-q_{2}\right\Vert &\leq &\left( 1+L\right) \frac{4R}{L}%
\left( \frac{\mu _{s,1}}{\xi _{cu,1,P}}\right) ^{n}.
\end{eqnarray*}
\end{lemma}

\begin{proof}
Our assumption $f^{j}(q_{i})\in J_{s}\left( f^{j}(z),1/L\right) $, $i=1,2$
and $j=0,1,\dots ,n$ implies that $f^{j}(q_{1}),f^{j}(q_{2}),f^j(z)$ are
contained in the same charts for $j=0,1,\dots ,n$.

By Corollary \ref{cor:stable-lip}, $f^{j}(q_{1})\notin
J_{s}(f^{j}(q_{2}),1/L)$ for $j=0,1,\dots ,n$, hence $f^{j}(q_{1})\in
J_{cu}(f^{j}(q_{2}),L)$. By Lemma \ref{lem:cu-cone-xi-cu} this implies 
\begin{equation}
\Vert \pi _{\theta }(f^{n}(q_{1})-f^{n}(q_{2}))\Vert \geq \xi _{cu,1,P}\Vert
\pi _{\theta }(f^{n-1}(q_{1})-f^{n-1}(q_{2}))\Vert \geq \ldots \geq \xi
_{cu,1,P}^{n}\Vert \pi _{\theta }(q_{1}-q_{2})\Vert .
\label{eq:q--iter-estm}
\end{equation}

We estimate $\Vert \pi _{\theta }\left( f^{n}(q_{1})-f^{n}(q_{2})\right)
\Vert $ using the contraction in the $y$-direction. By Lemma \ref%
{lem:Js-contraction}, for $i=1,2,$%
\begin{eqnarray*}
\left\Vert \pi _{y}\left( f^{n}(q_{i})-f^{n}(z)\right) \right\Vert &\leq
&\mu _{s,1}\left\Vert \pi _{y}\left( f^{n-1}(q_{i})-f^{n-1}(z)\right)
\right\Vert \\
&\leq &\mu _{s,1}^{n}\left\Vert \pi _{y}\left( q_{i}-z\right) \right\Vert
\leq \mu _{s,1}^{n}2R.
\end{eqnarray*}%
Since $f^{n}(q_{i})\in J_{s}\left( f^{n}(z),1/L\right) ,$ 
\begin{equation*}
\left\Vert \pi _{\theta }\left( f^{n}(q_{i})-f^{n}(z)\right) \right\Vert
\leq \frac{1}{L}\left\Vert \pi _{y}\left( f^{n}(q_{i})-f^{n}(z)\right)
\right\Vert \leq \mu _{s,1}^{n}\frac{2R}{L}
\end{equation*}%
which by the triangle inequality implies%
\begin{eqnarray*}
\left\Vert \pi _{\theta }\left( f^{n}(q_{1})-f^{n}(q_{2})\right) \right\Vert
&\leq &\left\Vert \pi _{\theta }\left( f^{n}(q_{1})-f^{n}(z)\right)
\right\Vert +\left\Vert \pi _{\theta }\left( f^{n}(q_{2})-f^{n}(z)\right)
\right\Vert \\
&\leq &\mu _{s,1}^{n}\frac{4R}{L}.
\end{eqnarray*}%
Combining the above with (\ref{eq:q--iter-estm}),%
\begin{equation*}
\Vert \pi _{\theta }(q_{1}-q_{2})\Vert \leq \frac{4R}{L}\left( \frac{\mu
_{s,1}}{\xi _{cu,1,P}}\right) ^{n}.
\end{equation*}

Since $q_{1}\in J_{cu}(q_{2},L)$, then 
\begin{equation*}
\left\Vert \pi _{y}\left( q_{1}-q_{2}\right) \right\Vert \leq L\left\Vert
\pi _{\theta }\left( q_{1}-q_{2}\right) \right\Vert ,
\end{equation*}%
hence%
\begin{eqnarray*}
\left\Vert q_{1}-q_{2}\right\Vert &\leq &\left\Vert \pi _{y}\left(
q_{1}-q_{2}\right) \right\Vert +\left\Vert \pi _{\theta }\left(
q_{1}-q_{2}\right) \right\Vert \\
&\leq &\left( 1+L\right) \frac{4R}{L}\left( \frac{\mu _{s,1}}{\xi _{cu,1,P}}%
\right) ^{n},
\end{eqnarray*}%
as required.
\end{proof}

\begin{lemma}
\label{lem:Wsz-lem1}For any $n\geq 0,$ $d_{n,z}$ is a vertical disc and the
sequence $d_{n,z}$ converges uniformly to horizontal disk $d_{z}$. Moreover, 
$W_{z}^{s}=\left\{ \left( w_{z}^{s}\left( y\right) ,y\right) :y\in \overline{%
B}_{s}(R)\right\} $, where $w_{z}^{s}:\overline{B}_{s}(R)\rightarrow \Lambda
\times \overline{B}_{u}(R)$ and $d_{z}(y)=\left( w_{z}^{s}\left( y\right)
,y\right) $
\end{lemma}

\begin{proof}
Let $y_{1},y_{2}\in \overline{B}_{s}(R)$. By construction, 
\begin{equation}
\pi _{\left( \lambda ,x\right) }f^{n}\left( d_{n,z}\left( y_{1}\right)
\right) =\pi _{\left( \lambda ,x\right) }f^{n}\left( z\right) =\pi _{\left(
\lambda ,x\right) }f^{n}\left( d_{n,z}\left( y_{2}\right) \right) ,
\label{eq:dnz-in-cone-fin}
\end{equation}%
and $f^{i}\left( d_{n,z}\left( y_{1}\right) \right) ,f^{i}\left(
d_{n,z}\left( y_{2}\right) \right) \in D$ for $i=0,\ldots ,n$. Since (\ref%
{eq:dnz-in-cone-fin}) implies that%
\begin{equation*}
f^{n}\left( d_{n,z}\left( y_{1}\right) \right) \in J_{s}\left( f^{n}\left(
d_{n,z}\left( y_{2}\right) \right) ,1/L\right) ,
\end{equation*}%
by the backward cone condition, 
\begin{equation*}
d_{n,z}\left( y_{1}\right) \in J_{s}\left( d_{n,z}\left( y_{2}\right)
,1/L\right) ,
\end{equation*}%
which means that $d_{n,z}$ is a vertical disc. Also, for any $y\in \overline{%
B}_{s}(R),$ since $f^{n}\left( d_{n,z}\left( y\right) \right) \in
J_{s}\left( f^{n}\left( z\right) ,1/L\right)$, by the backward cone
condition, 
\begin{equation}
f^{j}\left( d_{n,z}\left( y\right) \right) \in J_{s}\left( f^{j}\left(
z\right) ,1/L\right) \quad \text{for }j=0,\ldots ,n.
\label{eq:dnz-in-stbl-cone}
\end{equation}

Observe that since 
\begin{equation}
\pi _{y}\left( d_{n,z}\left( y\right) -d_{n+j,z}\left( y\right) \right) =0,
\label{eq:tmp-dzn-ineq1}
\end{equation}%
we have $\Vert \pi _{\theta }(d_{n,z}\left( y\right) -d_{n+j,z}\left(
y\right) )\Vert =\Vert d_{n,z}\left( y\right) -d_{n+j,z}\left( y\right)
\Vert $. Assume that $d_{n,z}\left( y\right) \neq d_{n+j,z}\left( y\right) .$
By (\ref{eq:tmp-dzn-ineq1}) we see that $d_{n,z}\left( y\right) \notin
J_{s}\left( d_{n+j,z}\left( y\right) ,1/L\right) .$ From Lemma~\ref%
{lem:stb-fib-conv}\label{proof-tu-fiber-contr} applied to $%
q_{1}=d_{n,z}\left( y\right) $, $q_{2}=d_{n+j,z}\left( y\right) $ and $z$ it
follows that 
\begin{equation}
\Vert d_{n,z}(y)-d_{n+j,z}(y)\Vert =\Vert \pi _{\theta
}(d_{n,z}(y)-d_{n+j,z}(y))\Vert \leq \frac{4R}{L}\left( \frac{\mu _{s,1}}{%
\xi _{cu,1,P}}\right) ^{n}.  \label{eq:tmp-Wsz-lem1-1}
\end{equation}%
Note that if $d_{n,z}\left( y\right) =d_{n+j,z}\left( y\right) ,$ then (\ref%
{eq:tmp-Wsz-lem1-1}) also holds. Since by our assumptions $\frac{\mu _{s,1}}{%
\xi _{cu,1,P}}<1$, therefore $d_{n,z}$ is a Cauchy sequence in supremum
norm. Let us denote the limit by $d_{z}$.

The $d_{n,z}$ are vertical discs. The properties (\ref{eq:ver-disc-prop})
are preserved when passing to the limit, hence $d_{z}$ is a vertical disc.

We show that for all $y\in \overline{B}_{s}(R),$ $d_{z}(y)\in W^{cs}$. By
construction, for any $i\geq 0$ and $n\geq i$, $f^{i}\left( d_{n,z}\left(
y\right) \right) \in D.$ Passing to the limit with $n$ to infinity gives $%
f^{i}\left( d_{z}\left( y\right) \right) \in D$, as required.

By (\ref{eq:dnz-in-stbl-cone}), passing to the limit with $n$ to infinity,
we see that for any $j\geq 0$%
\begin{equation*}
f^{j}\left( d_{n,z}\left( y\right) \right) \in J_{s}\left( f^{j}\left(
z\right) ,1/L\right) ,
\end{equation*}%
hence $\{d_{z}(y)\ |\ y\in \overline{B}_{s}(R)\}\subset W_{z}^{s}.$

To show that $W_{z}^{s}\subset \{d_{z}(y)\ |\ y\in \overline{B}_{s}(R)\},$
let us consider $p\in W^{cs}$, such that $f^{j}(p)\in J_{s}\left(
f^{j}\left( z\right) ,1/L\right) \ $for all $j\geq 0$. Let $y=\pi _{y}p$. We
will show that $p=d_{z}(y)$. From Lemma~\ref{lem:stb-fib-conv}, (taking $%
q_{1}=p,$ and $q_{2}=d_{z}(y),$) it follows that 
\begin{equation*}
\Vert p-d_{z}(y)\Vert =\Vert \pi _{\theta }(p-d_{z}(y))\Vert \leq \frac{4R}{L%
}\left( \frac{\mu _{s,1}}{\xi _{cu,1,P}}\right) ^{n}\rightarrow 0,\quad
n\rightarrow \infty .
\end{equation*}%
Therefore $p=d_{z}(y)$.

The function $w_{z}^{s}$ can be defined as $w_{z}^{s}(y)=\pi _{\theta
}d_{z}(y).$
\end{proof}

\begin{lemma}
\label{lem:Cm-norm-bound-fib-s}Let $m\leq k.$ Let $d_{n,z}$ be the sequence
of vertical discs defined in (\ref{eq:dnz-def-stbl}). Assume that $d_{n,z}$
are $C^{m}$ and that for any $n,$ $\left\Vert \pi _{\theta
}d_{n,z}\right\Vert _{C^{m}}<c_{m},$ with $c_{m}$ independent of $n$. If $f$
satisfies rate conditions of order $m$, then $\left\Vert \pi _{\theta
}d_{n,z}\right\Vert _{C^{m+1}}<c_{m+1}$ for a constant independent of $n$.
\end{lemma}

\begin{proof}
The proof follows from identical arguments to the proof of Lemma \ref%
{lem:Cm-norm-bound-cs}. The noticeable difference is that when we apply
Theorem \ref{thm:stb-jet-propagation}, we should choose coordinates $\mathrm{%
x}=\theta =(\lambda ,x),$ $\mathrm{y}=y$ and constants $\xi =\xi _{cu,2}$, $%
\mu =\mu _{s,1}$. Note that conditions (\ref{eq:rate-cond-0}), (\ref%
{eq:rate-cond-3}) imply (\ref{eq:full-stb-rate-cond}) for any $m\geq 0$.
\end{proof}

\begin{lemma}
\label{lem:Wsz-lem2}For any $z\in W^{cs}$ the manifold $W_{z}^{s}$ is $C^{k}$%
.
\end{lemma}

\begin{proof}
The functions $\pi _{\theta }d_{n,z}$ are $C^{k+1}$ and uniformly Lipschitz
with constant $1/L$. The fact that $C^{k}$ smoothness is preserved as we
pass to the limit follows from Lemma \ref{lem:Cm-norm-bound-fib-s} and
mirror arguments to the proof of Lemma~\ref{lem:cu-smooth}.
\end{proof}

\begin{lemma}
\label{lem:stb-fiber-C-conv} Let $z\in W^{cs}$. If $p\in W_{z}^{s}$, then
for $n\geq 0$ $f^{n}(p),f^{n}(z)$ are in the same chart and 
\begin{equation*}
\left\Vert f^{n}(p)-f^{n}(z)\right\Vert \leq \left( 1+\frac{1}{L}\right)
\Vert \pi _{y}(p-z)\Vert \mu _{s,1}^{n},\qquad n\geq 0.
\end{equation*}
\end{lemma}

\begin{proof}
We first observe that for any $q_{1},q_{2}\in D,$ such that $q_{1}\in
J_{s}\left( q_{2},1/L\right) ,$ holds%
\begin{equation*}
\left\Vert \pi _{\theta }\left( q_{1}-q_{2}\right) \right\Vert \leq \frac{1}{%
L}\left\Vert \pi _{y}\left( q_{1}-q_{2}\right) \right\Vert ,
\end{equation*}%
hence 
\begin{equation*}
\left\Vert q_{1}-q_{2}\right\Vert \leq \left\Vert \pi _{\theta }\left(
q_{1}-q_{2}\right) \right\Vert +\left\Vert \pi _{y}\left( q_{1}-q_{2}\right)
\right\Vert \leq \left( 1+\frac{1}{L}\right) \left\Vert \pi _{y}\left(
q_{1}-q_{2}\right) \right\Vert .
\end{equation*}

Let $p\in W_{z}^{s}$, which means that $f^{i}(p)\in J_{s}\left(
f^{i}(z),1/L\right) \cap D,$ for all $i\geq 0$. From Lemma~\ref%
{lem:Js-contraction} it follows for $n>0$ that 
\begin{equation*}
\left\Vert \pi _{y}\left( f^{n}(p)-f^{n}(z)\right) \right\Vert \leq (\mu
_{s,1})^{n}\left\Vert \pi _{y}\left( p-z\right) \right\Vert .
\end{equation*}%
Therefore 
\begin{eqnarray*}
\left\Vert f^{n}(p)-f^{n}(z)\right\Vert &\leq &\left( 1+\frac{1}{L}\right)
\left\Vert \pi _{y}\left( f^{n}(p)-f^{n}(z)\right) \right\Vert \\
&\leq &\mu _{s,1}^{n}\left( 1+\frac{1}{L}\right) \left\Vert \pi _{y}\left(
p-z\right) \right\Vert ,
\end{eqnarray*}%
as required.
\end{proof}

\begin{lemma}
\label{lem:Wsz-lem3}For any $z\in W^{cs}$, let us define a set $U=U(z)$ as 
\begin{eqnarray*}
U &=&\{p\in D:f^{n}(p)\in D\text{ for all }n\geq 0,\text{ and} \\
&&\exists n_{0}\geq 0,\,\exists C>0\text{ (which may depend on $p$, $z$)} \\
&&\text{s.t. for }n\geq n_{0},\text{ $f^{n}(p)$ and $f^{n}(z)$ are in the
same good chart and } \\
&&\left\Vert f^{n}(p)-f^{n}(z)\right\Vert \leq C\mu _{s,1}^{n}\}.
\end{eqnarray*}%
Then $W_{z}^{s}=U.$
\end{lemma}

\begin{proof}
From Lemma~\ref{lem:stb-fiber-C-conv} it follows that $W_{z}^{s}\subset U$.
It remains to prove that $U\subset W_{z}^{s}$.

For the proof by the contradiction let us consider $p\in U\setminus
W_{z}^{s} $. Observe that from the backward cone condition (Definition~\ref%
{def:back-cc}), since $p\notin W_{z}^{s},$ it follows that for $i\geq n_{0}$
holds $f^{i}(p)\notin J_{s}\left( f^{i}(z),1/L\right) $. Hence $f^{i}(p)\in
J_{cu}(f^{i}(z),L)$ for any $i\geq n_{0}$ (this makes sense because $%
f^{i}(p) $ and $f^{i}(z)$ are in the same good chart for $i\geq n_{0}$.)
Hence from Lemma~\ref{lem:cu-cone-xi-cu} it follows that 
\begin{equation*}
\left\Vert \pi _{\theta }\left( f^{n_{0}+i}(p)-f^{n_{0}+i}(z)\right)
\right\Vert \geq \xi _{cu,1,P}^{i}\left\Vert \pi _{\theta }\left(
f^{n_{0}}(p)-f^{n_{0}}(z)\right) \right\Vert ,
\end{equation*}%
and thus for any $n>n_{0}$%
\begin{equation}
\left\Vert \pi _{\theta }\left( f^{n}(p)-f^{n}(z)\right) \right\Vert \geq
\xi _{cu,1,P}^{n}\left( \xi _{cu,1,P}^{-n_{0}}\left\Vert \pi _{\theta
}\left( f^{n_{0}}(p)-f^{n_{0}}(z)\right) \right\Vert \right) .
\label{eq:s-mu-contr-2}
\end{equation}%
By our assumption 
\begin{equation}
\left\Vert f^{n}(p)-f^{n}(z)\right\Vert \leq C\mu _{s,1}^{n},\quad n\geq
n_{0}.  \label{eq:s-mu-contr-1}
\end{equation}

Since $\mu _{s,1}<\xi _{cu,1,P}$, conditions (\ref{eq:s-mu-contr-1}) and (%
\ref{eq:s-mu-contr-2}) contradict each other. This means that $p\in
W_{z}^{s},$ as required.
\end{proof}

Lemmas \ref{lem:Wsz-lem1}, \ref{lem:Wsz-lem2}, \ref{lem:Wsz-lem3} combined
prove the claims about $W_{z}^{u}$ from Theorem \ref{th:main}. We now prove
Theorem \ref{th:wsz}, which can be used to obtain tighter Lipschitz bounds
on $w_{z}^{s}$.

\begin{proof}[Proof of Theorem \protect\ref{th:wsz}]
Observe that by definition of $d_{n,z}$, for any $y_{1},y_{2}\in \overline{B}%
_{s}\left( R\right) $%
\begin{equation*}
\pi _{\left( \lambda ,x\right) }\left( f^{n}\left( d_{n,z}\left(
y_{1}\right) \right) \right) =\pi _{\left( \lambda ,x\right) }\left(
f^{n}\left( z\right) \right) =\pi _{\left( \lambda ,x\right) }f^{n}\left(
d_{n,z}\left( y_{2}\right) \right) ,
\end{equation*}%
hence for $y_{1}\neq y_{2}$ 
\begin{equation}
f^{n}\left( d_{n,z}\left( y_{2}\right) \right) \notin J_{s}^{c}\left(
f^{n}\left( d_{n,z}\left( y_{1}\right) \right) ,M\right) .
\label{eq:tmp-jets-M-prop-1}
\end{equation}%
By Theorem \ref{thm:Lip-stable-jet}, taking $\mathrm{x}=x$ and $\mathrm{y}%
=\left( \lambda ,y\right) $, for $i=1,\ldots ,n$,%
\begin{equation}
f(J_{s}^{c}(f^{i-1}\left( d_{n,z}\left( y_{1}\right) \right) ,M)\cap
D)\subset J_{s}^{c}(f^{i}\left( d_{n,z}\left( y_{1}\right) \right) ,M).
\label{eq:tmp-jets-M-prop-2}
\end{equation}%
This means that 
\begin{equation}
d_{n,z}\left( y_{2}\right) \notin J_{s}^{c}\left( d_{n,z}\left( y_{1}\right)
,M\right) .  \label{eq:proof-jet-M-prop}
\end{equation}
(Since otherwise from (\ref{eq:tmp-jets-M-prop-2}) and (\ref%
{eq:tmp-jets-M-prop-1}) we would get a contradiction.) By definition of $%
d_{n,z}$, 
\begin{equation*}
\pi _{y}d_{n,z}\left( y_{i}\right) =y_{i}\qquad \text{for }i=1,2,
\end{equation*}%
which combined with (\ref{eq:proof-jet-M-prop}) gives%
\begin{equation*}
\left\Vert \pi _{\left( \lambda ,x\right) }\left( d_{n,z}\left( y_{1}\right)
-d_{n,z}\left( y_{2}\right) \right) \right\Vert \leq M\left\Vert \pi
_{y}\left( d_{n,z}\left( y_{1}\right) -d_{y}\left( y_{2}\right) \right)
\right\Vert =M\left\Vert y_{1}-y_{2}\right\Vert ,
\end{equation*}%
as required.
\end{proof}

\begin{proposition}
Let $z\in W^{cs}$. Then the intersection $W_{z}^{s}\cap W^{cu}$ consists of
a single point and is transversal. The intersection $W_{z}^{s}\cap \Lambda^*$
consists of a single point.
\end{proposition}

\begin{proof}
The result follows from similar arguments to the proof of Proposition \ref%
{prop:fiber-inter}.
\end{proof}

%TCIDATA{Version=5.00.0.2606}
%TCIDATA{LaTeXparent=0,0,MMFedit.tex}

\section{Invariant manifolds for vector bundles\label{sec:bundles}}

The previous discussion was focused on the setting where $\Lambda$ was a
torus. We now generalize the result for $\Lambda$ which are compact
manifolds without boundaries.

\subsection{Vector bundles}

We start by recalling the definition of the vector bundle \cite{Hi}.

\begin{definition}
Let $B,E$ be topological spaces. Let $p:E\rightarrow B$ be a continuous map.
A \emph{vector bundle chart }on $(p,E,B)$ with \emph{domain } $U$ and \emph{%
dimension} $n$ is a homeomorphism $\varphi:p^{-1}(U)\rightarrow U\times%
\mathbb{R}^{n}$, where $U\subset B$ is open and such that 
\begin{equation*}
\pi_{1}\circ\varphi(z)=p(z),\qquad\mbox{for $z \in
p^{-1}(U)$}.
\end{equation*}
We will denote such bundle chart by a pair $(\varphi,U)$.

For each $\lambda\in U$ we define the homeomorphism $\varphi_{\lambda}$ to
be the composition 
\begin{equation*}
\varphi_{\lambda}:p^{-1}(\lambda)\overset{\varphi}{\rightarrow}\{\lambda
\}\times\mathbb{R}^{n}\rightarrow\mathbb{R}^{n}.
\end{equation*}
A \emph{vector bundle atlas $\Phi$} on $(p,E,B)$ is a family of vector
bundle charts on $(p,E,B)$ with the values in the same $\mathbb{R}^{n}$,
whose domains cover $B$ and such that whenever $(\varphi,U)$ and $(\psi,V)$
are in $\Phi$ and $\lambda\in U\cap V$, the homeomorphism $%
\psi_{\lambda}\varphi_{\lambda}^{-1}:\mathbb{R}^{n}\rightarrow\mathbb{R}^{n}$
is linear. The map 
\begin{equation*}
U\cap V\ni\lambda\mapsto\psi_{\lambda}\varphi_{\lambda}^{-1}\in GL(n)
\end{equation*}
is continuous for all pairs of charts in $\Phi$.

A maximal vector bundle atlas $\Phi $ is a \emph{vector bundle structure} on 
$(p,E,B)$. We then call $\gamma =(p,E,B,\Phi )$ a \emph{vector bundle}
having \emph{(fibre) dimension }n\emph{, projection $p$, total space }$E$
and \emph{base space }$B$.

The \emph{fibre} over $\lambda \in B$ is the space $p^{-1}(\lambda )=\gamma
_{\lambda }=E_{\lambda }$. $\gamma _{\lambda }$ has the vector space
structure.

If the $E,B$ are $C^{r}$ manifolds and all maps appearing in the above
definition are $C^{r}$, then we will say that the bundle $(p,E,B,\Phi)$ is a 
$C^{r}$-bundle.
\end{definition}

One can introduce the notion of subbundles, morphisms etc. (see \cite{Hi}
and references given there). The fibers can have a structure: for example a
scalar product, a norm, which depend continuously on the base point.

\begin{definition}
\label{def:Banachbundle} We say that the vector bundle $\gamma $ is a Banach
vector bundle with fiber being the Banach space $(F,\Vert \cdot \Vert )$, if
for each $\lambda \in B$ the fiber $\gamma _{\lambda }$ is a Banach space
with norm $\Vert \cdot \Vert _{\lambda }$ such that for each bundle chart $%
(\varphi ,U)$ the map $\varphi _{\lambda }:E_{\lambda }\rightarrow F$ is an
isometry ($\Vert \varphi _{\lambda }(v)\Vert =\Vert v\Vert _{\lambda }$).
\end{definition}

For vector bundles $\gamma _{1}$, $\gamma _{2}$ over the same base space one
can define $\gamma =\gamma _{1}\oplus \gamma _{2}$ by setting $\gamma
_{\lambda }=\gamma _{1,\lambda }\oplus \gamma _{2,\lambda }$. In the
following, points in $\gamma _{1}\oplus \gamma _{2}$ will be denoted by a
triple $(\lambda ,v_{1},v_{2})$, where $\lambda \in B$, $v_{1}\in \gamma
_{1,\lambda }$ and $v_{2}\in \gamma _{2,\lambda }$. If $\gamma _{1}$ and $%
\gamma _{2}$ are both Banach bundles, then $\gamma _{1}\oplus \gamma _{2}$
is also a Banach vector bundle with the norm on $\gamma _{1,\lambda }\oplus
\gamma _{2,\lambda }$ defined by $\Vert (v_{1},v_{2})\Vert _{\lambda }=\sqrt{%
\Vert v_{1}\Vert _{\lambda }^{2}+\Vert v_{2}\Vert _{\lambda }^{2}}$. We will
also always assume that the atlas on bundle $\gamma _{1}\oplus \gamma _{2}$
respects this structure, namely if $(\eta ,U)$ is a bundle chart for $\gamma
_{1}\oplus \gamma _{2}$, then its restriction (obtained through projection)
to $\gamma _{i}$ is also a bundle chart for $\gamma _{i}$ for $i=1,2$.

\subsection{Formulation of the result}

Assume that, we have Banach vector bundles $\gamma _{u}$, $\gamma _{s}$, $%
\gamma =\gamma _{u}\oplus \gamma _{s}$. Let $u$ and $s$ be the fiber
dimension of $\gamma _{u}$ and $\gamma _{s}$, respectively. Let the base
space for $\gamma $, denoted by $\Lambda $, be a $C^{k}$ compact manifold
without boundary of dimension $c$. We consider $D\subset \gamma _{u}\oplus
\gamma _{s}$ defined as:%
\begin{equation*}
D=\{(\lambda ,v_{1},v_{2})\in \gamma _{u}\oplus \gamma _{s}\ |\ \lambda \in
\Lambda ,v_{1}\in \gamma _{u,\lambda },v_{2}\in \gamma _{s,\lambda },\quad
\Vert v_{1}\Vert \leq R,\quad \Vert v_{2}\Vert \leq R\}.
\end{equation*}

Consider a finite open covering $\left\{ \mathcal{U}_{i}\right\} $ of $%
\Lambda$ and an atlas $\left\{ \left( \eta_{i},\mathcal{U}_{i}\right)
\right\} $, where 
\begin{equation*}
\eta_{i}:\mathcal{U}_{i}\rightarrow\eta_{i}\left( \mathcal{U}_{i}\right)
\subset\mathbb{R}^{c}
\end{equation*}
are charts. We assume that there exists a $R_\Lambda >0$ such that for any $%
\lambda\in\Lambda$ there exists an $i$ such that 
\begin{equation}
\overline{B}_{c}\left( \eta_{i}(\lambda),R_{\Lambda}\right) \subset\eta
_{i}\left( \mathcal{U}_{i}\right) .  \label{eq:good-chart-cond}
\end{equation}
Also, we assume that for any $\eta_{i}$ there exists a $\lambda$ such that (%
\ref{eq:good-chart-cond}) holds true.

\begin{definition}
We refer to a chart $\left( \eta_{i},\mathcal{U}_{i}\right) $ satisfying (%
\ref{eq:good-chart-cond}) as a good chart for $\lambda$.
\end{definition}

We assume that for each $\left( \eta _{i},\mathcal{U}_{i}\right) $ we have a
vector bundle chart for $\gamma $ of the form $\varphi _{i}=\left( p,\varphi
_{i}^{u},\varphi _{i}^{s}\right) :p^{-1}\left( \mathcal{U}_{i}\right)
\rightarrow \mathcal{U}_{i}\times \mathbb{R}^{s}\times \mathbb{R}^{u}$. We
define maps 
\begin{equation*}
\tilde{\eta}_{i}:p^{-1}\left( \mathcal{U}_{i}\right) \rightarrow \mathbb{R}%
^{c}\times \mathbb{R}^{u}\times \mathbb{R}^{s},
\end{equation*}%
as%
\begin{equation*}
\tilde{\eta}_{i}\left( \lambda ,v_{1},v_{2}\right) =\left( \eta _{i}\left(
\lambda \right) ,\varphi _{i,\lambda }^{u}\left( v_{1}\right) ,\varphi
_{i,\lambda }^{s}\left( v_{2}\right) \right) ,
\end{equation*}%
and sets 
\begin{align*}
\widetilde{\mathcal{U}}_{i}& =p^{-1}\left( \mathcal{U}_{i}\right) \cap D, \\
D_{i}& =\tilde{\eta}_{i}\left( \widetilde{\mathcal{U}}_{i}\right) =\eta
_{i}\left( \mathcal{U}_{i}\right) \times B_{u}\left( R\right) \times
B_{s}\left( R\right) .
\end{align*}

\begin{definition}
We say that $\tilde{\eta}_{i}$ is a good chart for $z\in \gamma _{u}\oplus
\gamma _{s}$ if $\eta _{i}$ is a good chart for $p\left( z\right) \in
\Lambda $.
\end{definition}

We use a notation $z=\left( \lambda ,v_{1},v_{2}\right) $ for points in $%
\gamma _{u}\oplus \gamma _{s}$ and $\left( \theta ,x,y\right) \in \mathbb{R}%
^{c}\times \mathbb{R}^{u}\times \mathbb{R}^{s}$ to make the distinction
between those on the bundle and those in local coordinates.

We fix a constant $L\in \mathbb{R}$ satisfying 
\begin{equation*}
L \in \left(\frac{2R}{R_\Lambda},1 \right).
\end{equation*}

\begin{remark}
Using mirror arguments to those in Remark \ref{rem:cone-in-chart} we see
that for any $M\leq\frac{1}{L}$ and any good chart $\tilde{\eta}_{i}$ around 
$z\in D$, holds%
\begin{align*}
J_{s}\left( \tilde{\eta}_{i}\left( z\right) ,M\right) \cap D_{i} & \subset%
\overline{B}_{c}\left( \pi_{\theta}\tilde{\eta}_{i}(z),R_{\Lambda }\right)
\times\overline{B}_{u}\left( R\right) \times\overline{B}_{s}\left( R\right) ,
\\
J_{u}\left( \tilde{\eta}_{i}\left( z\right) ,M\right) \cap D_{i} & \subset%
\overline{B}_{c}\left( \pi_{\theta}\tilde{\eta}_{i}(z),R_{\Lambda }\right)
\times\overline{B}_{u}\left( R\right) \times\overline{B}_{s}\left( R\right) .
\end{align*}
This is important for us, since it is one of the reasons why the proof
presented in previous sections will also work for the current setting. For
instance, one of the founding blocks of the proof was the Lemma \ref%
{lem:unstable-disc}, which states that images of horizontal discs are
horizontal discs. Horizontal discs are contained in cones, and here we see
that the relevant fragments of the cones will lie in local coordinates. This
will allow us to consider conditions defined locally.
\end{remark}

We consider a map%
\begin{equation*}
f:D\rightarrow \gamma _{u}\oplus \gamma _{s}.
\end{equation*}%
For any $z\in D,$ a good chart $\tilde{\eta}_{i}$ around $z$ and a good
chart $\tilde{\eta}_{j}$ around $f(z)$ we can define (locally around $\tilde{%
\eta}_{i}\left( z\right) $)%
\begin{equation*}
f_{ji}:=\tilde{\eta}_{j}\circ f\circ \tilde{\eta}_{i}^{-1}.
\end{equation*}

For a chart $\eta_{i}$ we define sets $D_{i}^{+},D_{i}^{-}\subset \mathbb{R}%
^{c}\times\mathbb{R}^{u}\times\mathbb{R}^{s}$ as%
\begin{align*}
D_{i}^{+} & =\eta_{i}\left( \mathcal{U}_{i}\right) \times\overline{B}%
_{u}(R)\times\partial B_{s}(R), \\
D_{i}^{-} & =\eta_{i}\left( \mathcal{U}_{i}\right) \times\partial \overline{B%
}_{u}(R)\times B_{s}(R).
\end{align*}

\begin{definition}
\label{def:covering-general}We say that $f$ satisfies covering conditions if
for any $z\in D$ the following conditions hold:

For any good chart $\tilde{\eta}_{i}$ around $z$, there exists a good chart $%
\tilde{\eta}_{j}$ around $f(z)$ such that the set%
\begin{equation*}
U=J_{u}(\tilde{\eta}_{i}\left( z\right) ,1/L)\cap D_{i}
\end{equation*}
is contained in the domain of $f_{ji}$. Additionally, for $%
\theta^{\ast}=\pi_{\theta}\tilde{\eta}_{j}\left( f(z)\right) $, there exists
a homotopy%
\begin{equation*}
h:\left[ 0,1\right] \times U\rightarrow B_{c}\left( \theta^{\ast
},R_{\Lambda}\right) \times\mathbb{R}^{u}\times\mathbb{R}^{s},
\end{equation*}
and a linear map $A:\mathbb{R}^{u}\rightarrow\mathbb{R}^{u}$ which satisfy:

\begin{enumerate}
\item $h_{0}=f_{ji}|_{U},$

\item for any $\alpha\in\left[ 0,1\right] $, 
\begin{align}
h_{\alpha}\left( U\cap D_{i}^{-}\right) \cap D_{j} & =\emptyset ,
\label{eq:homotopy-exit-general} \\
h_{\alpha}\left( U\right) \cap D_{j}^{+} & =\emptyset ,
\label{eq:homotopy-enter-general}
\end{align}

\item $h_{1}\left( \theta,x,y\right) =\left( \theta^{\ast},Ax,0\right) $,

\item $A\left( \partial B_{u}(R)\right) \subset\mathbb{R}^{u}\setminus 
\overline{B}_{u}(R).$
\end{enumerate}
\end{definition}

\begin{definition}
For $z\in D$, we refer to $\left( \tilde{\eta}_{j},\tilde{\eta}_{i}\right) $
which satisfy the conditions of Definition \ref{def:covering-general} as a
good charts pair for $z$.
\end{definition}

We assume that for any $z\in D$ and any good chart $\eta_{j}$ for $z$ there
exists an $i$ such that $\left( \tilde{\eta}_{j},\tilde{\eta}_{i}\right) $
is a good charts pair.

We use a notation 
\begin{equation*}
\mathfrak{C}(z)=\left\{ \left( j,i\right) :\left( \tilde{\eta}_{j},\tilde{%
\eta}_{i}\right) \text{ is a good charts pair for }z\right\} ,
\end{equation*}%
\begin{equation*}
\mathfrak{C}=\bigcup_{z\in D}\mathfrak{C}(z).
\end{equation*}

We now define the constants%
\begin{align*}
\xi _{u,1,P}& =\inf_{\left( j,i\right) \in \mathfrak{C}}m\left[ \frac{%
\partial \left( f_{ji}\right) _{x}}{\partial x}(\mathrm{dom}\left(
f_{ji}\right) )\right] -\frac{1}{L}\sup_{\left( j,i\right) \in \mathfrak{C}%
,z\in \mathrm{dom}\left( f_{ji}\right) }\left\Vert \frac{\partial \left(
f_{ji}\right) _{x}}{\partial \left( \lambda ,y\right) }(z)\right\Vert , \\
\xi _{cu,1,P}& =\inf_{\left( j,i\right) \in \mathfrak{C}}m\left[ \frac{%
\partial \left( f_{ji}\right) _{(\lambda ,x)}}{\partial (\lambda ,x)}(%
\mathrm{dom}\left( f_{ji}\right) )\right] -L\sup_{\left( j,i\right) \in 
\mathfrak{C},z\in \mathrm{dom}\left( f_{ji}\right) }\left\Vert \frac{%
\partial \left( f_{ji}\right) _{(\lambda ,x)}}{\partial y}(z)\right\Vert .
\end{align*}%
Similarly, we define constants which are analogues of $\mu _{\iota ,\kappa
}, $ $\xi _{\iota ,\kappa }$, (for $\iota \in \left\{ u,s,cu,cs\right\} $
and $\kappa \in \{1,2\}$) from section \ref{sec:main-results}, by changing
the conditions under the sup and inf, from \textquotedblleft $z\in D$" to
\textquotedblleft $\left( j,i\right) \in \mathfrak{C},z\in \mathrm{dom}%
\left( f_{ji}\right) $".

\begin{definition}
We say that $f$ satisfies cone conditions and rate conditions of order $%
k\geq1$ if the inequalities from Definition \ref{def:rate-conditions} are
satisfied.
\end{definition}

\begin{definition}
We say that $f$ satisfies backward cone conditions if for any $z\in D$ and a
good charts pair $\left( \tilde{\eta}_{j},\tilde{\eta}_{i}\right) $ for $z,$
the following condition is fulfilled:

If $z^{\prime}\in D$, $f(z^{\prime})\in\widetilde{\mathcal{U}}_{j}$ and $%
\tilde{\eta}_{j}\left( f(z^{\prime})\right) \in J_{s}\left( \tilde{\eta }%
_{j}(f(z)),1/L\right) $ then $z^{\prime}\in\widetilde{\mathcal{U}}_{i}$ and%
\begin{equation*}
\tilde{\eta}_{i}(z^{\prime})\in J_{s}\left( \tilde{\eta}_{i}\left( z\right)
,1/L\right) .
\end{equation*}
\end{definition}

\begin{definition}
For $z\in D$, we refer to $\left( \tilde{\eta}_{i_{n}},\tilde{\eta}%
_{i_{n-1}},\ldots,\tilde{\eta}_{i_{0}}\right) $ as a good charts sequence
for $z,$ if $\left( \tilde{\eta}_{i_{k+1}},\tilde{\eta}_{i_{k}}\right) $ is
a good charts pair for $f^{k}\left( z\right) $.
\end{definition}

For a good chart sequence $\left( \tilde{\eta}_{i_{n}},\tilde{\eta}%
_{i_{n-1}},\ldots,\tilde{\eta}_{i_{0}}\right) $ we use a notation%
\begin{equation*}
f_{i_{n},\ldots,i_{0}}=f_{i_{n}i_{n-1}}\circ\ldots\circ f_{i_{2}i_{1}}\circ
f_{i_{1}i_{0}}.
\end{equation*}
We can now formulate our result.

\begin{theorem}
\label{th:main-bundles}If $f$ is $C^{k+1}$ and satisfies covering
conditions, rate conditions of order $k$ and backward cone conditions, then $%
W^{cs},W^{cu}$ and $\Lambda ^{\ast }$ are $C^{k}$ manifolds. In local
coordinates given by a chart $\tilde{\eta}_{i}$, the manifolds are graphs of 
$C^{k}$ functions 
\begin{align*}
w_{i}^{cs}& :\eta _{i}\left( \mathcal{U}_{i}\right) \times \overline{B}%
_{s}(R)\rightarrow \overline{B}_{u}(R), \\
w_{i}^{cu}& :\eta _{i}\left( \mathcal{U}_{i}\right) \times \overline{B}%
_{u}(R)\rightarrow \overline{B}_{s}(R), \\
\chi _{i}& :\eta _{i}\left( \mathcal{U}_{i}\right) \rightarrow \overline{B}%
_{u}(R)\times \overline{B}_{s}(R),
\end{align*}%
meaning that 
\begin{align*}
\tilde{\eta}_{i}\left( W^{cs}\cap \widetilde{\mathcal{U}}_{i}\right) &
=\left\{ \left( \theta ,w_{i}^{cs}(\theta ,y),y\right) :\theta \in \eta
_{i}\left( \mathcal{U}_{i}\right) ,y\in \overline{B}_{s}(R)\right\} , \\
\tilde{\eta}_{i}\left( W^{cu}\cap \widetilde{\mathcal{U}}_{i}\right) &
=\left\{ \left( \theta ,x,w_{i}^{cu}(\theta ,y)\right) :\theta \in \eta
_{i}\left( \mathcal{U}_{i}\right) ,x\in \overline{B}_{u}(R)\right\} , \\
\tilde{\eta}_{i}\left( \Lambda ^{\ast }\cap \widetilde{\mathcal{U}}%
_{i}\right) & =\left\{ \left( \theta ,\chi (\theta )\right) :\theta \in \eta
_{i}\left( \mathcal{U}_{i}\right) \right\} .
\end{align*}%
Moreover, $f_{|W^{cu}}$ is an injection.

For any $z\in W^{cs}$ and any good chart $\tilde{\eta}_{i}$ around $z,$ 
\begin{equation*}
\tilde{\eta}_{i}\left( W^{cs}\cap\widetilde{\mathcal{U}}_{i}\right) \subset
J_{cs}\left( \tilde{\eta}_{i}\left( z\right) ,L\right) .
\end{equation*}
For any $z\in W^{cu}$ and any good chart $\tilde{\eta}_{i}$ around $z,$%
\begin{equation*}
\tilde{\eta}_{i}\left( W^{cu}\cap\widetilde{\mathcal{U}}_{i}\right) \subset
J_{cu}\left( \tilde{\eta}_{i}\left( z\right) ,L\right) ,
\end{equation*}
Also, for any $z\in\Lambda^{\ast}$ and any good chart $\tilde{\eta}_{i}$
around $z,$ for $M=\frac{\sqrt{2}L}{\sqrt{1-L^{2}}},$%
\begin{equation*}
\tilde{\eta}_{i}\left( \Lambda^{\ast}\cap\widetilde{\mathcal{U}}_{i}\right)
\subset\left\{ \left( \theta,x,y\right) :\left\Vert \left( x,y\right)
-\pi_{\left( x,y\right) }\tilde{\eta}_{i}\left( z\right) \right\Vert \leq
M\left\Vert \theta-\pi_{\theta}\tilde{\eta}_{i}\left( z\right) \right\Vert
\right\} .
\end{equation*}
The manifolds $W^{cs}$ and $W^{cu}$ intersect transversally, and $W^{cs}\cap
W^{cu}=\Lambda^{\ast}$.

The manifolds $W^{cs}$ and $W^{cu}$ are foliated by invariant fibers $%
W_{z}^{s}$ and $W_{z}^{u}$, which in local coordinates given by any good
chart $\tilde{\eta}_{i}$ around $z$ are graphs of $C^{k}$ functions 
\begin{align*}
w_{z,i}^{s} & :\overline{B}_{s}(R)\rightarrow\eta_{i}\left( \mathcal{U}%
_{i}\right) \times\overline{B}_{u}(R), \\
w_{z,i}^{u} & :\overline{B}_{u}(R)\rightarrow\eta_{i}\left( \mathcal{U}%
_{i}\right) \times\overline{B}_{s}(R),
\end{align*}
The functions $w_{z,i}^{s}$ and $w_{z,i}^{u}$ are Lipschitz with constants $%
1/L$. Moreover, 
\begin{align*}
W_{z}^{s} & =\{p\in D:f^{n}(p)\in D\text{ for all }n\geq0,\text{ and} \\
& \left. \exists\text{ }C>0\text{ (which can depend on }p\text{)}\right. \\
& \left. \text{s.t. for }n\geq0,\text{ and any good charts sequence }\left( 
\tilde{\eta}_{i_{n}},\ldots,\tilde{\eta}_{i_{0}}\right) \text{ for }z\right.
\\
& \left. \left\Vert
f_{i_{n},\ldots,i_{0}}(p)-f_{i_{n},\ldots,i_{0}}(z)\right\Vert \leq
C\mu_{s,1}^{n}\right\} ,
\end{align*}
and%
\begin{align*}
W_{z}^{u} & =\{p\in W^{cu}:\text{the unique backward trajectory }%
\{p_{i}\}_{i=-\infty}^{0}\text{ of }p\text{ in }D,\text{ and for} \\
& \left. \text{any such }\{p_{i}\},\text{ and the unique backward trajectory 
}\{z_{i}\}_{i=-\infty}^{0}\text{ of }z\text{ in }D\text{ }\right. \\
& \left. \exists\text{ }C>0\text{ (which can depend on }p\text{)}\right. \\
& \left. \text{s.t. for }n\geq0,\text{ and any good charts sequence }\left( 
\tilde{\eta}_{i_{0}},\ldots,\tilde{\eta}_{i_{-n}}\right) \text{ for }z_{-n}%
\text{ }\right. \\
& \left. \left\Vert \tilde{\eta}_{i_{-n}}\left( p_{-n}\right) -\tilde{\eta }%
_{i_{-n}}\left( z_{-n}\right) \right\Vert \leq C\xi_{u,1,P}^{-n}\right.
\end{align*}
\end{theorem}

\subsection{Outline of the proof}

The proof of the theorem follows from the same arguments as the proof of
Theorem \ref{th:main}. The only difference is that instead of investigating
compositions $f^{n}$, we consider good chart sequences $\left( \tilde{\eta }%
_{i_{n}},\ldots,\tilde{\eta}_{i_{0}}\right) $ and local maps $%
f_{i_{n},\ldots,i_{0}}$.

We shall now focus on the needed changes to perform the construction. We
first go over the construction of the center-unstable manifold (see section %
\ref{sec:wcu-exists}). The construction of the center-unstable manifold is
based on propagation of horizontal discs. In our context we modify the
definition of the horizontal disc as follows:

\begin{definition}
\label{def:hor-disc-general}We say that a set $b\subset D$ is a horizontal
disc if for any $z\in b$ there exists a good chart $\tilde{\eta}_{i}$ around 
$z$ such that $b\subset\widetilde{\mathcal{U}}_{i}$ and a continuous
function $b_{i}:\overline{B}_{u}(R)\rightarrow D_{i}$%
\begin{align}
\tilde{\eta}_{i}(b) & =b_{i}\left( \overline{B}_{u}(R)\right) ,  \notag \\
\pi_{x}b_{i}(x) & =x,  \notag \\
b_{i}\left( \overline{B}_{u}(R)\right) & \subset J_{u}\left(
b_{i}(x),1/L\right) \qquad\text{for any }x\in\overline{B}_{u}\left( R\right)
.  \label{eq:hor-disc-prop-general}
\end{align}
We say that $b$ is $C^{k}$ if $b_{i}$ are $C^{k}$.
\end{definition}

With such definition we have a mirror result to Lemma \ref{lem:unstable-disc}
. This is done in Lemma \ref{lem:unstable-disc-general}.

\begin{remark}
Lemmas \ref{lem:unstable-disc-general}, \ref%
{lem:center-unstable-disc-general} are the core of the construction of both $%
W^{cu}$ and $W^{cs}.$ For this reason we go into some degree of detail
outlining its proof, pointing out differences in approach when working with
local maps.
\end{remark}

\begin{lemma}
\label{lem:unstable-disc-general}Assume that $b\subset D$ is a horizontal
disc. If $f$ satisfies covering conditions and rate conditions of order $%
k\geq0$, then there exists a horizontal disc $b^{\ast}\subset D$ such that $%
f\left( b\right) \cap D=b^{\ast}$. Moreover if $b$ and $f$ are $C^{k}$ then
so is $b^{\ast}$.
\end{lemma}

\begin{proof}
The proof is a mirror argument to the proof of Lemma \ref{lem:unstable-disc}%
. We therefore restrict our attention to setting up the local maps needed
for the construction.

Let let us fix $z\in b.$ Let $\tilde{\eta}_{i}$ be a good chart around $z$,
for which conditions (\ref{eq:hor-disc-prop-general}) hold. Let $\tilde{\eta 
}_{j}$ be the good chart around $f(z_{0})$ from Definition \ref%
{def:covering-general}. Note that $\left( \tilde{\eta}_{j},\tilde{\eta }%
_{i}\right) $ is a good charts pair. Let $\theta^{\ast}=\pi_{\theta}\tilde{%
\eta}_{j}\left( f(z)\right) $. Existence and smoothness of $b_{j}^{\ast}:%
\overline{B}_{u}\left( R\right) \rightarrow D_{j}$ such that%
\begin{equation*}
f_{ji}\circ b_{i}\left( \overline{B}_{u}\left( R\right) \right) \cap
D_{i}=b_{j}^{\ast}\left( \overline{B}_{u}\left( R\right) \right)
\end{equation*}
follows from a mirror construction to the one from the proof of Lemma \ref%
{lem:unstable-disc}. We can define 
\begin{equation*}
b^{\ast}=\tilde{\eta}_{j}^{-1}\circ b_{j}^{\ast}\left( \overline{B}%
_{u}\left( R\right) \right) .
\end{equation*}
By construction, $b_{j}^{\ast}$ satisfies (\ref{eq:hor-disc-prop-general}).

Let us now take any $\hat{z}\in b^{\ast}$. We need to prove that we have a
good chart $\tilde{\eta}_{\hat{\jmath}}$ for $\hat{z},$ for which conditions
from Definition \ref{def:hor-disc-general} hold. By our construction $\hat {z%
}=f(\hat{z}_{0})$, for some $\hat{z}_{0}=\hat{z}_{0}\left( \hat{z}\right)
\in b$. Let $\tilde{\eta}_{\hat{\imath}}$ be the good chart around $\hat {z}%
_{0}$ for which conditions \ref{def:hor-disc-general} hold for $b$. Let $%
\tilde{\eta}_{\hat{\jmath}}$ be the good chart around $\hat{z}=f(\hat{z}%
_{0}) $ from Definition \ref{def:covering-general}. Once again, from the
same construction as in the proof of Lemma \ref{lem:unstable-disc} follows
the existence and smoothness of $b_{\hat{\jmath}}^{\ast}$.
\end{proof}

\begin{remark}
From the proof we see that for $z\in b$ such that $f(z)\in b^{\ast}\subset D$
and for a good charts pair $\left( \tilde{\eta}_{j},\tilde{\eta}_{i}\right) $
for $z$ we can construct $b_{j}$ satisfying (\ref{def:hor-disc-general}).
The chart $\tilde{\eta}_{j}$ is a good chart for $f(z)$.
\end{remark}

\begin{definition}
Let $D_{u}=\pi _{\gamma _{u}}D\subset \gamma _{u}$. We say that a function $%
b:D_{u}\rightarrow D$ is a center-horizontal disc if for any $\left( \lambda
,v_{1}\right) \in D_{u}$%
\begin{equation*}
\pi _{\gamma _{u}}b\left( \lambda ,v_{1}\right) =\left( \lambda
,v_{1}\right) ,
\end{equation*}%
and for any $z\in b\left( \gamma _{u}\right) $ and any good chart $\tilde{%
\eta}_{i}$ around $z$%
\begin{equation}
\tilde{\eta}_{i}\circ b\left( D_{u}\right) \cap D_{i}\subset J_{cu}\left( 
\tilde{\eta}_{i}(z),L\right) .  \label{eq:ch-cone-cond-general}
\end{equation}
\end{definition}

\begin{lemma}
\label{lem:center-unstable-disc-general}Assume that $b$ is a
center-horizontal disc. If $f$ satisfies covering conditions, backward cone
conditions and rate conditions of order $l\geq0$, then there exists a
center-horizontal disc $b^{\ast}$ such that 
\begin{equation*}
f\left( b\left( D_{u}\right) \right) \cap D=b^{\ast}(D_{u}).
\end{equation*}
Moreover, if $f$ and $b$ are $C^{k}$, then so is $b^{\ast}$.
\end{lemma}

\begin{proof}
The proof goes along the same lines as the proof of Lemma \ref%
{lem:center-unstable-disc}. We will outline the differences concerning the
choices of local maps.

We start by showing that%
\begin{equation}
f\circ b\left( D_{u}\right) \cap D\neq \emptyset .
\label{eq:chd-image-on-D-2}
\end{equation}%
To this end, we consider $\lambda \in \Lambda $ and define $b^{\lambda
}=\gamma _{\lambda }\cap b\left( D_{u}\right) $. By mirror arguments to the
ones from the proof of Lemma \ref{lem:center-unstable-disc} it follows that $%
b^{\lambda }$ is a horizontal disc. By Lemma \ref{lem:unstable-disc-general} 
$f(b^{\lambda })\cap D\neq \emptyset $, which implies (\ref%
{eq:chd-image-on-D-2}).

Using the same arguments as those from the proof of Lemma \ref%
{lem:center-unstable-disc} it follows that $\pi _{\gamma _{u}}\circ f\circ
b:D_{u}\rightarrow \gamma _{u}$ is an injective open map. By (\ref%
{eq:chd-image-on-D-2}) $\pi _{\gamma _{u}}\circ f\circ b\left( D_{u}\right)
\cap D_{u}\neq \emptyset $. If $\left( \lambda ,v_{1}\right) \in \partial
D_{u}$ then $\left\Vert v_{1}\right\Vert =R$. Let $z=\left( \lambda
,v_{1},v_{2}\right) =b\left( \lambda ,v_{1}\right) $. Let $\left( \tilde{\eta%
}_{j},\tilde{\eta}_{i}\right) $ be a good charts pair for $z$ and $%
U=J_{u}\left( \tilde{\eta}_{i}\left( z\right) ,1/L\right) $. Using the same
argument as in the proof of Lemma \ref{lem:center-unstable-disc} it follows
that $\pi _{\theta }f_{ji}\left( D_{i}^{-}\cap U\right) \cap D_{j}=\emptyset 
$. Thus $\pi _{\gamma _{u}}\circ f\circ b\left( \partial D_{u}\right) \cap
D_{u}=\emptyset $. This means that 
\begin{equation*}
\pi _{\gamma _{u}}\circ f\circ b\left( D_{u}\right) \cap D_{u}=D_{u},
\end{equation*}%
hence for any $\left( \lambda ^{\ast },v_{1}^{\ast }\right) \in D_{u}$ there
exists a $\left( \lambda ,v_{1}\right) \in D_{u}$ such that $\pi _{\gamma
_{u}}\circ f\circ b\left( \lambda ,v_{1}\right) =\left( \lambda ^{\ast
},v_{1}^{\ast }\right) .$ We can define $b^{\ast }\left( \lambda ^{\ast
},v_{1}^{\ast }\right) =f\circ b\left( \lambda ,v_{1}\right) $. All the
desired properties of $b^{\ast }$ follow from mirror arguments to the proof
of Lemma \ref{lem:center-unstable-disc}.
\end{proof}

For a center-horizontal disc $b$ we use the notation $\mathcal{G}\left(
b\right) $ for the center-horizontal disc $b^{\ast }$ from Lemma \ref%
{lem:center-unstable-disc-general}.

\begin{lemma}
\label{lem:Wcu-general}Let $b_{0}:D_{u}\rightarrow D$ be defined as $%
b_{0}\left( \lambda ,v_{1}\right) =\left( \lambda ,v_{1},0\right) $. If
assumptions of Theorem \ref{th:main-bundles} are satisfied, then $\mathcal{G}%
^{k}\left( b_{0}\right) $ converges to $W^{cu}$ as $k$ tends to infinity.
\end{lemma}

\begin{proof}
Let us fix $\left( \lambda ,v_{1}\right) \in D_{u}=\pi _{\gamma _{u}}D$. Let 
$k_{2}\geq k_{1}$ and let us define $q_{0}^{k_{1}}=\mathcal{G}%
^{k_{1}}b\left( \lambda ,v_{1}\right) $ and $q_{0}^{k_{2}}=\mathcal{G}%
^{k_{2}}b\left( \lambda ,v_{2}\right) $. By definition of $\mathcal{G}$,
there exist backward trajectories $\{q_{i}^{k_{1}}\}_{i=-k_{1}}^{0},$ $%
\{q_{i}^{k_{2}}\}_{i=-k_{2}}^{0}$ 
\begin{equation*}
f\left( q_{i}^{k_{1}}\right) =q_{i+1}^{k_{1}},\qquad f\left(
q_{i}^{k_{1}}\right) =q_{i+1}^{k_{1}}.
\end{equation*}%
Since $\pi _{\gamma _{u}}\mathcal{G}^{k_{i}}b\left( \lambda ,v_{1}\right)
=\left( \lambda ,v_{1}\right) $, we see that $q_{0}^{k_{1}}-q_{0}^{k_{2}}=%
\pi _{v_{2}}\left( q_{0}^{k_{1}}-q_{0}^{k_{2}}\right) $ hence we can compute%
\begin{equation*}
\left\Vert q_{0}^{k_{1}}-q_{0}^{k_{2}}\right\Vert =\left\Vert \pi
_{v_{2}}\left( q_{0}^{k_{1}}-q_{0}^{k_{2}}\right) \right\Vert ,
\end{equation*}%
and the norm is independent of the considered chart. Let us take a good
chart sequence $\left( \tilde{\eta}_{i_{k_{1}}},\tilde{\eta}%
_{i_{k_{1}-1}},\ldots ,\tilde{\eta}_{i_{0}}\right) $ for $%
q_{-k_{1}}^{k_{1}}. $ Since $f$ satisfies backward cone conditions, $\left( 
\tilde{\eta}_{i_{k_{1}}},\tilde{\eta}_{i_{k_{1}-1}},\ldots ,\tilde{\eta}%
_{i_{0}}\right) $ is also a good sequence for $q_{-k_{1}}^{k_{2}}.$ From
mirror computations to the ones from the proof of Lemma \ref{lem:w-cu} (see (%
\ref{eq:Wcu-back-orb-dist}))%
\begin{eqnarray*}
\left\Vert \mathcal{G}^{k_{1}}b\left( \lambda ,v_{1}\right) -\mathcal{G}%
^{k_{2}}b\left( \lambda ,v_{1}\right) \right\Vert &=&\Vert
q_{0}^{k_{1}}-q_{0}^{k_{2}}\Vert \\
&\leq &(1+1/L)2R(\mu _{s,1})^{k_{1}}.
\end{eqnarray*}%
We note that the estimate is independent of the choice of the good chart
sequence. Thus we obtain uniform convergence of $\mathcal{G}^{k}b\left(
\lambda ,v_{1}\right) $.

The proof of the fact that $\mathcal{G}^{k}b\left( \lambda ,v_{1}\right) $
converges to $W^{cu}$ follows from mirror arguments to the ones in the proof
of Lemma \ref{lem:w-cu}.
\end{proof}

Lemma \ref{lem:Wcu-general} establishes the existence of $W^{cu}$. The proof
of its smoothness follows from arguments identical to the proof of the
smoothness when $\Lambda $ was a torus (Lemma \ref{lem:cu-smooth}). All the
arguments in the proof are local, and can be performed using local maps
passing through good chart sequences.

We now move to outlining the method for the proof of the existence of $%
W^{cs} $. First we give two definitions.

\begin{definition}
We say that $b\subset D$ is a vertical disc if for any $z\in b$ there exists
a good chart $\tilde{\eta}_{i}$ around $z$ such that $b\subset \widetilde{%
\mathcal{U}}_{i}$ and a continuous function $b_{i}:\overline{B}%
_{s}(R)\rightarrow D_{i}$ 
\begin{eqnarray*}
\tilde{\eta}_{i}(b) &=&b_{i}\left( \overline{B}_{s}(R)\right) , \\
\pi _{y}b_{i}(y) &=&y, \\
b_{i}\left( \overline{B}_{s}(R)\right) &\subset &J_{s}\left(
b_{i}(y),1/L\right) \qquad \text{for any }y\in \overline{B}_{s}\left(
R\right) .
\end{eqnarray*}%
We say that $b$ is $C^{k}$ if $b_{i}$ are $C^{k}$.
\end{definition}

\begin{definition}
Let $D_{s}=\pi _{\gamma _{s}}D\subset \gamma _{s}$. We say that a continuous
function $b:D_{s}\rightarrow D$ is a center-vertical disc if for any $\left(
\lambda ,v_{2}\right) \in D_{s}$%
\begin{equation*}
\pi _{\gamma _{s}}b(\lambda ,v_{2})=\left( \lambda ,v_{2}\right)
\end{equation*}%
and for any $z\in b\left( \gamma _{s}\right) $ and any good chart $\tilde{%
\eta}_{i}$ around $z$%
\begin{equation*}
\tilde{\eta}_{i}\circ b\left( D_{u}\right) \cap D_{i}\subset J_{cs}\left( 
\tilde{\eta}_{i}(z),L\right) .
\end{equation*}
\end{definition}

The construction of $W^{cs}$ is analogous to the one from section \ref%
{sec:wcs-exists}: For any $i\in \mathbb{Z}_{+}$ and $(\lambda ,v_{2})\in
D_{s}$ we consider the following problem: Find $x$ such that 
\begin{equation*}
\pi _{v_{1}}f^{i}(\lambda ,v_{1},v_{2})=0
\end{equation*}%
under the constraint 
\begin{equation*}
f^{l}(\lambda ,v_{1},v_{2})\in D,\quad l=0,1,\dots ,i.
\end{equation*}%
From Lemma \ref{lem:unstable-disc-general} it follows that this problem has
a unique solution $v_{1,i}(\lambda ,v_{2})$ which is as smooth as $f$.

We consider $b_{i}:D_{s}\rightarrow D$ given by $b_{i}(\lambda
,v_{2})=(\lambda ,v_{1,i}(\lambda ,v_{2}),v_{2})$. Then, using mirror
arguments to the proof of Lemma \ref{lem:wcs-graph-lim}, $b_{i}$ is a center
vertical disc and the sequence $b_{i}$ converges uniformly to $W^{cs}$. The
proof of the smoothness of $W^{cs}$ follows from mirror arguments to the
proof of Lemma \ref{lem:cs-smooth}.

Intersection of $W^{cu}$ and $W^{cs}$ gives the center manifold $\Lambda
^{\ast }$.

Vertical and horizontal discs are contained in local charts. Thus the
arguments for the existence of $W_{z}^{s}$ and $W_{z}^{u}$ follow from
identical arguments as those from sections \ref{sec:unstb-fibers}, \ref%
{sec:stb-fibers}. The only difference is that instead of working with
compositions of $f$, we work with compositions of local maps passing through
good chart sequences.

%TCIDATA{Version=5.00.0.2606}
%TCIDATA{LaTeXparent=0,0,MMFedit.tex}

\section{Numerical example\label{sec:num}}

We consider a one dimensional torus (circle) $\Lambda $ and the rotating H%
\'{e}non map $F_{\varepsilon }:\Lambda \times \mathbb{R}^{2}\rightarrow
\Lambda \times \mathbb{R}^{2},$%
\begin{equation}
F_{\varepsilon }(\lambda ,q_{1},q_{2})=(\theta +c+\varepsilon q_{1}\cos
(2\pi \lambda )),1+q_{2}-aq_{1}^{2}+\varepsilon \cos (2\pi \lambda ),bq_{1}).
\label{eq:rotating-henon}
\end{equation}%
We take $a=0.68,$ $b=0.1$ and an arbitrary constant $c\in \mathbb{R}$. We
investigate the existence and smoothness of the NHIM and its associated
stable/unstable manifolds for a range of parameters $\boldsymbol{\varepsilon 
}=\left[ \varepsilon _{1},\varepsilon _{2}\right] .$

We consider the maps (\ref{eq:rotating-henon}) in local coordinates $\left(
\lambda ,x,y\right) $ given by the linear change%
\begin{equation*}
\left( \lambda ,q_{1},q_{2}\right) =C\left( \lambda ,x,y\right) +\left(
0,q_{1}^{\ast },q_{2}^{\ast }\right) ,
\end{equation*}%
where%
\begin{align*}
q_{1}^{\ast }& =\frac{-(1-b)-\sqrt{(1-b)^{2}+4a}}{2a}\approx -2.043\,3, \\
q_{2}^{\ast }& =bq_{1}^{\ast }\approx -0.204\,33.
\end{align*}%
and%
\begin{equation*}
C=\left( 
\begin{array}{lll}
1 & 0 & 0 \\ 
0 & 1 & -0.3553203857 \\ 
0 & 0.03553203857 & 1%
\end{array}%
\right) .
\end{equation*}%
Thus, in local coordinates $p=\left( \lambda ,x,y\right) $, we consider the
family of maps 
\begin{equation*}
f_{\varepsilon }\left( p\right) =F_{\varepsilon }\left( Cp+\left(
0,q_{1}^{\ast },q_{2}^{\ast }\right) \right) -\left( 0,q_{1}^{\ast
},q_{2}^{\ast }\right) .
\end{equation*}%
The choice of $\left( q_{1}^{\ast },q_{2}^{\ast }\right) $ is dictated by
the fact that this is a hyperbolic fixed point for the H\'{e}non map (with $%
\varepsilon =0$). The matrix $C$ diagonalizes (roughly) the linear part of $F
$ into a Jordan normal form.

For a fixed interval $\boldsymbol{\varepsilon }=\left[ \varepsilon
_{1},\varepsilon _{2}\right] $, we consider the set $D_{\boldsymbol{%
\varepsilon }}=\Lambda \times \overline{B}_{u=1}\left( R\right) \times 
\overline{B}_{s=1}\left( R\right) $, with $R=\varepsilon _{2}.$ Below we
take two examples of $\boldsymbol{\varepsilon }=\left[ 0,0.0001\right] $ and 
$\boldsymbol{\varepsilon }=\left[ 0.009,0.01\right] $. The bounds for $\left[
Df_{\boldsymbol{\varepsilon }}\left( D_{\boldsymbol{\varepsilon }}\right) %
\right] $ for these two intervals are:%
\begin{eqnarray}
\left[ Df_{\left[ 0,0.0001\right] }\left( D_{\left[ 0,0.0001\right] }\right) %
\right] &=&\left( 
\begin{array}{lll}
1_{-0.00129}^{+0.00129} & {0}_{-0.000100}^{+0.000101} & {0}%
_{-0.000036}^{+0.000036} \\ 
{0}_{-0.000621}^{+0.000621} & 2.814_{17}^{55} & {0}_{-0.000065}^{+0.000065}
\\ 
{0}_{-0.000023}^{+0.000023} & {0}_{-0.000007}^{+0.000007} & -0.0355_{29}^{35}%
\end{array}%
\right) ,  \label{eq:DF-bound1} \\
\left[ Df_{\left[ 0.009,0.01\right] }\left( D_{\left[ 0.009,0.01\right]
}\right) \right] &=&\left( 
\begin{array}{lll}
1_{-0.129\,24}^{+0.12924} & {0}_{-0.010001}^{+0.010001} & {0}%
_{-0.003554}^{+0.003554} \\ 
{0}_{-0.062049}^{+0.062049} & 2._{79615}^{83257} & {0}%
_{-0.006468}^{+0.006468} \\ 
{0}_{-0.002205}^{+0.002205} & {0}_{-0.000647}^{+0.000647} & 
-0.035_{302}^{762}%
\end{array}%
\right) .  \label{eq:DF-bound2}
\end{eqnarray}%
Above, by convention, $1_{-0.00129}^{+0.00129}$ stands for the interval $%
\left[ 1-0.00129,1+0.00129\right] $ and $2.814_{17}^{55}$ stands for $\left[
2.81417,2.81455\right] $. We choose $L=1-\frac{1}{100}$, and in Table \ref%
{Table:coefficients} display coefficients that were computed based on (\ref%
{eq:DF-bound1}) and (\ref{eq:DF-bound2}).

\begin{table}[h]
\begin{center}
\begin{tabular}{lll}
\hline\hline
& $\boldsymbol{\varepsilon}=\left[ 0,0.0001\right] $ & $\boldsymbol{%
\varepsilon }=\left[ 0.009,0.01\right] $ \\ \hline
$\xi_{u,1}$ & $2.81352$ & $2.7303$ \\ 
$\xi_{u,1,P}$ & $2.81352$ & $2.7303$ \\ 
$\xi_{u,2}$ & $2.81408$ & $2.78624$ \\ 
$\xi_{cu,1}$ & $0.997718$ & $0.748463$ \\ 
$\xi_{cu,2}$ & $0.997766$ & $0.753236$ \\ 
$\xi_{cu,1,P}$ & $0.997718$ & $0.748463$ \\ 
$\mu_{s,1}$ & $0.0355597$ & $0.0382945$ \\ 
$\mu_{s,2}$ & $0.0356074$ & $0.0430675$ \\ 
$\mu_{cs,1}$ & $1.0014$ & $1.14097$ \\ 
$\mu_{cs,2}$ & $1.00196$ & $1.19691$ \\ \hline
\end{tabular}%
\end{center}
\caption{Coefficients for the rate conditions computed from (\protect\ref%
{eq:DF-bound1}) and (\protect\ref{eq:DF-bound2}).}
\label{Table:coefficients}
\end{table}

In a similar fashion one can compute the coefficients for other intervals,
and based on these compute the order of the rate conditions. In Table \ref%
{Table:rate-cond} we show a sequence of intervals spanning from $\varepsilon
=0$ to $\varepsilon =\frac{1}{100}$, together with the established order. 
\begin{table}[h]
\begin{center}
\begin{tabular}{llllllll}
\hline\hline
$\boldsymbol{\varepsilon}$ & order & \hspace{0.6cm} & $\boldsymbol{%
\varepsilon} $ & order & \hspace{0.6cm} & $\boldsymbol{\varepsilon}$ & order
\\ \hline
$[0, 0.0001]$ & 737 &  & $[0.0005, 0.001]$ & 73 &  & $[0.005, 0.006]$ & 11
\\ 
$[0.0001, 0.0002]$ & 368 &  & $[0.001, 0.002]$ & 36 &  & $[0.006, 0.007]$ & 9
\\ 
$[0.0002, 0.0003]$ & 245 &  & $[0.002, 0.003]$ & 24 &  & $[0.007, 0.008]$ & 8
\\ 
$[0.0003, 0.0004]$ & 184 &  & $[0.003, 0.004]$ & 17 &  & $[0.008, 0.009]$ & 7
\\ 
$[0.0004, 0.0005]$ & 147 &  & $[0.004, 0.005]$ & 14 &  & $[0.009, 0.01]$ & 6
\\ \hline
\end{tabular}%
\end{center}
\caption{Rate conditions order for various parameters.}
\label{Table:rate-cond}
\end{table}

To establish the covering condition we have numerically verified that $\pi
_{y}f_{\boldsymbol{\varepsilon }}\left( D_{\boldsymbol{\varepsilon }}\right)
\subset \mathrm{int}\pi _{y}D_{\boldsymbol{\varepsilon }}$ and that for $D_{%
\boldsymbol{\varepsilon }}^{-,\mathrm{left}}=\Lambda \times \{-R\}\times 
\overline{B}_{s}\left( R\right) $ and $D_{\boldsymbol{\varepsilon }}^{-,%
\mathrm{right}}=\Lambda \times \{R\}\times \overline{B}_{s}\left( R\right) $
holds%
\begin{equation*}
\pi _{x}f_{\boldsymbol{\varepsilon }}\left( D_{\boldsymbol{\varepsilon }}^{-,%
\mathrm{left}}\right) <-R\qquad \text{and\qquad }\pi _{x}f_{\boldsymbol{%
\varepsilon }}\left( D_{\boldsymbol{\varepsilon }}^{-,\mathrm{right}}\right)
>R.
\end{equation*}

Now we show how we verified the backward cone conditions. Since $\lambda \in 
\mathbb{R}$ mod $2\pi $, we can take $R_{\Lambda }=1$. If $p_{1}\in
J_{s}\left( p_{2},1/L\right) ,$ then%
\begin{equation*}
\left\Vert \pi _{\lambda }\left( p_{1}-p_{2}\right) \right\Vert \leq \frac{1%
}{L}\left\Vert \pi _{y}\left( p_{1}-p_{2}\right) \right\Vert \leq \frac{1}{L}%
2\varepsilon _{2}.
\end{equation*}%
Let $U=\left[ -\frac{2}{L}\varepsilon _{2},\frac{2}{L}\varepsilon _{2}\right]
\times \overline{B}_{u}\left( R\right) \times \overline{B}_{s}\left(
R\right) ,$ then $p_{1}-p_{2}\in U$ and%
\begin{equation*}
\left\Vert \pi _{\lambda }\left( f^{-1}(p_{1})-f^{-1}(p_{2})\right)
\right\Vert \leq \max \left[ \left\vert \pi _{\lambda }\left( Df\left(
D\right) \right) ^{-1}U\right\vert \right] .
\end{equation*}%
We verify numerically that $\max \left[ \left\vert \pi _{\lambda }\left(
Df\left( D\right) \right) ^{-1}U\right\vert \right] <R_{\Lambda }$. This
means that%
\begin{equation*}
f^{-1}(p_{1})\in \overline{B}_{c}(\pi _{\lambda }f^{-1}(p_{2}),R_{\Lambda
})\times \overline{B}_{u}(R)\times \overline{B}_{s}(R),
\end{equation*}%
and the backward cone condition for $z_{1}= f^{-1}(p_{1})$, $z_{2}=
f^{-1}(p_{2})$ follows from Corollary \ref{cor:stable-lip}.

\begin{remark}
The smoothness established in Table \ref{Table:rate-cond} is not optimal.
The example serves only to demonstrate that our method is applicable. We
choose a single change of coordinates and use global estimates on the
derivative of the map. With a more careful choice of changes to local
coordinates and by a local treatment of the estimates on the derivatives one
could obtain better results.
\end{remark}

All computations were performed using the CAPD\footnote{%
Computer Assisted Proofs in Dynamics: \href{http://capd.ii.uj.edu.pl/}{%
http://capd.ii.uj.edu.pl/}} package.

%TCIDATA{Version=5.00.0.2606}
%TCIDATA{LaTeXparent=0,0,MMFedit.tex}

\appendix

%\section*{Appendix}

\section{An auxiliary lemma}

\begin{lemma}
\label{lem:Taylor-dziwny-exp} Let $U\subset \mathbb{R}^{u}\times \mathbb{R}%
^{s}$ be a convex bounded neighborhood of zero and assume that $%
f:U\rightarrow \mathbb{R}^{s}$ is a $C^{m+1}$ map satisfying $f(0)=0$ and 
\begin{eqnarray}
\Vert f(U)\Vert _{C^{m+1}} &\leq &c,  \label{eq:allDerEstm} \\
\frac{\partial ^{l}f}{\partial x^{l}}(0) &=&0,\qquad \text{for }|l|\leq m.
\label{eq:derToXZero}
\end{eqnarray}

Then 
\begin{equation*}
f(x,y)=\frac{\partial f}{\partial y} y+g_{2}(\mathrm{x},\mathrm{y}),
\end{equation*}
where 
\begin{eqnarray*}
g_{2}(x,y) &\leq &C( \|y\|^{2}+\| x\| \|y\| +\|x\|^{m+1}) ,
\end{eqnarray*}%
with $C$ depending on $c$, the diameter of $U$ and $m$.
\end{lemma}

\begin{proof}
Let us consider the Taylor formula with the integral remainder of order $%
(m+1)$ (here the convexity is used). We group the second or higher order
terms in this expansion in three groups. The first group contains only the
terms independent of $x$. The second group contains both $x$ and $y$. The
sums in both groups can be can be bounded by $C_1 \|y\|^2$ and $C_2 \|x\|
\|y\|$, respectively, where constants $C_1$ and $C_2$ depend on $c$, the
diameter of $U$ and $m$. The last group contains a single term coming from
the reminder 
\begin{equation*}
\frac{1}{m!} \int_0^1 D^{m+1}f(t(x,y))(x^{[m+1]})dt,
\end{equation*}
which bounded by $\frac{c}{(m+1)!} \|x\|^{m+1}$.
\end{proof}

\section{Proof of Lemma \protect\ref{lem:func-in-unstb-jet}\label%
{app:func-in-unstb-jet}}

\begin{proof}
For sufficiently small $\delta $, if $\Vert \mathrm{x}-\mathrm{x}_{0}\Vert
\leq \delta $ then $M>\Vert D^{m+1}g(\mathrm{x})\Vert $, hence%
\begin{eqnarray*}
&&\Vert g(\mathrm{x})-g(\mathrm{x}_{0})-\mathcal{P}_{m}(\mathrm{x}-\mathrm{x}%
_{0})\Vert \\
&=&\left\Vert R_{m+1,p}(\mathrm{x}-\mathrm{x}_{0})\right\Vert \\
&=&\left\Vert \int_{0}^{1}\frac{\left( 1-t\right) ^{m}}{m!}D^{m+1}g(\mathrm{x%
}_{0}+th)\left( \left( \mathrm{x}-\mathrm{x}_{0}\right) ^{[m+1]}\right)
dt\right\Vert \\
&\leq &\int_{0}^{1}\frac{\left( 1-t\right) ^{m}}{m!}M\left\Vert \mathrm{x}-%
\mathrm{x}_{0}\right\Vert ^{m+1}dt \\
&=&\frac{M}{\left( m+1\right) !}\left\Vert \mathrm{x}-\mathrm{x}%
_{0}\right\Vert ^{m+1}
\end{eqnarray*}%
Therefore for $\left\Vert \mathrm{x}-\mathrm{x}_{0}\right\Vert \leq \delta $
we have 
\begin{equation*}
(\mathrm{x},g(\mathrm{x}))=(\mathrm{x}_{0},g(\mathrm{x}_{0}))+(\mathrm{x}-%
\mathrm{x}_{0},\mathcal{P}_{m}(\mathrm{x})+\mathrm{y}),
\end{equation*}%
where $\mathrm{y}=g(\mathrm{x})-g(\mathrm{x}_{0})-\mathcal{P}_{m}(\mathrm{x}%
) $ satisfies $\Vert \mathrm{y}\Vert \leq \frac{M}{(m+1)!} \Vert \mathrm{x}-%
\mathrm{x}_{0}\Vert ^{m+1}$. Hence (\ref{eq:func-in-unstb-jet}) is satisfied.
\end{proof}

\section{Proof of Lemma \protect\ref{lem:inv-func-in-unstb-jet}\label%
{app:inv-func-in-unstb-jet}}

The proof of Lemma \ref{lem:inv-func-in-unstb-jet} is based on the following
result.

\begin{lemma}
\label{lem:SymMultFormCoeff} Let $\|\cdot\|$ be an euclidean norm on $%
\mathbb{R}^n$. Let $B:\mathbb{R}^{n}\times \mathbb{R}^{n}\times \cdots
\times \mathbb{R}^{n}\rightarrow \mathbb{R}$ be $k$-linear symmetric form.
Assume that $M>0$ is such that for all $h\in \mathbb{R}^{n}$ holds 
\begin{equation*}
|B(h^{\left[ k\right] })|\leq M\| h\| ^{k}.
\end{equation*}
Let $\{e_{i}\}_{i=1,\dots ,n}\in \mathbb{R}^{n}$ be an orthonormal basis.

Then there exists $c=k^k$ such that for all $(i_{1},i_{2},\dots ,i_{k})\in
\{1,\dots ,n\}^{k}$ 
\begin{equation*}
|B(e_{i_{1}},e_{i_{2}},\dots ,e_{i_{k}})|\leq cM.
\end{equation*}
\end{lemma}

\begin{proof}
We now introduce some notations. For any set $Z$ by $\#Z$ we will denote its
number of elements. To deal with symmetric multiindices we define a set $%
S_{n,k}\subset \{1,\dots ,n\}^{k}$ by 
\begin{equation*}
S_{n,k}=\{i\in \{1,\dots ,n\}^{k}\ |\ i_{m}\leq i_{m+1},\quad m=1,\dots
,k-1\}.
\end{equation*}%
For any $i\in \{1,\dots ,n\}^{k}$ by $z=S(i)$ we denote a unique element in $%
S_{n,k}$, such that for each $j\in \{1,\dots ,n\}$ holds $\#\{m\ |\
i_{m}=j\}=\#\{m\ |\ z_{m}=j\}$. Hence $S(i)$ is an `ordered' $i$. For $i\in
S_{n,k}$ we define a multiplicity of $i$, denoted by $m(i)$, 
\begin{equation*}
m(i)=\#\{S^{-1}(i)\}.
\end{equation*}%
For $i \in S_{n,k}$ and $j\in \{1,\dots ,n\}$ we define a multiplicity of $j$
in $i$ by 
\begin{equation*}
m(j,i)=\#\{m\ |\ i_{m}=j\}.
\end{equation*}%
It is easy to see that 
\begin{equation*}
m(i)=\frac{k!}{\Pi _{j=1,\dots ,n}(m(j,i)!)}.
\end{equation*}%
For $i\in \{1,\dots ,n\}^{k}$ we write $x^{i}=x_{i_{1}}x_{i_{2}}\dots
x_{i_{k}}$.

Let us denote by $\mathcal{D}$ the diagonal of $B$, i.e. $\mathcal{D}:%
\mathbb{R}^{n}\rightarrow \mathbb{R}$, $\mathcal{D}(h)=B(h^{\left[ k\right]
})$. Let us consider the following polynomial of degree $k$ of $n$ variables 
$x_{1},\dots ,x_{n}$ 
\begin{eqnarray*}
P(x_{1},\dots ,x_{n}) &=&\mathcal{D}\left( \sum_{i=1}^{n}x_{i}e_{i}\right) \\
&=&\sum_{i\in \{1,\dots ,n\}^{k}}x_{i_{1}}x_{i_{2}}\dots
x_{i_{k}}B(e_{i_{1}},e_{i_{2}},\dots ,e_{i_{k}}) \\
&=&\sum_{l\in S_{n,k}}m(l)x_{l_{1}}x_{l_{2}}\dots
x_{l_{k}}B(e_{l_{1}},e_{l_{2}},\dots ,e_{l_{k}}).
\end{eqnarray*}

Now our task can is reduced to the following one: given bounds on $%
P(x_{1},\dots ,x_{n})$ can we produce bounds for its coefficients.

First of all we will develop a formula for each coefficient. To shorten some
expressions let us denote coefficients of $P$ by $p_{l}$, that is, 
\begin{equation}
P(x_{1},\dots ,x_{n})=\sum_{l\in S_{n,k}}p_{l}x^{l},\quad
p_{l}=m(l)B(e_{l_{1}},e_{l_{2}},\dots ,e_{l_{k}}).  \label{eq:polP}
\end{equation}
Each coefficient $p_{l}$ can be computed by finite differences as follows.

For any polynomial $W(x_{1},\dots ,x_{n})$ and $i=1,\dots ,n$ we define a
finite difference operator $\Delta _{i}$ as 
\begin{equation*}
(\Delta _{i}W)(x_{1},\dots ,x_{n}):=W(x_{1},\dots ,x_{i}+1/2,\dots
,x_{n})-W(x_{1},\dots ,x_{i}-1/2,\dots ,x_{n}).
\end{equation*}%
It is easy to see that $\Delta _{i}W$ is a polynomial, whose degree with
respect to variable $x_{i}$ decreases by $1$ (if it is nonzero). It is easy
to check that $\Delta _{i}\Delta _{j}=\Delta _{j}\Delta _{i}$. For $l\in
\{1,\dots ,n\}^{a}$ we set 
\begin{equation*}
\Delta ^{l}=\Delta _{l_{1}}\Delta _{l_{2}}\dots \Delta _{l_{s}}.
\end{equation*}%
We shall use the fact that for any polynomial $W(x_{1},\dots,x_{n})=%
\sum_{l}w_{l}x^{l}$, any $k$ and $l\in S_{n,k}$ we have 
\begin{equation*}
(\Delta ^{l}W)(0,\dots ,0)=\left( \Pi _{i=1,\dots ,n}(m(i,l)!)\right) w_{l}.
\end{equation*}

Observe that for polynomial $P$ given by (\ref{eq:polP}) and $l\in S_{n,k}$ $%
\Delta ^{l}P$ is a constant polynomial. Therefore from the above formula we
obtain 
\begin{equation*}
p_{l}=\left( \Pi _{i=1,\dots ,n}m(i,l)\right) ^{-1}\Delta ^{l}P.
\end{equation*}%
Now we are ready to estimate $p_{l}$. We set $(x_{1},\dots ,x_{n})=0$.
Observe that $\Delta ^{l}P$ will involve $2^{k}$ terms of the form $\pm
P(j_{1},\dots ,j_{n})$, where $j_{r}\in \{-k/2,\dots ,k/2\}$ and $%
\sum_{r=1}^n |j_{r}|\leq k/2$.

Hence 
\begin{equation*}
|p_l| \leq | \Delta ^{l}P| \leq 2^{k}\max_{\| x\| \leq k/2}|P(x_{1},\dots
,x_{n})|\leq M k^{k}
\end{equation*}
Therefore 
\begin{equation*}
|B(e_{l_{1}},\dots ,e_{l_{k}})|\leq \frac{Mk^{k}}{m(l)}\leq M k^{k}.
\end{equation*}
\end{proof}

We are now ready to prove Lemma \ref{lem:inv-func-in-unstb-jet}:

\begin{proof}
Using the Taylor formula 
\begin{equation}
\left\Vert R_{m+1,\mathrm{x}_{0}}(h)\right\Vert =\left\Vert g\left( \mathrm{x%
}_{0}+h\right) -g(\mathrm{x}_{0})-\mathcal{P}_{m}(h)\right\Vert \leq
M\left\Vert h\right\Vert ^{m+1}.  \label{eq:tmp-rem-1}
\end{equation}
Let $S^u$ denotes the sphere of radius $1$ in $\mathbb{R}^u$. Let $e\in
S^{u} $ and let $h=\eta e$ for $\eta \in [0,1]$. Then 
\begin{eqnarray}
R_{m+1,\mathrm{x}_{0}}(h) &=&\int_{0}^{1}\frac{\left( 1-t\right) ^{m}}{m!}%
D^{m+1}g(\mathrm{x}_{0}+th)\left( h^{[m+1]}\right) dt  \notag \\
&=&\eta^{m+1}\int_{0}^{1}\frac{\left( 1-t\right) ^{m}}{m!}D^{m+1}g(\mathrm{x}%
_{0}+t\eta e)\left( e^{[m+1]}\right) dt  \notag \\
&=&\eta^{m+1}\int_{0}^{1}\frac{\left( 1-t\right) ^{m}}{m!}D^{m+1}g(\mathrm{x}%
_{0})\left( e^{[m+1]}\right) dt+\eta ^{m+1}\varepsilon (\mathrm{x}%
_{0},e,\eta )  \notag \\
&=&\eta ^{m+1}\frac{D^{m+1}g(\mathrm{x}_{0})}{\left( m+1\right) !}\left(
e^{[m+1]}\right) +\eta ^{m+1}\varepsilon (\mathrm{x}_{0},e,\eta ),
\label{eq:tmp-rem-2}
\end{eqnarray}
where 
\begin{equation*}
\varepsilon (\mathrm{x}_{0},e,\eta )=\int_{0}^{1}\frac{( 1-t)^{m}}{m!}%
\left(D^{m+1}g(\mathrm{x}_{0}+t\eta e)\left( e^{[m+1]}\right) - D^{m+1}g(%
\mathrm{x}_{0})\left( e^{[m+1]}\right) \right)dt
\end{equation*}
Since $D^{m+1}g$ is continuous, $\varepsilon (\mathrm{x}_{0},e,\eta
)\rightarrow 0$ as $\eta \rightarrow 0$. Combining (\ref{eq:tmp-rem-1}) and (%
\ref{eq:tmp-rem-2}) we obtain 
\begin{equation*}
\eta ^{m+1}\left\Vert \frac{D^{m+1}g(\mathrm{x}_{0})}{\left( m+1\right) !}%
\left( e^{[m+1]}\right) \right\Vert -\eta ^{m+1}\left\Vert \varepsilon (%
\mathrm{x}_{0},e,\eta )\right\Vert \leq M\eta ^{m+1}.
\end{equation*}
Dividing by $\eta ^{m+1}$ and passing with $\eta $ to zero gives 
\begin{equation*}
\left\Vert \frac{D^{m+1}g(\mathrm{x}_{0})}{\left( m+1\right) !}\left(
e^{[m+1]}\right) \right\Vert \leq M.
\end{equation*}
This by Lemma \ref{lem:SymMultFormCoeff} gives 
\begin{eqnarray*}
\left\Vert \frac{\partial ^{m+1}g(\mathrm{x}_{0})}{\partial \mathrm{x}%
_{i_{1}}\ldots \partial \mathrm{x}_{i_{m+1}}}\right\Vert &\leq &\left\vert 
\frac{\partial ^{m+1}g_{1}(\mathrm{x}_{0})}{\partial \mathrm{x}%
_{i_{1}}\ldots \partial \mathrm{x}_{i_{m+1}}}\right\vert +\ldots +\left\vert 
\frac{\partial ^{m+1}g_{s}(\mathrm{x}_{0})}{\partial \mathrm{x}%
_{i_{1}}\ldots \partial \mathrm{x}_{i_{m+1}}}\right\vert \\
&\leq &s\left( m+1\right) !cM,
\end{eqnarray*}
which concludes our proof.
\end{proof}

\section{Proof of Theorem \protect\ref{thm:Lip-unstable-jet}\label%
{app:Lip-unstable-jet}}

\begin{proof}
If $(\mathrm{x},\mathrm{y})\in J_{u}(0,\mathcal{P}_{0},M)$ then $\left\Vert 
\mathrm{y}\right\Vert \leq M\left\Vert \mathrm{x}\right\Vert $. Since 
\begin{equation*}
f(\mathrm{x},\mathrm{y})=f(0)+\int_{0}^{1}Df(t(\mathrm{x},\mathrm{y}))dt(%
\mathrm{x},\mathrm{y})\in \left[ Df(U)\right] (\mathrm{x},\mathrm{y})
\end{equation*}%
by (\ref{eq:Lip-jet-cond-1}) we obtain 
\begin{eqnarray*}
\Vert \pi _{\mathrm{x}}f(\mathrm{x},\mathrm{y})\Vert &\geq &m\left( \left[ 
\frac{\partial f_{\mathrm{x}}}{\partial \mathrm{x}}(U)\right] \right) \Vert 
\mathrm{x}\Vert -\sup_{z\in U}\left\Vert \frac{\partial f_{\mathrm{x}}}{%
\partial \mathrm{y}}\left( z\right) \right\Vert \Vert \mathrm{y}\Vert \\
&\geq &\left( m\left( \left[ \frac{\partial f_{\mathrm{x}}}{\partial \mathrm{%
x}}(U)\right] \right) -M\sup_{z\in U}\left\Vert \frac{\partial f_{\mathrm{x}}%
}{\partial \mathrm{y}}\left( z\right) \right\Vert \right) \Vert \mathrm{x}%
\Vert \geq \xi \Vert \mathrm{x}\Vert .
\end{eqnarray*}

Using (\ref{eq:Lip-jet-cond-2}) in the last inequality, we have 
\begin{eqnarray*}
\Vert \pi _{\mathrm{y}}f(\mathrm{x},\mathrm{y})\Vert &=&\left\Vert
\int_{0}^{1}\frac{\partial f_{\mathrm{y}}}{\partial \mathrm{x}}(t(\mathrm{x},%
\mathrm{y}))\mathrm{x}+\frac{\partial f_{\mathrm{y}}}{\partial {\mathrm{y}}}%
(t(\mathrm{x},\mathrm{y})))\mathrm{y}dt\right\Vert \\
&\leq &\int_{0}^{1}\left\Vert \frac{\partial f_{\mathrm{y}}(t(\mathrm{x},%
\mathrm{y}))}{\partial \mathrm{x}}\right\Vert \Vert \mathrm{x}\Vert
+\left\Vert \frac{\partial f_{\mathrm{y}}(t(\mathrm{x},\mathrm{y}))}{%
\partial \mathrm{y}}\right\Vert \Vert \mathrm{y}\Vert dt \\
&\leq &M\int_{0}^{1}\left( \frac{1}{M}\left\Vert \frac{\partial f_{\mathrm{y}%
}(t(\mathrm{x},\mathrm{y}))}{\partial \mathrm{x}}\right\Vert \left\Vert 
\mathrm{x}\right\Vert +\left\Vert \frac{\partial f_{\mathrm{y}}(t(\mathrm{x},%
\mathrm{y}))}{\partial \mathrm{y}}\right\Vert \left\Vert \mathrm{x}%
\right\Vert \right) dt \\
&\leq &M\mu \left\Vert \mathrm{x}\right\Vert .
\end{eqnarray*}%
From the above estimates and (\ref{eq:Lip-jet-cond-3}), if $(\mathrm{x},%
\mathrm{y})\neq 0$ we obtain%
\begin{equation*}
\frac{\left\Vert \pi _{\mathrm{x}}f(\mathrm{x},\mathrm{y})\right\Vert }{%
\left\Vert \pi _{\mathrm{y}}f(\mathrm{x},\mathrm{y})\right\Vert }\geq \frac{%
\xi \left\Vert \mathrm{x}\right\Vert }{M\mu \left\Vert \mathrm{x}\right\Vert 
}>\frac{1}{M}.
\end{equation*}%
This implies 
\begin{equation*}
\left\Vert \pi _{\mathrm{y}}f(\mathrm{x},\mathrm{y})\right\Vert <M\left\Vert
\pi _{\mathrm{x}}f(\mathrm{x},\mathrm{y})\right\Vert ,
\end{equation*}%
hence, $f(\mathrm{x},\mathrm{y})\in \mathrm{int}J_{u}(0,\mathcal{R}%
_{0}=0,M), $ as required.
\end{proof}

\section{Proof of Theorem \protect\ref{thm:unstb-jet-propagation}\label%
{app:unstb-jet-propagation}}

We start by proving the theorem with an additional assumption that $\mathcal{%
P}_{m}=\mathcal{Q}_{m}=0$, and that $\frac{\partial ^{l}f_{\mathrm{y}}}{%
\partial \mathrm{x}^{l}}(0)=0$, for $|l|\leq m.$ We formulate this as a
lemma:

\begin{lemma}
\label{lem:unstbjet-in-good-frame} Let $U\subset \mathbb{R}^{u}\times 
\mathbb{R}^{s}$ be a convex bounded neighborhood of zero and assume that $%
f:U\rightarrow \mathbb{R}^{u}\times \mathbb{R}^{s}$ is a $C^{m+1}$ map
satisfying $f(0)=0$ and 
\begin{eqnarray}
\Vert f(U)\Vert _{C^{m+1}} &\leq &c,  \label{eq:unstjet-cond2} \\
\frac{\partial ^{l}f_{\mathrm{y}}}{\partial \mathrm{x}^{l}}(0) &=&0,\qquad 
\text{for }|l|\leq m.  \label{eq:unstjet-cond3}
\end{eqnarray}

If for $\xi >0$, and $\rho <1$%
\begin{equation}
\begin{array}{rr}
m\left( \frac{\partial f_{\mathrm{x}}}{\partial \mathrm{x}}(0)\right) \geq & 
\xi ,\medskip \\ 
\left\Vert \frac{\partial f_{\mathrm{x}}}{\partial \mathrm{y}}(0)\right\Vert
\leq & B,\medskip \\ 
\left\Vert \frac{\partial f_{\mathrm{y}}}{\partial \mathrm{y}}(0)\right\Vert
\leq & \mu ,%
\end{array}
\label{eq:unstjet-cond1}
\end{equation}%
and%
\begin{equation}
\frac{\mu }{\xi ^{m+1}}<\rho ,  \label{eq:unstbratecond}
\end{equation}%
then there exists a constant $M^{\ast }=M^{\ast }\left( c,B,1/\xi \right) $,
such that for any $M>M^{\ast }$ there exists\textbf{\ }$\delta =\delta
(M,c,B,1/\xi )$ such that 
\begin{equation*}
f(J_{u}(0,\mathcal{P}_{m}=0,M,\delta )\cap U)\subset J_{u}(0,\mathcal{P}%
_{m}=0,M).
\end{equation*}%
Moreover, if for some $K>0$ holds $c,B,\frac{1}{\xi }\in \left[ 0,K\right] $%
, then $M^{\ast }$ depends only on $K$ and $\rho $.
\end{lemma}

\begin{proof}
Let us introduce the following notations 
\begin{equation*}
D_{11}=\frac{\partial f_{\mathrm{x}}}{\partial \mathrm{x}}(0),\quad D_{12}=%
\frac{\partial f_{\mathrm{x}}}{\partial \mathrm{y}}(0),\quad D_{22}=\frac{%
\partial f_{\mathrm{y}}}{\partial \mathrm{y}}(0),
\end{equation*}%
then since $\frac{\partial ^{l}f_{\mathrm{y}}}{\partial \mathrm{x}^{l}}(0)=0$
for $\left\vert l\right\vert \leq m$ 
\begin{equation*}
f(\mathrm{x},\mathrm{y})=(D_{11}\mathrm{x}+D_{12}\mathrm{y}+g_{1}(\mathrm{x},%
\mathrm{y}),D_{22}\mathrm{y}+g_{2}(\mathrm{x},\mathrm{y})),
\end{equation*}%
where by the Taylor formula and Lemma~\ref{lem:Taylor-dziwny-exp} 
\begin{eqnarray*}
g_{1}(\mathrm{x},\mathrm{y}) &\leq &C(\Vert \mathrm{x}\Vert ^{2}+\Vert 
\mathrm{y}\Vert ^{2}), \\
g_{2}(\mathrm{x},\mathrm{y}) &\leq &C\left( \Vert \mathrm{y}\Vert
^{2}+\left\Vert \mathrm{x}\right\Vert \left\Vert \mathrm{y}\right\Vert
+\left\Vert \mathrm{x}\right\Vert ^{m+1}\right) ,
\end{eqnarray*}%
for $(x,y)\in U$, $\Vert (x,y)\Vert \leq 1$ with $C$ depending on $c$.

Let $(\mathrm{x},\mathrm{y})\in J_{u}(0,0,M)\cap U$. Then $\left\Vert 
\mathrm{y}\right\Vert \leq M\left\Vert \mathrm{x}\right\Vert ^{m+1}$. Let $(%
\mathrm{x}_{1},\mathrm{y}_{1})=f(\mathrm{x},\mathrm{y})$. We have 
\begin{eqnarray*}
\Vert \mathrm{x}_{1}\Vert &\geq &m(D_{11})\Vert \mathrm{x}\Vert -\Vert
D_{12}\Vert \cdot \Vert \mathrm{y}\Vert -C(\Vert \mathrm{x}\Vert ^{2}+\Vert 
\mathrm{y}\Vert ^{2}) \\
&\geq &\xi \Vert \mathrm{x}\Vert -BM\left\Vert \mathrm{x}\right\Vert
^{m+1}-C\Vert \mathrm{x}\Vert ^{2}(1+M^{2}\left\Vert \mathrm{x}\right\Vert
^{2m})
\end{eqnarray*}%
It is apparent that there exists $\delta =\delta (M,c,B,1/\xi )>0$, such
that if $\Vert x\Vert \leq \delta $, then $x_{1}$ is positive. Observe that $%
\delta $ is decreasing with respect to all of its arguments.

We now compute 
\begin{eqnarray*}
\Vert \mathrm{y}_{1}\Vert &\leq &\Vert D_{22}\Vert \Vert \mathrm{y}\Vert
+C\left( \Vert \mathrm{y}\Vert ^{2}+\Vert \mathrm{x}\Vert \Vert \mathrm{y}%
\Vert +\Vert \mathrm{x}\Vert ^{m+1}\right) \\
&\leq &M\left\Vert \mathrm{x}\right\Vert ^{m+1}\left( \mu +C\left(
M\left\Vert \mathrm{x}\right\Vert ^{m+1}+\Vert \mathrm{x}\Vert +\frac{1}{M}%
\right) \right) .
\end{eqnarray*}%
By further decreasing $\delta $ if necessary we obtain for $\Vert \mathrm{x}%
\Vert \leq \delta $ the following inequalities 
\begin{equation*}
\left\Vert \mathrm{x}\right\Vert ^{m}\leq M^{-2},\qquad \left\Vert \mathrm{x}%
\right\Vert ^{m+1}\leq M^{-2},\qquad \left\Vert \mathrm{x}\right\Vert
^{2m}\leq M^{-2}
\end{equation*}%
hence by (\ref{eq:unstbratecond}), for sufficiently large $M$ 
\begin{eqnarray*}
\frac{\Vert \mathrm{y}_{1}\Vert }{\Vert \mathrm{x}_{1}\Vert ^{m+1}} &\leq &%
\frac{M\left\Vert \mathrm{x}\right\Vert ^{m+1}\left( \mu +C\left(
M\left\Vert \mathrm{x}\right\Vert ^{m+1}+\Vert \mathrm{x}\Vert +\frac{1}{M}%
\right) \right) }{\left( \xi \Vert \mathrm{x}\Vert -BM\left\Vert \mathrm{x}%
\right\Vert ^{m+1}-C\Vert \mathrm{x}\Vert ^{2}(1+M^{2}\left\Vert \mathrm{x}%
\right\Vert ^{2m})\right) ^{m+1}} \\
&=&M\frac{\mu +C\left( M\left\Vert \mathrm{x}\right\Vert ^{m+1}+\Vert 
\mathrm{x}\Vert +\frac{1}{M}\right) }{\xi ^{m+1}\left( 1-\frac{1}{\xi }%
BM\left\Vert \mathrm{x}\right\Vert ^{m}-\frac{1}{\xi }C\Vert \mathrm{x}\Vert
(1+M^{2}\left\Vert \mathrm{x}\right\Vert ^{2m})\right) ^{m+1}} \\
&\leq &M\frac{\mu +\frac{1}{M}C\left( 2+\frac{1}{M}\right) }{\xi
^{m+1}\left( 1-\frac{1}{M}\frac{1}{\xi }\left( B+2C\right) \right) ^{m+1}} \\
&\leq &M\frac{\rho +\frac{1}{M\xi ^{m+1}}C\left( 2+\frac{1}{M}\right) }{%
\left( 1-\frac{1}{M}\frac{1}{\xi }\left( B+2C\right) \right) ^{m+1}} \\
&\leq &M.
\end{eqnarray*}%
The choice of the size of $M$ depends on $C,B,\rho $ and $\frac{1}{\xi }$.
Since $\frac{\Vert \mathrm{y}_{1}\Vert }{\Vert \mathrm{x}_{1}\Vert ^{m+1}}%
\leq M$ we have shown that $(\mathrm{x}_{1},\mathrm{y}_{1})\in J_{u}(0,0,M),$
as required.
\end{proof}

We are now ready to prove Theorem \ref{thm:unstb-jet-propagation}:

\begin{proof}
We would like to change the coordinates around zero, so that the map $f$ in
these coordinates will satisfy the assumptions of Lemma~\ref%
{lem:unstbjet-in-good-frame}.

Let $(\mathrm{x}_{0},\mathrm{y}_{0})$ be the new coordinates in the
neighborhood of zero given by $(\mathrm{x},\mathrm{y})=\Phi
_{z_{0}\rightarrow z}(\mathrm{x}_{0},\mathrm{y}_{0})$ 
\begin{eqnarray*}
\mathrm{x} &=&\mathrm{x}_{0}, \\
\mathrm{y} &=&\mathrm{y}_{0}+\mathcal{P}_{m}(\mathrm{x}_{0}),
\end{eqnarray*}%
We denote the inverse transformation as $\Phi_{z\rightarrow z_{0}}=\Phi
_{z_{0}\rightarrow z}^{-1}$.

Analogously, let us also consider coordinates $(\mathrm{x}_{1},\mathrm{y}%
_{1})$ given by $\Phi _{z_{1}\rightarrow z}(\mathrm{x}_{1},\mathrm{y}_{1})=(%
\mathrm{x},\mathrm{y})$ 
\begin{eqnarray*}
\mathrm{x} &=&\mathrm{x}_{1}, \\
\mathrm{y} &=&\mathrm{y}_{1}+\mathcal{R}_{m}(\mathrm{x}_{1}).
\end{eqnarray*}%
and denote $\Phi _{z\rightarrow z_{1}}=\Phi _{z_{1}\rightarrow z}^{-1}$.

Observe that both inverse transformations $\Phi _{z_{0}\rightarrow z}^{-1}$
and $\Phi _{z_{1}\rightarrow z}^{-1}$ are polynomial: 
\begin{eqnarray*}
\Phi _{z_{0}\rightarrow z}^{-1}(\mathrm{x},\mathrm{y}) &=&(\mathrm{x},%
\mathrm{y}-\mathcal{P}_{m}(\mathrm{x})) \\
\Phi _{z_{1}\rightarrow z}^{-1}(\mathrm{x},\mathrm{y}) &=&(\mathrm{x},%
\mathrm{y}-\mathcal{R}_{m}(\mathrm{x}))
\end{eqnarray*}%
and satisfy he same bound on the coefficients.

Now let $\tilde{f}(\mathrm{x}_{0},\mathrm{y}_{0})=\Phi _{z\rightarrow
z_{1}}(f(\Phi _{z_{0}\rightarrow z}(\mathrm{x}_{0},\mathrm{y}_{0})))$, i.e.
we express $f$ in new coordinates.

Observe that in coordinates $(\mathrm{x}_{0},\mathrm{y}_{0})$ the set $%
J_{u}(0,\mathcal{P}_{m},M)$ is just $J_{u}(0,0,M)$, i.e. $\Phi_{z\rightarrow
z_{0}}(J_{u}(0,\mathcal{P}_{m},M))=J_{u}(0,0,M)$. Analogously, $%
\Phi_{z\rightarrow z_{1}}(J_{u}(0,\mathcal{R}_{m},M))=J_{u}(0,0,M)$.

Now we compute the derivative of $\tilde{f}$. We have 
\begin{eqnarray*}
D\tilde{f}(0) &=&\left[ 
\begin{array}{cc}
I & 0 \\ 
-D\mathcal{R}_m(0) & I%
\end{array}%
\right] \cdot \left[ 
\begin{array}{cc}
Df_{11} & Df_{12} \\ 
Df_{21} & Df_{22}%
\end{array}%
\right] \cdot \left[ 
\begin{array}{cc}
I & 0 \\ 
D\mathcal{P}_m(0) & I%
\end{array}%
\right] \\
&=&\left[ 
\begin{array}{cc}
D\tilde{f}_{11} & D\tilde{f}_{12} \\ 
D\tilde{f}_{21} & D\tilde{f}_{22}%
\end{array}%
\right] ,
\end{eqnarray*}%
hence%
\begin{eqnarray*}
D\tilde{f}_{11} &=&Df_{11}+Df_{12}D\mathcal{P}(0), \\
D\tilde{f}_{12} &=&Df_{12}, \\
D\tilde{f}_{22} &=&Df_{22}-D\mathcal{R}(0)Df_{12}.
\end{eqnarray*}%
By (\ref{eq:full-unstb-cond-1}--\ref{eq:full-unstb-rate-cond}) we see that
assumptions (\ref{eq:unstjet-cond1}--\ref{eq:unstbratecond}) of Lemma \ref%
{lem:unstbjet-in-good-frame} are satisfied.

We now show that assumption (\ref{eq:unstjet-cond2}) from Lemma \ref%
{lem:unstbjet-in-good-frame} is satisfied with a common constant $c$ for 
\underline{all} polynomials $\mathcal{P}_{m}$ and $\mathcal{R}_{m}$
satisfying our assumptions. The fact that $\Vert \tilde{f}(D)\Vert
_{C^{m+1}} $ is bounded follows from the fact that $\Vert f(D)\Vert
_{C^{m+1}}\leq C$ and since $\Phi _{z_{0}\rightarrow z},$ $\Phi
_{z\rightarrow z_{1}}$ are polynomial changes of coordinates. We assumed
that the coefficients of $\mathcal{P}$, $\mathcal{R}$ are bounded by $C,$
hence $\Vert \tilde{f}(D)\Vert _{C^{m+1}}$ can be bounded by a constant
dependent only on $C$, $m$ and the size of the set $D$.

What remains is to verify that condition (\ref{eq:unstjet-cond3}) holds for $%
\tilde{f}$. From the definition of $\Phi _{z\rightarrow z_{1}}$ we see that%
\begin{equation}
\pi _{\mathrm{y}}\circ \Phi _{z\rightarrow z_{1}}\circ \left( \mathrm{x},%
\mathcal{R}_{m}(\mathrm{x})\right) =0.  \label{eq:tmp-jet-graph-cond}
\end{equation}%
By (\ref{eq:full-unstb-graph-cond}), for any $\mathrm{x}$%
\begin{equation*}
T_{f\circ \left( \mathrm{id},\mathcal{P}_{m}\right) ,m,0}(\mathrm{x})=(%
\mathrm{x}^{\prime },\mathcal{R}_{m}(\mathrm{x}^{\prime }))
\end{equation*}%
for some $\mathrm{x}^{\prime }\in \mathbb{R}^{u}$, hence from (\ref%
{eq:tmp-jet-graph-cond}) it follows that%
\begin{equation*}
\pi _{\mathrm{y}}\circ \Phi _{z\rightarrow z_{1}}\circ T_{f\circ \left( 
\mathrm{id},\mathcal{P}_{m}\right) ,m,0}(\mathrm{x})=0.
\end{equation*}%
The Taylor expansion of $\pi _{\mathrm{y}}\circ \Phi _{z\rightarrow
z_{1}}\circ T_{f\circ \left( \mathrm{id},\mathcal{P}_{m}\right) ,m,0}$ up to
order $m$ is equal to the Taylor expansion of $\pi _{\mathrm{y}}\circ \Phi
_{z\rightarrow z_{1}}\circ f\circ \left( \mathrm{id},\mathcal{P}_{m}\right) $
up to order $m$. This means that%
\begin{equation*}
T_{\pi _{\mathrm{y}}\circ \Phi _{z\rightarrow z_{1}}\circ f\circ \left( 
\mathrm{id},\mathcal{P}_{m}\right) ,m,0}(\mathrm{x})=0.
\end{equation*}%
Since $\left( \mathrm{id},\mathcal{P}_{m}\right) (\mathrm{x})=\Phi
_{z_{0}\rightarrow z}(\mathrm{x},0),$ above implies that 
\begin{equation*}
T_{\pi _{\mathrm{y}}\circ \tilde{f}\circ (\mathrm{id},0),m,0}(\mathrm{x})=0,
\end{equation*}%
hence%
\begin{equation*}
\frac{\partial ^{l}\tilde{f}_{\mathrm{y}}}{\partial \mathrm{x}^{l}}(0)=0.
\end{equation*}%
We have thus shown (\ref{eq:unstjet-cond3}), which concludes the proof.
\end{proof}

\section{Proof of Theorem \protect\ref{thm:Lip-stable-jet}\label%
{app:Lip-stable-jet}}

\begin{proof}
We first observe that%
\begin{equation*}
J_{s}^{c}(0,\mathcal{Q}_{0}=0,M)=\left\{ \left\Vert \mathrm{x}\right\Vert
>M\left\Vert \mathrm{y}\right\Vert \right\} =\mathrm{int}J_{u}\left( 0,0,%
\frac{1}{M}\right) .
\end{equation*}%
Conditions (\ref{eq:Lip-jet-cond-4}--\ref{eq:Lip-jet-cond-6}) imply that
assumptions of Theorem \ref{thm:Lip-unstable-jet} are satisfied (for $1/M$
in place of $M$). This means that 
\begin{eqnarray*}
f\left( \overline{J_{s}^{c}(0,\mathcal{Q}_{0}=0,M)}\cap U\right) &=&f\left(
J_{u}\left( 0,0,\frac{1}{M}\right) \cap U\right) \\
&\subset &\mathrm{int}J_{u}\left( 0,0,\frac{1}{M}\right) \cup \{0\} \\
&=&J_{s}^{c}(0,\mathcal{R}_{0}=0,M)\cup \{0\},
\end{eqnarray*}%
as required.
\end{proof}

\section{Proof of Theorem \protect\ref{thm:stb-jet-propagation}\label%
{app:stb-jet-propagation}}

The proof goes along the same lines as the proof of Theorem \ref%
{thm:unstb-jet-propagation}. There are some differences though in the needed
estimates.

Similarly to the proof of Theorem \ref{thm:unstb-jet-propagation}, we start
by proving the theorem with an additional assumption that $\mathcal{Q}_{m}=%
\mathcal{R}_{m}=0,$ and that $\frac{\partial ^{l}f_{\mathrm{x}}}{\partial 
\mathrm{y}^{l}}(0)=0,$ for $|l|\leq m.$ We formulate this as a lemma:

\begin{lemma}
\label{lem:stbjet-in-good-frame} Let $U\subset \mathbb{R}^{u}\times \mathbb{R%
}^{s}$ be a convex bounded neighborhood of zero and assume that $%
f:U\rightarrow \mathbb{R}^{u}\times \mathbb{R}^{s}$ is a $C^{m+1}$ map
satisfying $f(0)=0$ and 
\begin{eqnarray}
\Vert f(U)\Vert _{C^{m+1}} &\leq &c,  \label{eq:stjet-cond2} \\
\frac{\partial ^{l}f_{\mathrm{x}}}{\partial \mathrm{y}^{l}}(0) &=&0,\qquad 
\text{for }|l|\leq m.  \label{eq:stjet-cond3}
\end{eqnarray}

If for $\xi >0$, and $\rho <1$%
\begin{equation}
\begin{array}{rr}
m\left( \frac{\partial f_{\mathrm{x}}}{\partial \mathrm{x}}(0)\right) \geq & 
\xi ,\medskip \\ 
\left\Vert \frac{\partial f_{\mathrm{y}}}{\partial \mathrm{x}}(0)\right\Vert
\leq & B,\medskip \\ 
\left\Vert \frac{\partial f_{\mathrm{y}}}{\partial \mathrm{y}}(0)\right\Vert
\leq & \mu ,%
\end{array}
\label{eq:stjet-cond1}
\end{equation}%
and%
\begin{equation}
\frac{\mu ^{m+1}}{\xi }<\rho ,  \label{eq:stbratecond}
\end{equation}%
then there exists a constant $M^{\ast }=M^{\ast }\left( c,B,1/\xi ,\rho
\right) $, such that for any $M>M^{\ast }$ there exists $\delta =\delta
(M,c,B,1/\xi )$ such that 
\begin{equation*}
f(J_{s}^{c}(0,\mathcal{P}_{m}=0,M,\delta )\cap U)\subset J_{s}^{c}(0,%
\mathcal{P}_{m}=0,M).
\end{equation*}%
Moreover, if for some $K>0$ holds $c,B,\frac{1}{\xi }\in \lbrack 0,K]$, then 
$M^{\ast }$ depends only on $K$ and $\rho $.
\end{lemma}

\begin{proof}
Let us introduce the following notations 
\begin{equation*}
D_{11}=\frac{\partial f_{\mathrm{x}}}{\partial \mathrm{x}}(0),\quad D_{21}=%
\frac{\partial f_{\mathrm{y}}}{\partial \mathrm{x}}(0),\quad D_{22}=\frac{%
\partial f_{\mathrm{y}}}{\partial \mathrm{y}}(0),
\end{equation*}%
then since $\frac{\partial ^{l}f_{\mathrm{x}}}{\partial \mathrm{y}^{l}}(0)=0$
for $\left\vert l\right\vert \leq m$ 
\begin{equation*}
f(\mathrm{x},\mathrm{y})=(D_{11}\mathrm{x}+g_{1}(\mathrm{x},\mathrm{y}%
),D_{21}\mathrm{x}+D_{22}\mathrm{y}+g_{2}(\mathrm{x},\mathrm{y})),
\end{equation*}%
where by the Taylor formula and Lemma~\ref{lem:Taylor-dziwny-exp} 
\begin{eqnarray*}
\Vert g_{1}(\mathrm{x},\mathrm{y})\Vert &\leq &C\left( \Vert \mathrm{x}\Vert
^{2}+\left\Vert \mathrm{x}\right\Vert \left\Vert \mathrm{y}\right\Vert
+\left\Vert \mathrm{y}\right\Vert ^{m+1}\right) , \\
\Vert g_{2}(\mathrm{x},\mathrm{y})\Vert &\leq &C(\Vert \mathrm{x}\Vert
^{2}+\Vert \mathrm{y}\Vert ^{2}),
\end{eqnarray*}%
for $(x,y)\in U\cap B(0,1)$ with $C$ depending on $c$.

Let $(\mathrm{x},\mathrm{y})\in J_{s}^{c}(0,0,M)\cap U$. Then $\left\Vert 
\mathrm{x}\right\Vert >M\left\Vert \mathrm{y}\right\Vert ^{m+1}$. Let $(%
\mathrm{x}_{1},\mathrm{y}_{1})=f(\mathrm{x},\mathrm{y})$. From (\ref%
{eq:stjet-cond1}) we have 
\begin{eqnarray}
\Vert \mathrm{x}_{1}\Vert &\geq &m(D_{11})\Vert \mathrm{x}\Vert -C\left(
\Vert \mathrm{x}\Vert ^{2}+\left\Vert \mathrm{x}\right\Vert \left\Vert 
\mathrm{y}\right\Vert +\left\Vert \mathrm{y}\right\Vert^{m+1}\right)  \notag
\\
&\geq &\Vert \mathrm{x}\Vert \left( \xi -C\left( \Vert \mathrm{x}\Vert
+\left( \frac{1}{M}\left\Vert \mathrm{x}\right\Vert \right) ^{\frac{1}{m+1}}+%
\frac{1}{M}\right) \right) .  \label{eq:tmp-x1-stb}
\end{eqnarray}%
It is apparent that taking $M$ sufficiently large and sufficiently small $%
\left\Vert \mathrm{x}\right\Vert $, the lower bound for $\Vert \mathrm{x}%
_{1}\Vert $ is positive.

We now compute 
\begin{equation}
\Vert \mathrm{y}_{1}\Vert \leq \left\Vert D_{21}\right\Vert \left\Vert 
\mathrm{x}\right\Vert +\Vert D_{22}\Vert \Vert \mathrm{y}\Vert +C(\Vert 
\mathrm{x}\Vert ^{2}+\Vert \mathrm{y}\Vert ^{2}).  \label{eq:tmp-y1-stb}
\end{equation}%
Taking $\left\Vert \mathrm{x}\right\Vert \leq M^{-m}$ we see that $\left( 
\frac{1}{M}\left\Vert \mathrm{x}\right\Vert \right) ^{\frac{1}{m+1}}\leq
M^{-1},$ hence for $(\mathrm{x},\mathrm{y})\in J_{s}^{c}(0,0,M),$ 
\begin{eqnarray*}
\left\Vert \mathrm{y}\right\Vert &<&\left( \frac{1}{M}\left\Vert \mathrm{x}%
\right\Vert \right) ^{\frac{1}{m+1}}, \\
\left\Vert \mathrm{y}\right\Vert ^{2} &<&M^{-1}\left( \frac{1}{M}\left\Vert 
\mathrm{x}\right\Vert \right) ^{\frac{1}{m+1}}.
\end{eqnarray*}

If $M>1$ and $\Vert \mathrm{x}\Vert \leq M^{-4m}$ then 
\begin{eqnarray*}
\left\Vert \mathrm{x}\right\Vert &=&\left\Vert \mathrm{x}\right\Vert ^{\frac{%
1}{m+1}}\left\Vert \mathrm{x}\right\Vert ^{\frac{m}{m+1}} \\
&\leq &\left\Vert \mathrm{x}\right\Vert ^{\frac{1}{m+1}}\left(
M^{-4m}\right) ^{\frac{m}{m+1}} \\
&=&\frac{1}{M}\left( \frac{1}{M}\left\Vert \mathrm{x}\right\Vert \right) ^{%
\frac{1}{m+1}}M^{\frac{m+2-4m^{2}}{m+1}} \\
&<&\frac{1}{M}\left( \frac{1}{M}\left\Vert \mathrm{x}\right\Vert \right) ^{%
\frac{1}{m+1}},
\end{eqnarray*}%
and%
\begin{equation}
\left\Vert \mathrm{x}\right\Vert ^{2}\leq M^{-1}\left( \frac{1}{M}\left\Vert 
\mathrm{x}\right\Vert \right) ^{\frac{1}{m+1}}.  \notag
\end{equation}

This by (\ref{eq:tmp-y1-stb}) and (\ref{eq:stjet-cond1}) means that for $M>1$
and $\left\Vert \mathrm{x}\right\Vert \leq M^{-4m}\boldsymbol{\ }$ 
\begin{eqnarray*}
\Vert \mathrm{y}_{1}\Vert &\leq &\left( \frac{1}{M}\left\Vert \mathrm{x}%
\right\Vert \right) ^{\frac{1}{m+1}}\left[ \left\Vert D_{21}\right\Vert
M^{-1}+\Vert D_{22}\Vert +2CM^{-1})\right] \\
&\leq &\left( \frac{1}{M}\left\Vert \mathrm{x}\right\Vert \right) ^{\frac{1}{%
m+1}}\left[ \mu +\frac{1}{M}\left( B+2C\right) )\right] ,
\end{eqnarray*}%
which combined with (\ref{eq:tmp-x1-stb}) gives 
\begin{eqnarray*}
\frac{\Vert \mathrm{y}_{1}\Vert ^{m+1}}{\Vert \mathrm{x}_{1}\Vert } &\leq &%
\frac{\frac{1}{M}\left\Vert \mathrm{x}\right\Vert \left[ \mu +\frac{1}{M}%
\left( B+2C\right) )\right] ^{m+1}}{\Vert \mathrm{x}\Vert \left( \xi
-C\left( \Vert \mathrm{x}\Vert +\left( \frac{1}{M}\left\Vert \mathrm{x}%
\right\Vert \right) ^{\frac{1}{m+1}}+\frac{1}{M}\right) \right) } \\
&\leq &\frac{\frac{1}{M}\left[ \mu +\frac{1}{M}\left( B+2C\right) )\right]
^{m+1}}{\xi -3C/M} \\
&=&\frac{1}{M}\frac{\left( \frac{\mu }{\xi ^{1/m+1}}+\frac{B+2C}{M\xi
^{1/m+1}}\right) ^{m+1}}{(1-\frac{3C}{\xi M})} \\
&<&\frac{1}{M}\frac{\left( \rho ^{1/m+1}+\frac{B+2C}{M\xi ^{1/m+1}}\right)
^{m+1}}{(1-\frac{3C}{\xi M})}
\end{eqnarray*}%
Since $\rho <1$, taking sufficiently large $M$ (the choice of $M$ depends on 
$C,B$, $1/\xi $ and $\rho $), we see that $\frac{\Vert \mathrm{y}_{1}\Vert
^{m+1}}{\Vert \mathrm{x}_{1}\Vert }<\frac{1}{M}$, hence $(\mathrm{x}_{1},%
\mathrm{y}_{1})\in J_{s}^{c}\left( 0,0,M\right) $, as required.
\end{proof}

We are now ready to give the proof of Theorem \ref{thm:stb-jet-propagation}.

\begin{proof}
The claim of Theorem \ref{thm:stb-jet-propagation} follows from Lemma \ref%
{lem:stbjet-in-good-frame}, by considering $f$ in suitable local
coordinates. The proof follows from a mirror argument to the proof of
Theorem \ref{thm:unstb-jet-propagation}, with the only difference that we
need to swap the roles of the coordinates $\mathrm{x}$ and $\mathrm{y}$.
\end{proof}

\section{Proof of Lemma \protect\ref{lem:Ju-expansion}\label%
{app:Ju-expansion}}

\begin{proof}
By Remark \ref{rem:cone-in-chart}, $z_1$ and $z_2$ are contained in the same
chart.

Since $z_{1}\in J_{u}\left( z_{2},1/L\right) $%
\begin{equation*}
\left\Vert \pi _{\left( \lambda ,y\right) }\left( z_{1}-z_{2}\right)
\right\Vert \leq \frac{1}{L}\left\Vert \pi _{x}\left( z_{1}-z_{2}\right)
\right\Vert .
\end{equation*}%
We have 
\begin{equation*}
f(z_{1})-f(z_{2})=\int_{0}^{1}Df(z_{2}+t(z_{1}-z_{2}))dt(z_{1}-z_{2}),
\end{equation*}%
hence 
\begin{eqnarray*}
&&\left\Vert \pi _{x}\left( f(z_{1})-f(z_{2})\right) \right\Vert \\
&\geq &m\left( \frac{\partial f_{x}}{\partial x}(P(z_{1}))\right) \left\Vert
\pi _{x}\left( z_{1}-z_{2}\right) \right\Vert -\sup_{z\in D}\left\Vert \frac{%
\partial f_{x}}{\partial (\lambda ,y)}z)\right\Vert \left\Vert \pi _{\left(
\lambda ,y\right) }\left( z_{1}-z_{2}\right) \right\Vert \\
&\geq &\inf_{z\in D}m\left( \frac{\partial f_{x}}{\partial x}(P(z))\right)
\left\Vert \pi _{x}\left( z_{1}-z_{2}\right) \right\Vert -\frac{1}{L}%
\sup_{z\in D}\left\Vert \frac{\partial f_{x}}{\partial \left( \lambda
,y\right) }\left( z\right) \right\Vert \left\Vert \pi _{x}\left(
z_{1}-z_{2}\right) \right\Vert \\
&\geq &\xi _{u,1,P}\Vert \pi _{x}(z_{1}-z_{2})\Vert ,
\end{eqnarray*}
\end{proof}

\section{Proof of Lemma \protect\ref{lem:Js-contraction}\label%
{app:Js-contraction}}

\begin{proof}
By Remark \ref{rem:cone-in-chart}, $z_{1}$ and $z_{2}$ are contained in the
same chart.

Since $z_{1}\in J_{s}\left( z_{2},1/L\right) ,$%
\begin{equation*}
\left\Vert \pi _{\theta }\left( z_{1}-z_{2}\right) \right\Vert
<1/L\left\Vert \pi _{y}\left( z_{1}-z_{2}\right) \right\Vert .
\end{equation*}%
We have 
\begin{equation*}
f(z_{1})-f(z_{2})=\int_{0}^{1}Df(z_{2}+t(z_{1}-z_{2})))dt(z_{1}-z_{2}).
\end{equation*}

This implies that 
\begin{eqnarray*}
&&\left\Vert \pi _{y}f(z_{1})-\pi _{y}f(z_{2})\right\Vert \\
&=&\left\Vert \int_{0}^{1}\frac{\partial f_{y}}{\partial y}%
(z_{2}+t(z_{1}-z_{2}))\pi _{y}(z_{1}-z_{2})+\frac{\partial f_{y}}{\partial
\theta }(z_{2}+t(z_{1}-z_{2}))\pi _{\theta }(z_{1}-z_{2})dt\right\Vert \\
&\leq &\int_{0}^{1}\left\Vert \frac{\partial f_{y}}{\partial y}%
(z_{2}+t(z_{1}-z_{2}))\right\Vert \left\Vert \pi
_{y}(z_{1}-z_{2})\right\Vert +\left\Vert \frac{\partial f_{y}}{\partial
\theta }(z_{2}+t(z_{1}-z_{2}))\right\Vert \left\Vert \pi _{\theta
}(z_{1}-z_{2})\right\Vert dt \\
&\leq &\sup_{z\in D}\left( \left\Vert \frac{\partial f_{y}}{\partial y}%
(z)\right\Vert +\frac{1}{L}\left\Vert \frac{\partial f_{y}}{\partial \theta }%
(z)\right\Vert \right) \left\Vert \pi _{y}(z_{1}-z_{2})\right\Vert \\
&\leq &\mu _{s,1}\left\Vert \pi _{y}(z_{1}-z_{2})\right\Vert ,
\end{eqnarray*}%
as required.
\end{proof}

\section{Proof of Lemma \protect\ref{lem:cu-cone-xi-cu}\label%
{app:cu-cone-xi-cu}}

\begin{proof}
Since $z_{1}\in J_{cu}\left( z_{2},L\right) $, 
\begin{equation*}
\left\Vert \pi _{y}\left( z_{1}-z_{2}\right) \right\Vert \leq L\left\Vert
\pi _{\left( \lambda ,x\right) }\left( z_{1}-z_{2}\right) \right\Vert .
\end{equation*}%
We have 
\begin{equation*}
f(z_{1})-f(z_{2})=\int_{0}^{1}Df(z_{2}+t(z_{1}-z_{2}))dt(z_{1}-z_{2}).
\end{equation*}%
This gives%
\begin{eqnarray*}
&&\left\Vert \pi _{\left( \lambda ,x\right) }\left( f(z_{1})-f(z_{2})\right)
\right\Vert \\
&\geq &m\left( \frac{\partial f_{(\lambda ,x)}}{\partial \left( \lambda
,x\right) }(P(z_{1}))\right) \Vert \pi _{\lambda ,x}(z_{1}-z_{2})\Vert
-\sup_{z\in D}\left\Vert \frac{\partial f_{(\lambda ,x)}}{\partial y}%
(z)\right\Vert \Vert \pi _{y}(z_{1}-z_{2})\Vert \\
&\geq &\left( \inf_{z\in D}m\left( \frac{\partial f_{(\lambda ,x)}}{\partial
(\lambda ,x)}(P(z))\right) -L\sup_{z\in D}\left\Vert \frac{\partial
f_{(\lambda ,x)}}{\partial y}(z)\right\Vert \right) \left\Vert \pi _{\left(
\lambda ,x\right) }\left( z_{1}-z_{2}\right) \right\Vert \\
&\geq &\xi _{cu,1,P}\left\Vert \pi _{\left( \lambda ,x\right) }\left(
z_{1}-z_{2}\right) \right\Vert ,
\end{eqnarray*}
\end{proof}

\section{Proof of Lemma \protect\ref{lem:cs-cone-mu-cs}\label%
{app:cs-cone-mu-cs}}

\begin{proof}
Since $z_{1}\in J_{cs}\left( z_{2},L\right) $, 
\begin{equation*}
\left\Vert \pi _{x}\left( z_{1}-z_{2}\right) \right\Vert \leq L\left\Vert
\pi _{(\lambda ,y)}\left( z_{1}-z_{2}\right) \right\Vert .
\end{equation*}%
We have 
\begin{equation*}
f(z_{1})-f(z_{2})=\int_{0}^{1}Df(z_{2}+t(z_{1}-z_{2})))dt(z_{1}-z_{2}).
\end{equation*}%
hence 
\begin{eqnarray*}
&&\Vert \pi _{(\lambda ,y)}\left( f(z_{1})-f(z_{2})\right) \Vert \\
&=&\left\Vert \int_{0}^{1}Df_{(\lambda
,y)}(z_{2}+t(z_{1}-z_{2})))(z_{1}-z_{2})dt\right\Vert \\
&\leq &\int_{0}^{1}\left\Vert \frac{\partial f_{(\lambda ,y)}}{\partial
\left( \lambda ,y\right) }(z_{2}+t(z_{1}-z_{2}))\right\Vert \left\Vert \pi
_{\left( \lambda ,y\right) }\left( z_{1}-z_{2}\right) \right\Vert \\
&&+\left\Vert \frac{\partial f_{\left( \lambda ,y\right) }}{\partial x}%
(z_{2}+t(z_{1}-z_{2}))\right\Vert \left\Vert \pi _{x}\left(
z_{1}-z_{2}\right) \right\Vert dt \\
&\leq &\sup_{z\in D}\left( \left\Vert \frac{\partial f_{\left( \lambda
,y\right) }}{\partial \left( \lambda ,y\right) }(z)\right\Vert +L\left\Vert 
\frac{\partial f_{\left( \lambda ,y\right) }}{\partial x}(z)\right\Vert
\right) \left\Vert \pi _{\left( \lambda ,y\right) }\left( z_{1}-z_{2}\right)
\right\Vert \\
&\leq &\mu _{cs,1}\left\Vert \pi _{\left( \lambda ,y\right) }\left(
z_{1}-z_{2}\right) \right\Vert ,
\end{eqnarray*}%
as required.
\end{proof}

\section{Proof of Lemma \protect\ref{lem:unstable-disc}\label%
{app:unstable-disc}}

\begin{proof}
Let $z=b(0)$ and $\lambda ^{\ast }\in \Lambda $ be the point from Definition~%
\ref{def:covering} for $z$. Note that since $b$ is a horizontal disc, $b(%
\overline{B}_{u}(R))\subset J_{u}\left( z,1/L\right) $. This also means that%
\begin{equation}
\left\Vert \pi _{\left( \lambda ,y\right) }\left( b\left( x_{1}\right)
-b\left( x_{2}\right) \right) \right\Vert \leq \frac{1}{L}\left\Vert \pi
_{x}\left( b\left( x_{1}\right) -b\left( x_{2}\right) \right) \right\Vert =%
\frac{1}{L}\left\Vert x_{1}-x_{2}\right\Vert .  \label{eq:Lip-b-unst}
\end{equation}
From Definition~\ref{def:covering} follows that 
\begin{equation*}
f(b(\overline{B}_{u}(R)))\subset B_{c}(\lambda ^{\ast },R_{\Lambda })\times 
\mathbb{R}^{u}\times \mathbb{R}^{s},
\end{equation*}%
hence 
\begin{equation*}
f(b(\overline{B}_{u}(R)))\cap D\subset D_{\lambda ^{\ast }}.
\end{equation*}%
Observe that by Remark \ref{rem:cone-in-chart}, $h_{t}$ maps the disk $b$ in
a set contained in a single chart.

We start by showing that for any $\hat{x}\in \overline{B}_{u}(R)$ there
exists $x=x( \hat{x})$ such that 
\begin{equation}
\pi_x f(b(x))=\hat{x}.  \label{eq:x-hat-target-new}
\end{equation}
and then disk $b^\ast$ will be defined by $b^\ast(\hat{x})=f(b(x(\hat{x})))$.

Let us fix $\hat{x}\in \overline{B}_{u}(R)$ and consider a function 
\begin{equation*}
F: \overline{B}_u(R) \to \mathbb{R}^u
\end{equation*}%
defined as%
\begin{equation*}
F(x)= \pi_x f(b(x))- \hat{x}.
\end{equation*}%
Our objective is to show that there exists a unique $x$ such that $F(x)=0$.

Let $h_{\alpha }$ be the homotopy from Definition \ref{def:covering}. Let us
define a homotopy 
\begin{equation*}
H:[0,1]\times \overline{B}_{u}(R)\rightarrow \mathbb{R}^{u}
\end{equation*}%
\begin{equation*}
H_{\alpha }(x)=\pi _{x}h_{\alpha }(b(x))-\hat{x}.
\end{equation*}%
Note that $H_{0}=F$. We will start by showing that 
\begin{equation}
\forall \alpha \in \lbrack 0,1]\quad \forall x\in \partial B_{u}(R)\qquad
H_{\alpha }(q)\neq 0.  \label{eq:pixF-nonzero-bnd}
\end{equation}%
To prove (\ref{eq:pixF-nonzero-bnd}) let us take $x\in \partial B_{u}(R)$.
Since $b(x)\in J_{u}\left( z,1/L\right) \cap D_{\pi _{\lambda }(z)}^{-}$, by
condition (\ref{eq:homotopy-exit}) from Definition \ref{def:covering} $%
h_{\alpha }(b(x))\notin D_{\lambda ^{\ast }}$, which means that $h_{\alpha
}(b(x))\neq \hat{x}$, implying $H_{\alpha }\left( q\right) \neq 0$.

Let $U\subset \mathbb{R}^{n}$ be a set and $q\in \mathbb{R}^{n}$ be a point.
We use the notation $\deg \left( F,U,q\right) $ for the Brouwer degree of $F 
$ with respect to the set $D$ at $q$ . From condition (\ref%
{eq:pixF-nonzero-bnd}) by the homotopy property of the Brouwer degree (see 
\cite{Ll}), we obtain 
\begin{equation}
\deg (F,B_{u}(R),0)=\deg (H_{\alpha },B_{u}(R),0)=\deg (H_{1},B_{u}(R),0).
\label{eq:homotop-F-new}
\end{equation}

Our next step is to show that $\deg (H_{1},B_{u}(R),0)\neq 0$. Since $%
h_{1}(x)=Ax$ we see that 
\begin{equation*}
H_{1}(x)=(Ax,0)-\hat{x}.
\end{equation*}

By point 4. from Definition~\ref{def:covering} it follows that $\det(A) \neq
0$ and $A^{-1}\hat{x} \in B_u(R)$. Therefore equation $H_1(q)=0$ has a
unique solution in $B_u(R)$ and by the degree property for affine maps 
\begin{equation*}
\deg ( H_{1},B_u(R),0) =\mathrm{sgn} \det A = \pm 1
\end{equation*}

By (\ref{eq:homotop-F-new}), this gives 
\begin{equation*}
\deg \left( F,B_u(R),0\right) =\deg \left( H_{1},B_u(R),0\right) \neq 0.
\end{equation*}
This means that there exists an $x\in B_u(R)$ such that $F(x)=0.$ This
finishes the proof of (\ref{eq:x-hat-target-new}).

We now define the candidate for $b^{\ast }(\hat{x})$ as 
\begin{equation}
b^{\ast }\left( \hat{x}\right) =f\circ b(x(\hat{x})).
\label{eq:b-str-tmp-new}
\end{equation}%
By construction, $\pi _{x}b^{\ast }(\hat{x})=\hat{x}$. We need to show that $%
b^{\ast }(\hat{x})$ is well defined (meaning that the choice of $x(\hat{x})$
is unique), and that it is a horizontal disc.

Let $x_{1}\neq x_{2}$. By Lemma \ref{lem:Ju-expansion} we have%
\begin{equation}
\Vert \pi _{x}(f\circ b(x_{1})-f\circ b(x_{2}))\Vert \geq \xi _{u,1,P}\Vert
\pi _{x}(b(x_{1})-b(x_{2}))\Vert =\xi _{u,1,P}\Vert x_{1}-x_{2}\Vert \neq 0.
\label{eq:fbx1x2-below-estm}
\end{equation}%
Hence, $b^{\ast }(\hat{x})$ is well defined.

Observe that (\ref{eq:fbx1x2-below-estm}) can be rewritten as 
\begin{equation*}
\Vert \hat{x}_{1}-\hat{x}_{2}\Vert =\Vert \pi _{x}(f\circ b(x(\hat{x_{1}}%
))-f\circ b(x(\hat{x_{2}})))\Vert \geq \xi _{u,1,P}\Vert x(\hat{x}_{1})-x(%
\hat{x}_{2})\Vert .
\end{equation*}%
Therefore $x(\hat{x})$ is Lipschitz, hence $b^{\ast }$ is continuous.

We will now show that for any $\hat{x}\in \overline{B}_{u}(R)$ we have $%
b^{\ast }\left( \overline{B}_{u}\right) \subset J_{u}\left( b^{\ast }(\hat{x}%
),1/L\right) $. By Corollary \ref{cor:unstable-lip}%
\begin{equation*}
f\left( J_{u}\left( b(x),1/L\right) \cap D\right) \subset J_{u}\left( f\circ
b(x),1/L\right) .
\end{equation*}%
Since for any $x$ we have $b(\overline{B}_{u})\subset J_{u}\left(
b(x),1/L\right) ,$ we obtain%
\begin{equation*}
f\circ b(\overline{B}_{u})\subset J_{u}\left( f\circ b(x),1/L\right) ,
\end{equation*}%
which by the definition of $b^{\ast }$ from (\ref{eq:b-str-tmp-new}) implies%
\begin{equation*}
b^{\ast }(\overline{B}_{u})\subset J_{u}\left( b^{\ast }(\hat{x}),1/L\right)
,
\end{equation*}%
as required.

We now need to show that if $f,b$ are $C^{k}$, for $k\geq 1,$ then so is $%
b^{\ast }$. Let us introduce the notation%
\begin{eqnarray*}
g &:&\overline{B}_{u}(R)\rightarrow \mathbb{R}^{u}, \\
g(x) &=&\pi _{x}f\circ b(x).
\end{eqnarray*}%
We can rewrite the definition of $b^{\ast }$ using $g$ as%
\begin{equation*}
b^{\ast }\left( x^{\ast }\right) =f\circ b\circ g^{-1}(x^{\ast }).
\end{equation*}%
To show that $b^{\ast }$ is $C^{k},$ it is sufficient for $g^{-1}$ to be $%
C^{k}$. From (\ref{eq:Lip-b-unst}) we see that $\pi _{\left( \lambda
,y\right) }b$ is Lipschitz with the constant $1/L$, hence 
\begin{eqnarray*}
m(Dg(x)) &=&m(D\pi _{x}f\circ b(x)) \\
&=&m\left( \frac{\partial f_{x}}{\partial x}(b(x))+\frac{\partial f_{x}}{%
\partial \left( \lambda ,y\right) }(b(x))\frac{\partial \pi _{\left( \lambda
,y\right) }b}{\partial \left( \lambda ,y\right) }(x)\right) \\
&\geq &m\left( \frac{\partial f_{x}}{\partial x}(b(x))\right) -\left\Vert 
\frac{\partial f_{x}}{\partial \left( \lambda ,y\right) }(b(x))\frac{%
\partial \pi _{\left( \lambda ,y\right) }b}{\partial \left( \lambda
,y\right) }(x)\right\Vert \\
&\geq &m\left( \frac{\partial f_{x}}{\partial x}(b(x))\right) -\frac{1}{L}%
\left\Vert \frac{\partial f_{x}}{\partial \left( \lambda ,y\right) }%
(b(x))\right\Vert \geq \xi _{u,1}>0.
\end{eqnarray*}%
and by the inverse function theorem $g^{-1}$ is $C^{k}$; as required.
\end{proof}

\section{Proof of Lemma \protect\ref{lem:center-unstable-disc}\label%
{app:center-unstable-disc}}

\begin{proof}
First, we will prove that 
\begin{equation}
\left( f\circ b(\Lambda \times \overline{B}_{u}(R))\right) \cap D\neq
\emptyset .  \label{eq:chd-image-on-D}
\end{equation}%
For any $\lambda \in \Lambda $ let us consider $b^{\lambda }:\overline{B}%
_{u}(R)\rightarrow D$ given by $b^{\lambda }(x)=b(\lambda ,x)$. We will
argue that $b^{\lambda }$ is a horizontal disc. We first observe that 
\begin{equation*}
\pi _{x}b^{\lambda }(x)=\pi _{x}b(\lambda ,x)=\pi _{x}\pi _{\left( \lambda
,x\right) }b(\lambda ,x)=\pi _{x}(\lambda ,x)=x.
\end{equation*}%
We need to show that 
\begin{equation}
b^{\lambda }(\overline{B}_{u}(R))\subset J_{u}\left( b^{\lambda }\left(
x\right) ,1/L\right) .  \label{eq:tmp-hor-disc-in-cone}
\end{equation}%
Since $b$ is a center-horizontal disc, by definition, for any $%
x_{1},x_{2}\in \overline{B}_{u}(R),$%
\begin{equation*}
b(\lambda ,x_{2})\in J_{u}\left( b\left( \lambda ,x_{2}\right) ,L\right) ,
\end{equation*}%
hence%
\begin{equation*}
\left\Vert \pi _{y}b(\lambda ,x_{1})-\pi _{y}b(\lambda ,x_{2})\right\Vert
<L\left\Vert \pi _{\left( \lambda ,x\right) }b(\lambda ,x_{1})-\pi _{\left(
\lambda ,x\right) }b(\lambda ,x_{2})\right\Vert .
\end{equation*}%
Since $\pi _{\left( \lambda ,x\right) }b(\lambda ,x_{i})=(\lambda ,x_{i}),$
this gives (remember that $L<1$) 
\begin{eqnarray*}
\left\Vert \pi _{(\lambda ,y)}b^{\lambda }(x_{1})-\pi _{(\lambda
,y)}b^{\lambda }(x_{2})\right\Vert &=&\left\Vert \pi _{(\lambda
,y)}b(\lambda ,x_{1})-\pi _{(\lambda ,y)}b(\lambda ,x_{2})\right\Vert \\
&=&\left\Vert \left( \lambda ,\pi _{y}b(\lambda ,x_{1})\right) -\left(
\lambda ,\pi _{y}b(\lambda ,x_{2})\right) \right\Vert \\
&=&\left\Vert \pi _{y}b(\lambda ,x_{1})-\pi _{y}b(\lambda ,x_{2})\right\Vert
\\
&<&L\left\Vert \pi _{\left( \lambda ,x\right) }b(\lambda ,x_{1})-\pi
_{\left( \lambda ,x\right) }b(\lambda ,x_{2})\right\Vert \\
&=&L\left\Vert x_{1}-x_{2}\right\Vert \\
&=&L\left\Vert \pi _{x}b^{\lambda }(x_{1})-\pi _{x}b^{\lambda
}(x_{2})\right\Vert \\
&\leq &1/L\left\Vert \pi _{x}b^{\lambda }(x_{1})-\pi _{x}b^{\lambda
}(x_{2})\right\Vert ,
\end{eqnarray*}%
which implies (\ref{eq:tmp-hor-disc-in-cone}). We have thus shown that $%
b^{\lambda }$ is a horizontal disc.

From Lemma \ref{lem:unstable-disc} it follows that $f\circ b^{\lambda }(%
\overline{B}_{u}(R))\cap D$ is a horizontal disk in $D$, in particular this
implies (\ref{eq:chd-image-on-D}).

In the remainder of the proof we will use notation $\theta =(\lambda ,x)$.

We will now show that $\pi _{\theta }f\circ b$ is an open map, in fact it is
continuous and locally injective.

Let us fix $\theta _{1}$ and let us take $U$, an convex open neighborhood
contained in a single chart and such that $f(b(U))$ is contained in a single
chart. From Lemma~\ref{lem:cu-cone-xi-cu} it follows that 
\begin{equation*}
\Vert \pi _{\theta }f\circ b(\theta _{1})-\pi _{\theta }f\circ b(\theta
_{2})\Vert \geq \xi _{cu,1,P}\Vert \theta _{1}-\theta _{2}\Vert .
\end{equation*}%
Therefore $\pi _{\theta }f\circ b:U\rightarrow \Lambda \times \mathbb{R}^{u}$
is continuous and injective, hence by the Brouwer open map theorem we know $%
\pi _{\theta }f\circ b(U)$ is an open set. This means that $\pi _{\theta
}f\circ b$ is an open map, and therefore $\pi _{\theta }f\circ b(\Lambda
\times B_{u}(R))$ is an open set.

From the covering relation (Definition \ref{def:covering}) we know that the
points $b(\theta )$ for $\theta \in \Lambda \times \partial \overline{B}%
_{u}(R)$ are mapped by $f$ out of the set $D$ 
\begin{equation*}
\pi _{\theta }f\circ b(\Lambda \times B_{u}(R))\cap (\Lambda \times 
\overline{B}_{u}(R))=\pi _{\theta }f\circ b(\Lambda \times \overline{B}%
_{u}(R))\cap (\Lambda \times \overline{B}_{u}(R)).
\end{equation*}%
Therefore the set $\pi _{\theta }f\circ b(\Lambda \times B_{u}(R))\cap
(\Lambda \times \overline{B}_{u}(R))$ is both open and closed in $\Lambda
\times \overline{B}_{u}(R)$ and since it is also nonempty and $\Lambda
\times \overline{B}_{u}(R)$ is connected, we infer that 
\begin{equation}
\pi _{\theta }f\circ b(\Lambda \times B_{u}(R))\cap (\Lambda \times 
\overline{B}_{u}(R))=\Lambda \times \overline{B}_{u}(R).
\label{eq:ch-graph-transform-def-tmp}
\end{equation}

We need to show that the map $\pi _{\theta }f\circ b$ is an injection on $%
(\pi _{\theta }f\circ b)^{-1}(\Lambda \times \overline{B}_{u}(R))$. This is
a direct consequence of the backward cone condition (see Definition \ref%
{def:back-cc}). To show this, assume that there exists $\theta _{1}\neq
\theta _{2}$ in $\Lambda \times \overline{B}_{u}(R)$ such that 
\begin{equation*}
\pi _{\theta }f(b(\theta _{1}))=\pi _{\theta }f(b(\theta _{2})).
\end{equation*}%
Then%
\begin{equation*}
f(b(\theta _{1}))\in J_{s}\left( f(b(\theta _{2})),1/L\right) ,
\end{equation*}%
therefore the backward cone condition implies that%
\begin{equation*}
b(\theta _{1})\in J_{s}\left( b(\theta _{2}),1/L\right) ,
\end{equation*}%
which contradicts condition (\ref{eq:ch-cone-cond}) required of
center-horizontal disks.

We have shown (\ref{eq:ch-graph-transform-def-tmp}), which means that for
any $\theta ^{\ast }\in \Lambda \times \overline{B}_{u}(R)$ there exists an $%
\theta $ such that 
\begin{equation*}
\pi _{\theta }f\circ b(\theta )=\theta ^{\ast }.
\end{equation*}%
Such $\theta $ is unique due to the fact that $\pi _{\theta }f\circ b$ is
injective. We can therefore define%
\begin{equation*}
b^{\ast }(\theta ^{\ast })=f\circ b(\theta ).
\end{equation*}%
From the construction of $b^{\ast }$ it follows that $\pi _{\theta }b^{\ast
}(\theta ^{\ast })=\theta ^{\ast }$. Condition (\ref{eq:ch-cone-cond}) is a
consequence of backward cone conditions, and follows from a mirror argument
to the one used for the proof of injectivity of $\pi _{\theta }f\circ b$,
which was done in the preceding paragraph.

What is left is to show that if $f,b$ are $C^{k}$, for $k\geq 1,$ then so is 
$b^{\ast }$. This follows from mirror arguments to the proof of $C^{k}$
smoothness in Lemma \ref{lem:unstable-disc}.
\end{proof}

%\section*{References}

\end{document}